    \RequirePackage{fix-cm}
    \documentclass[smallextended]{svjour3}       
    \smartqed  

\usepackage{lipsum}
\usepackage[utf8]{inputenc}
\usepackage{amssymb}
\usepackage{amsmath}
\usepackage{graphicx} 
\usepackage{xspace}
\usepackage{multirow}
\usepackage{bbm}
\usepackage{algorithm}
\usepackage{algpseudocode}
\usepackage{graphicx}
\usepackage{enumitem}
\usepackage{multirow}
\usepackage{adjustbox}
\usepackage{bm}
\usepackage{xcolor}
\usepackage{soul}
\usepackage{hyperref}

\usepackage{geometry}
\newcommand{\igor}[1]{\textcolor{black}{#1}}
\newcommand{\mn}[1]{\textcolor{black}{#1}}
\DeclareMathOperator{\diag}{diag}

\begin{document}

\title{Residual viscosity stabilized RBF-FD methods for solving nonlinear conservation laws}


\author{Igor Tominec         \and
        Murtazo Nazarov 
}

\institute{I. Tominec \at
              Uppsala University, Department of Information Technology, Division of Scientific Computing \\
              \email{igor.tominec@it.uu.se}           
           \and
           M. Nazarov \at
           Uppsala University, Department of Information Technology, Division of Scientific Computing \\
           \email{murtazo.nazarov@it.uu.se}
}
\date{}

\maketitle

\begin{abstract}
    In this paper, we solve nonlinear conservation laws using the radial basis function generated finite difference (RBF-FD) method. Nonlinear conservation laws have solutions that entail
    strong discontinuities and shocks, which give rise to numerical instabilities when the solution
    is approximated by a numerical method. We introduce a residual-based artificial viscosity (RV)
    stabilization framework adjusted to the RBF-FD method, where the residual of the conservation law adaptively
    locates discontinuities and shocks. The RV stabilization framework is applied to the collocation RBF-FD
    method and the oversampled RBF-FD method. Computational tests confirm that the stabilized methods 
    are reliable and accurate in solving scalar conservation laws and conservation law systems such as compressible Euler equations.


\end{abstract}


\section{Introduction}
\begin{figure}[h!]
    \centering
\includegraphics[width=0.44\linewidth]{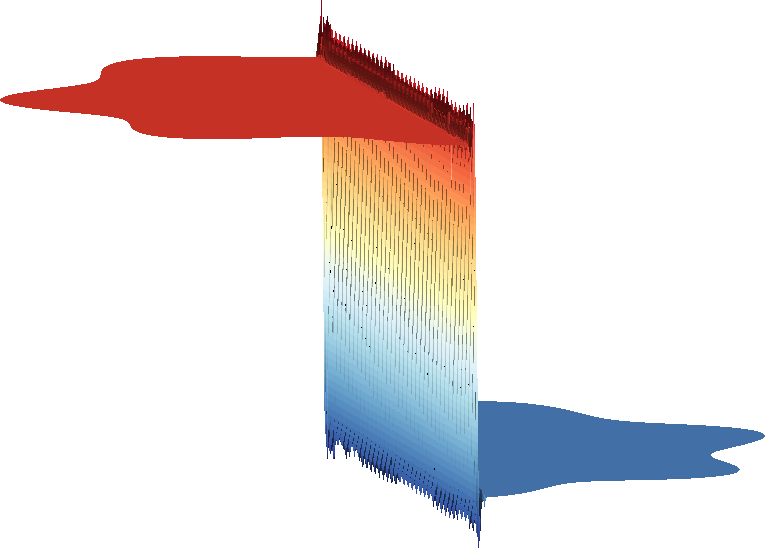}
\includegraphics[width=0.44\linewidth]{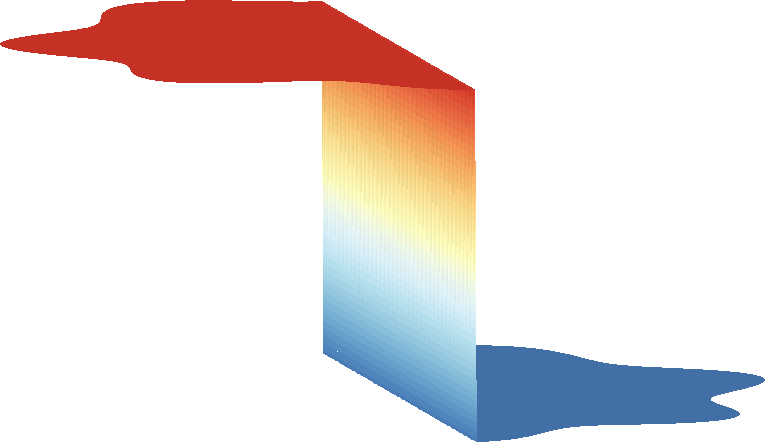}
\caption{A discontinuous function drawn over a star shaped domain. Left: a numerical approximation exhibiting oscillations (Gibbs phenomenon) obtained using the RBF-FD method. Right: the exact function without oscillations.}
\label{fig:intro:stepFunction}
\end{figure}

We are interested in solving the following system of nonlinear conservation laws\igor{:}
\begin{equation}
    \label{eq:model:conservationlaw}
\partial_t \bm U(y,t) = - \nabla \cdot \bm F(\bm U(y,t)) \quad (y,t) \in \Omega \times \mathbb R^+\igor{,}
\end{equation}
with appropriate initial and boundary conditions, 
where $\Omega \subset \mathbb{R}^2$ is an open and bounded domain, $\bm U \in \mathbb{R}^{k}$ is the solution and $\bm F ( \bm U) \in \mathcal C^1(\mathbb{R}^k; \mathbb{R}^{2})$ is a given smooth flux function. 
Since the flux is a nonlinear function of $\bm U$, the solutions of \eqref{eq:model:conservationlaw} lead to discontinuities and shocks in finite time. Numerical approximation of these discontinuities leads to non-physical oscillations, 
the so-called Gibbs phenomenon, see, e.g., Figure \ref{fig:intro:stepFunction}, which make the numerical scheme unstable. Therefore, adding additional stabilization terms to the scheme is vital to suppressing non-physical numerical oscillations.

State-of-the-art numerical methods for solving \eqref{eq:model:conservationlaw} use finite difference, finite volume, and discontinuous Galerkin approximations. One well-known approach to control and suppress the Gibbs phenomenon is the artificial viscosity method proposed by von Neumann and Richtmyer \cite{NeumannRichtmyer} in the 1950s in the context of finite difference methods. This method has been successfully used to solve many complex conservation laws, but the artificial viscosity term in \cite{NeumannRichtmyer} was not consistent with the PDE in its design. 

This consistency problem was later resolved by invoking the least-squares argument in the finite element community, see, e.g., \cite{Hughes_et_al_2010} and references in it. The so-called Galerkin Least Squares (GLS) method is consistent but adds artificial viscosity only in the streamline direction, which is insufficient to prove convergence. The authors of \cite{SDandRV} proved the convergence of the GLS approximations by supplementing the method with the residual-based artificial viscosity (GLS+RV) method.
The GLS+RV method is challenging to implement and requires a coupled space-time discretization. 
Recently in \cite{Nazarov13} we proved that the RV term is the main convergence mechanism, 
so we proposed to simplify the method by altogether abandoning the GLS terms. 
The RV method has also been used to stabilize systems of conservation laws using other discretization techniques in \cite{NazarovHoffman13,stiernstrom2021,MarrasNazarov15,Lu_spectral_rv}. Stabilization of finite element approximations \igor{only utilizing} artificial viscosity conditions is not new: the well-known entropy viscosity method of Guermond et al. \cite{Guermond11} constructs the artificial viscosity coefficient using the entropy residual of the PDE. 

\igor{The aforementioned stabilization techniques have not been extended to the radial basis function (RBF) methods yet}. A standard way of stabilizing hyperbolic time-dependent problems within the RBF community is to augment the collocation scheme using a hyperviscosity term that damps high-frequency oscillations \cite{Shankar_hypervi1,Shankar_hypervi2,FornbergLehto,Flyer12}. \igor{The hyperviscosity term acts on the higher derivatives of the numerical solution to establish a stable eigenvalue spectrum of the discretized advection operator. Hyperviscosity is effective in stabilizing the high-frequency mode numerical oscillations, however, as shown later in this article, for example, in Section \ref{sec:experiments:kpp}, it is ineffective to resolve the Gibbs phenomenon resulting from the presence of shocks and discontinuities. Thus, the RBF approximation of nonlinear conservation laws requires additional shock-capturing terms.} 


In this paper we base our stabilization approach on the artificial viscosity method, where we add a parabolic term to the RBF-FD discretization \igor{in addition to hyperviscosity}. The objective of the RV stabilization is to detect the position where the shock starts to form, utilizing the magnitude of the numerical residual in \eqref{eq:model:conservationlaw}. In the regions of $\Omega$ where the numerical solution starts to produce an overshoot over a shock, the residual is expected to be large, whereas, in the regions where the solution is free of oscillations, the residual is expected to be small. 

Without loss of generality, we center our discussions on the RV stabilization, when the stabilization is augmenting an oversampled RBF-FD discretization. More traditional collocation RBF-FD schemes can be seen as a special case of oversampling. \mn{The idea can easily be applied to different RBF approximation techniques.}

The paper is organized as follows. In Section \ref{sec:discretization} we formulate an oversampled RBF-FD method for solving a scalar conservation law, and the Euler system of equations that models compressed gas dynamics. In Section \ref{sec:residualviscosity} 
we derive a residual viscosity stabilization in the context of the RBF-FD method, for scalar conservation laws and Euler's system of PDEs. 
In Section \ref{sec:experiments:linearadvection} we make a comparison between the oversampled and the collocation RBF-FD methods when solving 
a linear conservation law (linear advection problem) and numerically verify the RV stabilization consistency. 
In Section \ref{sec:experiments:burgers} and Section \ref{sec:experiments:kpp} we present numerical results for solving 
Burger's and the Kurganov-Petrova-Popov equations respectively. In Section  \ref{sec:experiments:Euler} we solve several benchmark problems for 
the Euler equations.
Final remarks are given in Section \ref{sec:finaldiscussion}.

\section{An oversampled/collocated RBF-FD discretization}
\label{sec:discretization}
In this section we describe the oversampled RBF-FD discretization and the collcated RBF-FD discretization 
that we use to solve 
\eqref{eq:model:conservationlaw}. 
In particular, we describe our choice of point sets spread over $\Omega$, 
the formation of evaluation and differentiation matrices using the RBF-FD method, 
and the discretization of a time-dependent PDE where an explicit method is used to advance the solution in time. 

\subsection{Point sets}
The domain $\Omega$ is discretized using two point sets:
\begin{itemize}
\item the interpolation point set $X = \{x_i\}_{i=1}^N$ for generating the cardinal functions (left plot in Figure \ref{fig:discretization:rbffd_stencil}),
\item the evaluation point set $Y = \{y_j\}_{j=1}^M$ for sampling the PDE \eqref{eq:model:conservationlaw} (right plot in Figure \ref{fig:discretization:rbffd_stencil}).
\end{itemize}
In the present work, we generate the $X$ point set such that the mean distance between the points is set to a given value $h$. Depending on the application, 
the algorithms which we use are the 2D node generator from \cite{FBF15_nodes}, Gmsh \cite{Gmsh} or DistMesh \cite{Distmesh}.
The $Y$ point set is generated by placing $q$ points in each Voronoi region centered around every $x_i \in X$, $i=1,..,N$, 
so that the relation between the number of $X$ and $Y$ points is $M \approx \lceil q N \rceil$, 
where we define $q$ as the oversampling parameter. \igor{When $q=1$ and $Y=X$, then the point sets are ready for a collocation type of discretization. 
When $q>1$, then the discretization is oversampled.}

\subsection{Evaluation and differentiation matrices}
\label{sec:discretization:matrices}
\begin{figure}
    \centering
\includegraphics[width=0.4\linewidth]{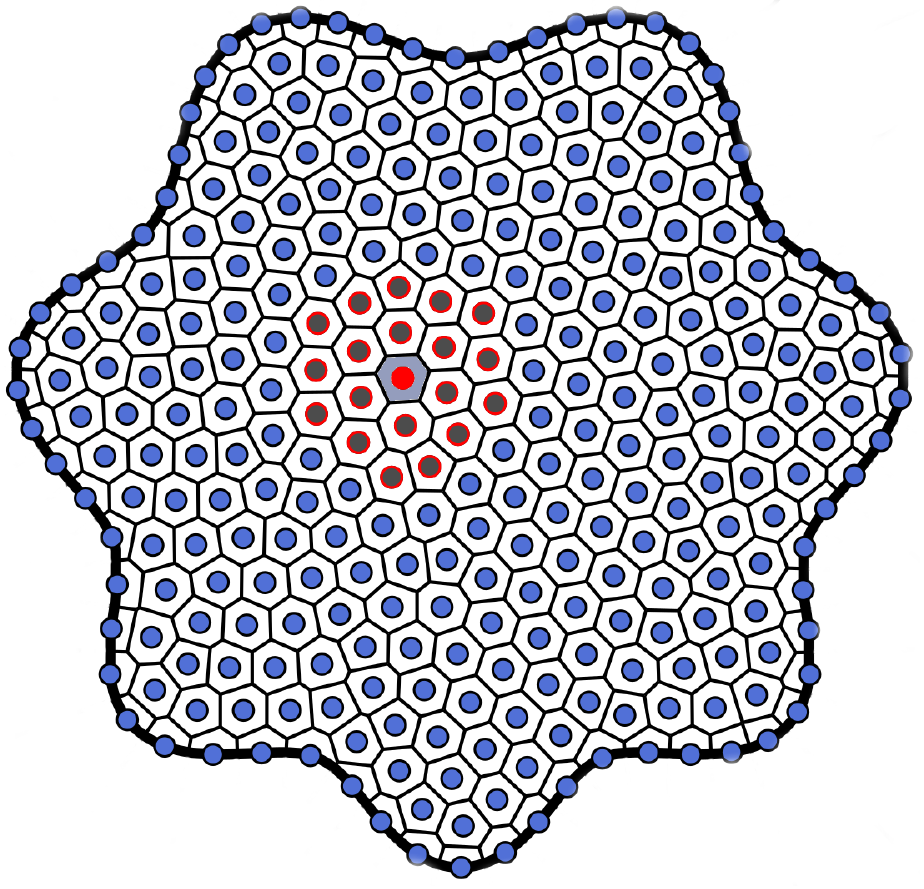}
\includegraphics[width=0.4\linewidth]{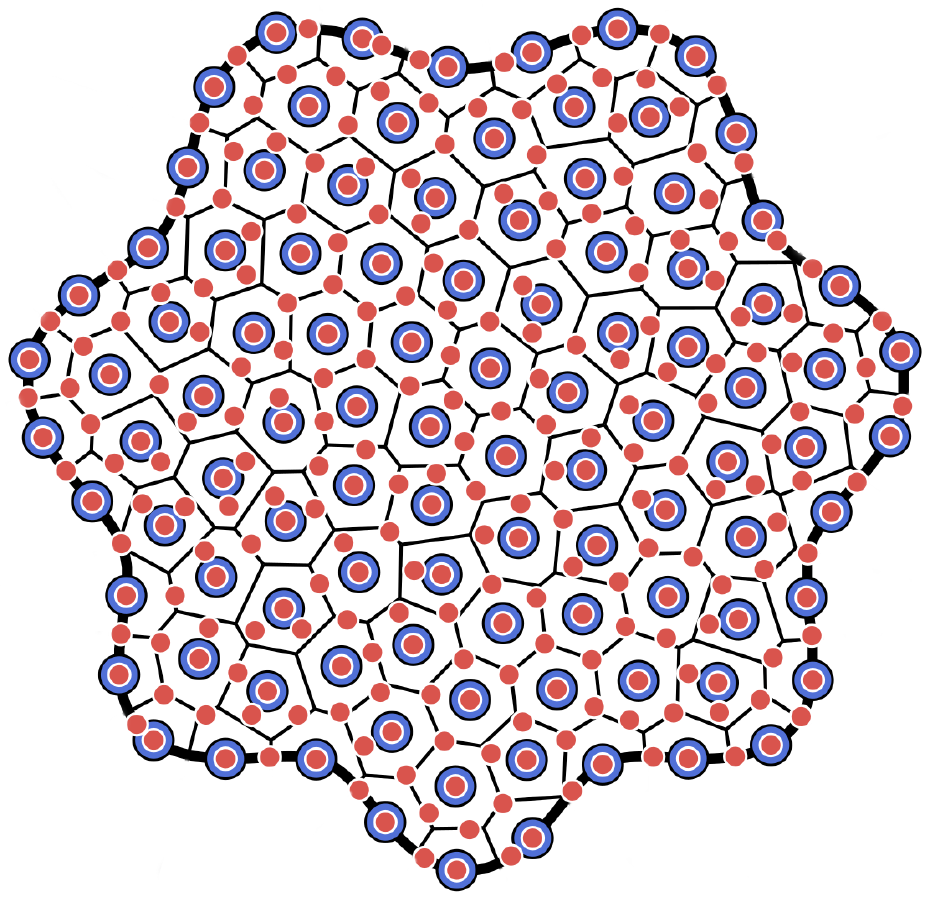}
\caption{Left: an example of $X$ points (blue points with black edges) that discretize the computational domain $\Omega$, 
together with one stencil (black points with a red edge) centered around the point filled with a red color. Right: 
an example of $Y$ points (small red points) added on top of $X$ points, together with a Voronoi diagram, where each cell 
is centered around a corresponding $X$ point.}
\label{fig:discretization:rbffd_stencil}
\end{figure}

We use the RBF-FD method in order to generate the RBF-FD trial space, which is spanned by 
a set of compactly supported cardinal basis functions $\{\Psi_i(y)\}_{i=1}^N$. These functions are constructed 
through a sequence of interpolation problems over stencils. 
An example of a stencil over $\Omega$ is given in Figure \ref{fig:discretization:rbffd_stencil}.

The cardinal basis functions are used to compute a PDE solution $u_h(Y, t)$ and a derivative $\mathcal{L} u_h(Y, t)$, 
in the form:
\begin{equation}
    \label{eq:discretization:ansatz}
    u_h(y, t) = \sum_{i=1}^N u_h(x_i,t) \Psi_i(y),\quad \mathcal{L}u_h(Y, t) = \sum_{i=1}^N u_h(x_i,t)\, \mathcal{L}\Psi_i(y), 
\end{equation}
where $\mathcal{L}$ is a linear differential operator and $\{u_h(x_i,t)\}_{i=1}^N$ are the solution nodal values which are the unknowns.

Using $X = \{x_i\}_{i=1}^N$, we reformulate the relations in \eqref{eq:discretization:ansatz} using a semi-discrete matrix notation:
\begin{equation}
    \label{eq:discretization:semimatrixansatz}
    u_h(y, t) = E(y,X) u_h(X,t),\quad u_h^{\mathcal{L}}(y, t) = D^{\mathcal{L}}(y,X) u_h(X,t).
\end{equation}
Here $u_h(X,t)=[u_h(X_1,t), ..., u_h(X_N,t)]$ is the vector of nodal values and the terms 
$E(y,X) = [\Psi_1(y), ..., \Psi_N(y)]$ and $D^{\mathcal{L}}(y,X) = [\mathcal{L}\Psi_1(y), ..., \mathcal{L}\Psi_N(y)]$ 
are semi-discrete matrices constructed using the RBF-FD evaluation/differentiation weights, where 
each set of weights is generated specifically for a point $y$. 
An algorithmic discussion on how to generate the evaluation and differentiation weights is given in \cite{tominecDiaphragm2D}. 
A mathematical formulation for computing these weights is given in \cite{tominec2020unfitted,ToLaHe21}. 
In this paper we use these approaches to generate a vector of weights for every $y \in Y$, such that the weights are exact for a cubic polyharmonic spline basis and a 
monomial basis of degree $p$, which is a concept introduced in \cite{Barnett15} and studied in \cite{FFBB16,BFFB17,Bayona19,Jancic21}.
In our computations we use the stencil size $n = 2\binom{p+2}{2}$ for $p\geq 2$, and $n=15$ for $p=1$.
To form evaluation and differentiation matrices, we sample \eqref{eq:discretization:semimatrixansatz} in each $y_k \in Y$, $k=1,..,M$. 
This gives the following two relations:
\begin{equation}
    \label{eq:discretization:matrixansatz}
    u_h(Y, t) = E(Y,X) u_h(X,t),\quad u_h^{\mathcal{L}}(Y, t) = D^{\mathcal{L}}(Y,X) u_h(X,t),
\end{equation}
where $E(Y,X)$ is a rectangular evaluation matrix of size $M \times N$ and $D^{\mathcal{L}}(Y,X)$ is a rectangular differentiation matrix of size $M \times N$. 
The two matrices can easily be computed in MATLAB by using the code provided in \cite{Tominec_rbffdcode2021}. \igor{We note that when $M=N$ and $Y=X$, then the 
matrices are square, and the underlying discretization is of the collocation type. Otherwise, the underlying discretization is of the oversampled type.}

\subsection{Discretization of a conservation law}
\label{sec:discretization:advection}
Here we discretize \eqref{eq:model:conservationlaw}, where we have $k$ unknown functions that form 
a set $\bm U = \left ( u_1, u_2, .., u_k \right)$. 
Each unknown function in $\bm U$ is represented using the RBF-FD ansatz \eqref{eq:discretization:ansatz} such that: 
\begin{eqnarray*}
\sum_{i=1}^N \partial_t \bm U_h(x_i,t)\, \Psi_i(y) = - \sum_{i=1}^N F_1(\bm U_h(x_i,t))\, \nabla_1 \Psi_i(y) - \sum_{i=1}^N F_2( \bm U_h(x_i,t))\, \nabla_2 \Psi_i(y),
\end{eqnarray*}
where $\bm U_h(x_i, t) = \left( u_{h,1}(x_i, t), ... , u_{h,k}(x_i, t) \right )$ are the $i$-th unknown nodal values of the functions that we are solving the conservation law for, 
and where we used an interpolant of the flux $\bm F_l(\bm U_h(y,t)) = \sum_{i=1}^N \bm F_l(\bm U_h(x_i,t)) \Psi_i(y)$, $l=1,2$, to compute the divergence 
term in \eqref{eq:model:conservationlaw}.
Now we sample the equation above in every $y \in Y$ in order to obtain a system of $M$ equations with $N$ unknowns:
\begin{eqnarray}
    \label{sec:discretization:advection_system}
    \sum_{i=1}^N \partial_t \bm U_h(x_i,t)\, \Psi_i(y_1) &=& - \sum_{i=1}^N F_1(\bm U_h(x_i,t))\, \nabla_1 \Psi_i(y_1) - \sum_{i=1}^N F_2(\bm U_h(x_i,t))\,\, \nabla_2 \Psi_i(y_1) \nonumber \\
    &\vdots&  \\
    \sum_{i=1}^N \partial_t \bm U_h(x_i,t)\, \Psi_i(y_M) &=& - \sum_{i=1}^N F_1(\bm U_h(x_i,t))\, \nabla_1 \Psi_i(y_M) - \sum_{i=1}^N F_2(\bm U_h(x_i,t))\, \nabla_2 \Psi_i(y_M). \nonumber
\end{eqnarray}
The matrix-vector formulations for the two problems that we consider in this paper are collected in the subsections below.
\subsection{Scalar conservation law in matrix-vector format}
\label{sec:discretization:scalarLaw}
Discussion in this section is kept general in the sense that it applies to $\bm F$ for scalar conservation laws. 
In this case we have that $\bm U_h = u_h$ and $\bm F(u_h) = [F_1(u_h), F_2(u_h)]$. 
Below, we rewrite the system \eqref{sec:discretization:advection_system} using the matrices defined in \eqref{eq:discretization:matrixansatz} and \eqref{eq:discretization:semimatrixansatz}:
\begin{eqnarray}
    \label{eq:discretization:advection_tmp}
    E(Y,X)\, \partial_t u_h(X,t) &=& -D^{\nabla_1}(Y,X) F_1(u_h(X,t)) - D^{\nabla_2}(Y,X)  F_2(u_h(X,t)).
\end{eqnarray}
Here $E(Y,X)$, $D^{\nabla_1}(Y,X)$ and $D^{\nabla_2}(Y,X)$ are the rectangular evaluation and differentiation matrices of size $M \times N$ defined in Section \ref{sec:discretization:matrices}. 
The vectors $F_1$ and $F_2$ of size $M \times 1$ correspond to the components of the (non)linear flux, thus their elements are 
$F_{1,i}(u_h(X,t)) = F_1(u_h(x_i, t))$ and $F_{2,i}(u_h(X,t)) = F_2(u_h(x_i, t))$ 
respectively. 
Now we multiply \eqref{eq:discretization:advection_tmp} with $h_y$, the approximate mean distance between the $Y$ points 
to arrive at:
\begin{equation}
    \label{eq:discretization:advection}
    \bar E(Y,X)\, \partial_t u_h(X,t) = -\bar D^{\nabla_1}(Y,X)F_1(u_h(X,t)) - \bar D^{\nabla_2}(Y,X) F_2(u_h(X,t)).
\end{equation}
where $\bar E(Y,X) = h_y\, E(Y,X)$ and $\bar D(Y,X) = h_y\, D(Y,X)$.
The multiplication by $h_y$ is a
norm scaling introduced in \cite{ToLaHe21}, which makes it possible to couple our discrete 
problem to an equivalent continuous problem \cite{ToLaHe21}.

When employing the oversampling approach, the system of equations \eqref{eq:discretization:advection} is rectangular and 
can not be solved using an explicit time-stepping technique in the form it has, 
as the time derivative on the left hand side is not explicitly defined. 
For this reason we project the residual of \eqref{eq:discretization:advection} onto the column space of $\bar E(Y,X)$ by 
multiplication with $\bar E^T(Y,X)$, 
where for simplicity, we also drop the $(Y,X)$ notation and obtain:
\begin{equation}
    \label{eq:discretization:advection_tmp1}
    \bar E^T \bar E\, \partial_t u_h(X,t) = \bar E^T \left( -\bar D^{\nabla_1}F_1(u_h(X,t)) - \bar D^{\nabla_2} F_2(u_h(X,t)) \right).
\end{equation}
Looking at \eqref{eq:discretization:advection_tmp1}, the components of the left hand side matrix $\bar E^T \bar E$ are:
\begin{equation}
    \label{eq:discretization:discreteinner}
(\bar E^T \bar E)_{ij} = h_y^2 \sum_{k=1}^M \Psi_i(y_k) \Psi_j(y_k).
\end{equation}
In other words, every component of $\bar E^T \bar E$ is a discrete inner product $(\Psi_i(y_k), \Psi_j(y_k))$, implying that $\bar E^T \bar E$ is a mass matrix in this inner product. 
The same type of inner product also applies to the components of the right-hand-side matrix products. Thus, our projection using $\bar E^T$ has a 
tight connection with a classical 
Galerkin projection. A key difference is that our type of discretization uses discrete (summation) instead of a continuous (exact integration) inner product. 

Now we solve \eqref{eq:discretization:advection_tmp1} for the time derivative $\partial_t\, u_h(X,t)$ 
by inverting the mass matrix $\mathcal{M} = \bar E^T \bar E$, in order to get a system of ODEs ready to be advanced in time:
\begin{equation}
    \label{eq:discretization:advection_tmpp}
\begin{aligned}
    \partial_t u_h(X,t) &= \mathcal{M}^{-1}\, \bar E^T \left( -\bar D^{\nabla_1}F_1(u_h(X,t)) - \bar D^{\nabla_2} F_2(u_h(X,t)) \right) \\
    &\equiv \mathcal{M}^{-1} \bar E^T \bar D(F(u_h(X,t))).
\end{aligned}
\end{equation}
\igor{To reduce the computational cost, we, instead of directly inverting $\mathcal{M}$, solve the system for $\partial_t u_h(X,t)$ iteratively,}
using the conjugate gradient method (function \texttt{pcg()} in MATLAB). \igor{Note that we have $\mathcal{M} = I$ in the collocation case, which implies that 
the system of equations does not have to be solved in each time step.}

To \eqref{eq:discretization:advection_tmpp} we also add two stabilization terms; (i) the term $P_1$ to stabilize the system of ODEs in time 
using a hyperviscosity operator, 
(ii) the term $P_2$ to treat the discontinuities which appear when using a nonlinear flux $\bm F(u_h)$. The scheme is then:
\begin{equation}
    \label{eq:discretization:advection_final}
    \partial_t u_h(X,t) = \mathcal{M} ^{-1}\, \Big [\bar E^T \bar D(F(u_h(X,t))) + P_1 u(X,t) + P_2 u(X,t) \Big].
\end{equation}
This is a discretization of the scalar conservation law in conservative form.
In this work we formulate the hyperviscosity term $P_1 \approx \Delta^2$ as:
\begin{equation}
    \label{eq:discretization:hypervi}
    \begin{aligned}
    P_1 &= \gamma\, h_y^2\, \left[D^{\Delta_{11}}(Y,X) + D^{\Delta_{22}}(Y,X)\right]^T\left[D^{\Delta_{11}}(Y,X) + D^{\Delta_{22}}(Y,X)\right] \\
    \gamma &= h^{4.5},
    \end{aligned}
\end{equation}
where $h_y^2$ is the norm scaling and $\gamma$ is the hyperviscosity scaling. Our choice of $\gamma$ is justified by 
making the $\gamma P_1$ term consistent with the PDE up to the $4$-th order, 
and then increasing that number to $4.5$ in order to allow a larger time step when using an explicit time-stepping algorithm.
Term $P_2$ is constructed using the residual viscosity concept described in Section \ref{sec:residualviscosity}.

Note that when the flux $\bm F(u_h)$ is linear, then the conservative and the non-conservative forms of a conservation laws 
exhibit the same numerical properties, and the ODE system \eqref{eq:discretization:advection_final} transforms to:
\begin{equation}
    \label{eq:discretization:advection_linear_final}
    \partial_t u_h(X,t) = \mathcal{M}^{-1}\, \bar E^T \left( -F_1' \bar D^{\nabla_1} - F_2' \bar D^{\nabla_2}  \right)u_h(X,t) - (P_1+P_2) u_h(X,t) \equiv \bar D u_h(X,t).
\end{equation}
Here $F_1'$ and $F_2'$ correspond to the components of the velocity field $\bm F'(t)$. More specifically, 
they represent rectangular identity matrices of size $M \times N$ with components $F'_{1_{ii}} = F'_1(t)$, $F'_{2_{ii}} = F'_1(t)$, $i=1,..,N$. 

At this point we can use an explicit time-stepping algorithm to advance the solution of 
\eqref{eq:discretization:advection_final} from $t=0$ to $t=T$. 
Throughout this paper we use an explicit classical Runge-Kutta 4 time-stepping algorithm.
In the numerical experiment section of the paper we solve problems with Dirichlet-type boundary conditions on the inflow boundary. 
To impose these conditions we use a so-called injection method, 
that is, we step the discretization \eqref{eq:discretization:advection_final} 
and then in every time step overwrite the correct subset of boundary elements of $u_h(X,t)$ with the corresponding values of the Dirichlet boundary condition.
After the nodal solution $u(X,T)$ at the final time $t=T$ is obtained, we evaluate the solution at $Y$ points as:
\begin{equation}
    \label{eq:discretization:finalEvaluation}
    u_h(Y,T) = E(Y,X)\, u_h(X,T).
\end{equation}

\subsection{System of conservation laws in matrix-vector format: compressible Euler equations}
Consider compressible Euler equations in two spatial dimensions.
The unknowns $\bm U$ and the flux $\bm F(\bm U)$ are in this case given by:
\begin{equation}
    \label{eq:discretization:Euler_variables}
 \bm U = 
 \begin{pmatrix}
\rho \\
\bm m \\
\mathcal{E}
\end{pmatrix},
\quad 
\mathbf{F(U)}
=
\begin{pmatrix}
\bm m \\
\bm m \otimes \bm v + p\, \bm I \\
\bm v (\mathcal{E} + p)
\end{pmatrix}.
\end{equation}
Here $\bm I$ is an identity matrix, $\rho$ is the density, $\bm m = (m_1, m_2)$ is the momentum, 
$\mathcal{E}$ is the total energy and $p$ is the pressure of the fluid.  The 
relation between the momentum and the velocity field is $\bm v = \bm m / \rho$.
When forming a system of equations we have $5$ unknown functions, but only $4$ equations. 
We close the system of equations using an ideal gas condition, 
by defining the pressure variable as:
$p = \rho T = (c_{\text{ad}} -1)(\mathcal{E} - \rho |\bm v |^2/2),$
where $T$ is the temperature of the fluid and $c_{\text{ad}}$ is the adiabatic gas constant.
In this case the expanded form of the nonlinear conservation law \eqref{eq:model:conservationlaw} is:
\begin{equation}
\begin{aligned}
    \partial_t
\begin{pmatrix}
    \rho \\
    m_1 \\
    m_2 \\
    \mathcal{E}
\end{pmatrix}
=
\begin{pmatrix}
- \nabla_1 m_1       - \nabla_2 m_2 \\
- \nabla_1 (m_1 v_1) - \nabla_2 (m_2 v_1) - \nabla_1 p \\
- \nabla_1 (m_1 v_2) - \nabla_2 (m_2 v_2) - \nabla_2 p \\
- \nabla_1 (v_1 (\mathcal{E} + p)) - \nabla_2 (v_2 (\mathcal{E}+p))
\end{pmatrix}
\end{aligned}
\end{equation}
We discretize the above system by introducing unknown nodal values for all unknown functions $\rho(X,t), m_1(X,t), m_2(X,t), E(X,t)$. 
The oversampled RBF-FD discretization analogous to \eqref{eq:discretization:advection_final} is then:
\begin{equation}
    \label{eq:discretization:Euler_final}
    \begin{aligned}
        \partial_t
        \begin{pmatrix}
            \rho(X,t) \\
            m_1 (X,t) \\
            m_2 (X,t) \\
            \mathcal{E} (X,t) 
        \end{pmatrix}
         =
         \mathcal{M}^{-1}\,         
        \begin{pmatrix}
            E^T G_\rho(Y,t) + P \rho(X,t) \\
            E^T G_{m_1}(Y,t) + P m_1 (X,t) \\
            E^T G_{m_2}(Y,t) + P m_2 (X,t) \\
            E^T G_\mathcal{E}(Y,t) + P \mathcal{E}(X,t) 
        \end{pmatrix}       
    \end{aligned}
\end{equation}
where:
\begin{equation}
    \label{eq:discretization:Euler_RHS_flux}
\begin{aligned}
    G_\rho(Y,t) &=& - D^{\nabla_1} m_1(X,t) - D^{\nabla_2} m_2(X,t) \\
    G_{m_1}(Y,t) &=& - D^{\nabla_1} \big (m_1(X,t)\, v_1(X,t) \big ) - D^{\nabla_2} \big (m_2(X,t)\, v_1(X,t) \big ) - D^{\nabla_1} p \\
    G_{m_2}(Y,t) &=& - D^{\nabla_1} \big (m_1(X,t)\, v_2(X,t) \big ) - D^{\nabla_2} \big (m_2(X,t)\, v_2(X,t) \big ) - D^{\nabla_2} p  \\
    G_\mathcal{E}(Y,t) &=& - D^{\nabla_1} \big ( v_1(\mathcal{E}+p) \big ) -  D^{\nabla_2} \big ( v_2(\mathcal{E}+p) \big ) \\
    P &=& - P_1 - P_2,
\end{aligned}
\end{equation}
with differentiation matrices $D^{\nabla_1}, D^{\nabla_2}$ defined in Section \ref{sec:discretization:matrices}, and stabilizers $P_1$ and $P_2$ described in Section \ref{sec:discretization:scalarLaw}.
The Dirichlet-type boundary conditions for a given unknown function are imposed using the injection method, just as described for the scalar conservation law 
in Section \ref{sec:discretization:scalarLaw}. After the system of ODEs is advanced to the final time $t=T$, each unknown set of nodal values is evaluated 
at the $Y$ points, 
analogously to \eqref{eq:discretization:finalEvaluation}.

\section{Residual based artificial viscosity stabilization of shocks for nonlinear conservation laws}
\label{sec:residualviscosity}
Solutions to nonlinear conservation laws \eqref{eq:model:conservationlaw} are discontinuous, when the flux $\bm F (\bm U)$ does not include any physical viscous forces. 
Numerical approximation of discontinuous functions suffers from the Gibbs phenomenon. 
Based on \cite{Nazarov13,NazarovHoffman13,NazarovLarcher17} we formulate a residual based viscosity 
stabilization of the solution in the 
context of the RBF-FD method.

The numerical scheme is stabilized using a viscosity term $P_2 \approx \nabla \cdot (\varepsilon \nabla u)$, where $\varepsilon$ is 
a spatially variable coefficient, computed such that the viscous term is prevalently active in the regions of discontinuities.
The viscosity term is below discretized by using the differentiation matrices introduced in Section \ref{sec:discretization:matrices}:
\begin{equation}
    \label{eq:resviscosity:P2}
P_2 = (D^{\nabla_1})^T\, \diag (E\, \varepsilon(X, t_n))\, D^{\nabla_1}  + (D^{\nabla_2})^T\, \diag (E\, \varepsilon(X, t_n))\, D^{\nabla_2}.
\end{equation}
The definition of $\varepsilon$ for scalar conservation laws and the Euler system, is given in the following two subsections. 

\subsection{Definition of the viscosity coefficient for scalar conservation laws}
The discrete problem which we aim to stabilize is \eqref{eq:discretization:advection_final}. The discussion 
is centered around the stabilization term $P_2$ defined in \eqref{eq:resviscosity:P2}.
For each node $x_i \in X$, $i=1,..,N$, 
we define the viscosity coefficient vector $\varepsilon$ included inside the definition of $P_2$ as:
\begin{equation}
    \label{eq:resviscosity:epsilon}
    \begin{aligned}
    \varepsilon(x_i, t_{n+1}) &=& \min \left(\varepsilon_{\text{RV}}(x_i, t_{n}),\, \varepsilon_{\text{UW}}(x_i, t_{n}) \right), &\quad n=2,3,... \\
    \varepsilon(x_i, t_{n+1}) &=& \varepsilon_{\text{UW}}(x_i, t_{n}), &\quad n=1.
    \end{aligned}
\end{equation} 
Here $\varepsilon_{\text{RV}}$ and $\varepsilon_{\text{UW}}$ are the residual and the upwind coefficients respectively. Let
 $K_i$ be a Voronoi cell  with a corresponding Voronoi center $x_i$, and $y_j \in Y$ an evaluation point (see the right plot in Figure \ref{fig:discretization:rbffd_stencil}). Then the definitions of $\varepsilon_{\text{RV}}$ and $\varepsilon_{\text{UW}}$ 
 are \cite{Nazarov13}:
\begin{equation}
    \label{eq:resviscosity:epsilon_rv_uw}
    \begin{aligned}
\varepsilon_{\text{RV}}(x_i, t_{n+1}) &= C_{\text{RV}}\, h_{\text{loc}}^2(x_i) \max_{y_j \in K_i} |R(y_j, t_n)|\, \frac{1}{n(x_i)}, \\  
\varepsilon_{\text{UW}}(x_i,t_{n+1}) &= \frac{1}{2}\, h_{\text{loc}}(x_i)\, \max_{y_j \in K_i} \sqrt{F'_1(u_h(y_j,t_{n}))^2 + F'_2(u_h(y_j,t_{n}))^2 }.
    \end{aligned}
\end{equation}
Here we define $h_{\text{loc}}(x_i)$ as the minimum pairwise distance between points in a patch centered around $x_i$, where the patch consists of $5$ points closest to $x_i$ \cite{stiernstrom2021}.
The temporal solution $u_h(y_j,t_{n})$ at a point $y_j \in Y$  comes from the relation $u_h(Y,t_n) = E(Y,X)\, u_h(X,t_n)$. 
The user-defined constant $C_{\text{RV}}$ is of size $\mathcal{O}(1)$ and is independent of $h_{\text{loc}}(x_i)$.  
As $C_{\text{RV}} \to \infty$, then $\varepsilon_{\text{RV}} \to \infty$, and $\varepsilon \to \varepsilon_{\text{UW}}$. On the other hand, 
as $C_{\text{RV}} \to 0$, then $\varepsilon_{\text{RV}} \to 0$, and $\varepsilon \to 0$.
The term $n(x_i, t_n)$ is the residual normalization defined by:
\begin{equation}
    \label{eq:resviscosity:normalization}
    \begin{aligned}
n(x_i, t_n) &= \Big | u_{\text{loc}} - \| u_h(X,t_n) - \bar u_h(X,t_n) \|_\infty \Big |, \\
u_{\text{loc}} &= \max_{y_j \in K_i}  u_h(y_j,t_n) - \min_{y_j \in K_i}  u_h(y_j,t_n),
\end{aligned}
\end{equation}
where $\bar u_h(X,t_n)$ is the mean of the temporal solution $u_h (X,t_n)$. The role of the normalization is to unify the physical units of $\varepsilon_{\text{RV}}$ and $\varepsilon_{\text{UW}}$. 
The residual $R(Y,t_n)$ is defined as:
\begin{equation}
    \label{eq:resviscosity:residual}
    R(Y,t_n) = E \left(D^{\partial_t} u_h(X,t_n) \right ) +  D^{\nabla_1} F_1(u_h(X,t_n)) + D^{\nabla_2} F_2(u_h(X,t_n)), 
\end{equation}
where $D^{\nabla_1} F_1(u_h(X,t_n)) + D^{\nabla_2} F_2(u_h(X,t_n))$ is the flux divergence of the solution at time $t_n$ computed using the RBF-FD differentiation matrices introduced in 
Section \ref{sec:discretization:matrices}. Furthermore, $E$ is an evaluation matrix also introduced in Section \ref{sec:discretization:matrices}, and the term 
$D^{\partial_t} u_h(X,t_n)$ is an approximation of the time derivative at $t_n$ using the already computed solution from previous time points $t_n, t_{n-1}, t_{n-2}, ...$. 
We construct $D^{\partial_t}$ in two different ways: (i) when the solution is advanced in time using a constant $\Delta t$ we use the backward-differentiation formulae (BDF), 
(ii) when the solution is advanced in time by a variable $\Delta t_i$, $i=1,2,..$, then we construct an approximation to a time derivative using the polynomial basis, in each time step. 
In any case, the approximation order in (i) or (ii) has to match the approximation order of the spatial discretization.  At the beginning of the simulation we can only 
use a few points from the past temporal solutions. When a sufficient amount of temporal solutions have been computed we are allowed to use a high-order formula. 
We then write that:
\begin{equation}
    \label{eq:resviscosity:BDF}
D^{\partial_t} \big |_{t_n} 
=   
\begin{cases}
    D^{\partial_t}_{1} & n=2 \\
    D^{\partial_t}_{2} & n=3 \\
    D^{\partial_t}_{3} & n=4 \\
    D^{\partial_t}_{4}  & n \geq 5, \\
\end{cases}
\end{equation}
where $D^{\partial_t}_{1}$ is a differentiation operator of order $1$, $D^{\partial_t}_{2}$ of order two, and so on.
In the case (i) we use the BDF4 formula to define $D^{\partial_t}_{\bar k}$, for $\bar k=1,...,4$ consecutively as:
\begin{eqnarray}
    \label{eq:resviscosity:BDF_detailed}
\frac{1}{\Delta t} \left(   U(X,t_n) - U(X,t_{n-1})   \right ), \\
\frac{1}{\Delta t} \frac{3}{2} \left(   U(X,t_n) - \frac{4}{3} U(X,t_{n-1}) + \frac{1}{3} U(X,t_{n-2})   \right ), \nonumber \\
\frac{1}{\Delta t} \frac{11}{6} \left(   U(X,t_n) - \frac{18}{11} U(X,t_{n-1}) + \frac{9}{11} U(X,t_{n-2})  - \frac{2}{11} U(X,t_{n-3})  \right ), \nonumber \\
\frac{1}{\Delta t} \frac{25}{12} \left(   U(X,t_n) - \frac{48}{25} U(X,t_{n-1}) + \frac{36}{25} U(X,t_{n-2})  - \frac{16}{25} U(X,t_{n-3}) + \frac{3}{25} U(X,t_{n-4})  \right ). \nonumber
\end{eqnarray}
In the case (ii) we construct $D^{\partial_t}_{\bar k}$, $\bar k=1,...,4$ by computing a set of differentiation weights through a 1D polynomial interpolant $u(t) = \sum_{j=1}^{\bar k+1} \xi_i t^{j-1}$, 
where $\bar k+1$ is the number of the time points in which we already know the solution 
$u(t_1), u(t_2), ..., u(t_{k+1})$, 
and $\bar k$ is the desired order of the approximation. Then we sample $u(t)$ in $t_n, ..., t_{n-\bar k}$ to obtain a system of equations 
$\underline{u} = A\, \underline{\xi}$, where $A_{ij} = t_i^{j-1}$ and $\xi_i$, $j=1,..,\bar k + 1$, $i=1,.., \bar k + 1$  are the interpolation matrix and 
the unknown coefficients respectively. 
A time-derivative of $u(t)$ is $\partial_t u(t) = \sum_{j=1}^{\bar k + 1} (j-1) t^{j-2} \xi_j$, an equivalent vector formulation is 
$\partial_t u(t) = \underline{b}(t)\, \underline{\xi}$, 
where $b_j(t) = (j-1) t^{j-2}$. Plugging the computed $\underline{\xi}$ into $\partial_t u(t)$ we have 
$\partial_t u(t) = \underline{b}(t) A^{-1} \underline u \equiv \underline{w}(t) \underline u$. Here $\underline{w}(t)$ is a 
final vector of $k+1$ weights for computing a derivative in time when the time grid is non-uniform. The MATLAB code to compute $\underline w$ is provided 
in Appendix \ref{sec:appendix:timederivatives}.

Note that the time derivative approximated using the approaches (i) and (ii), is used only with the purpose to compute the residual, once the solution at $t_n$ has already been 
computed. The time stepping scheme that is actually advancing the solution in time 
is chosen differently.

\subsection{Definition of the viscosity coefficient for the Euler system}
Here we formulate the viscosity coefficient $\varepsilon$ for the Euler system of equations \eqref{eq:discretization:Euler_final}. 
The $\varepsilon$ coefficient is computed using the same relation as in \eqref{eq:resviscosity:epsilon}, however, the definitions of $\varepsilon_{\text{UW}}$ 
and $\varepsilon_{\text{RV}}$ in \eqref{eq:resviscosity:epsilon} are different. The upwind viscosity coefficient at time $t_{n+1}$ and a point $x_i \in X$ is defined as:
\begin{equation}
\varepsilon_{\text{UW}}(x_i, t_{n+1}) = \frac{1}{2}\, h_{\text{loc}}(x_i)\, \left (\sqrt{v_1^2(x_i,t_n) + v_2^2(x_i,t_n)} + \sqrt{c_{\text{ad}} T(x_i,t_n)} \right ),
\end{equation}
where the term in parenthesis is the local wave speed computed from the eigenvalues of $\frac{\partial \bm F}{\partial \bm U}$, 
and where $v_1(x_i,t_n) = m_1(x_i,t_n)/\rho(x_i,t_n)$, 
$v_2(x_i,t_n) = m_2(x_i,t_n)/\rho(x_i,t_n)$ are the velocities in horizontal and vertical directions respectively, and $T = p(x_i,t)/\rho(x_i,t)$ is the temperature of the fluid, 
and the internodal distance $h_{\text{loc}(x_i)}$ is defined in the scope of \eqref{eq:resviscosity:epsilon_rv_uw}. 
The coefficient $\varepsilon_{\text{RV}}$ is defined by:
\begin{equation}
    \label{eq:residualviscosity:euler_ep_rv}
\varepsilon_{\text{RV}} = C_{\text{RV}}\, h^2_{\text{loc}}(x_i)\, \max \left( \frac{|R_\rho|}{n_\rho}, \frac{|R_{m_1}|}{n_{m_1}}, \frac{|R_{m_2}|}{n_{m_2}}, \frac{|R_\mathcal{E}|}{n_\mathcal{E}} \right ),
\end{equation}
where the residuals $R_{*}$ are defined for each equation in \eqref{eq:discretization:Euler_final}, analogously to the formulation for the scalar conservation law in \eqref{eq:resviscosity:residual}. 
We have:
\begin{equation}
    \begin{aligned}
R_\rho(Y,t_n) &= E \left(D^{\partial_t}_{\text{BDF}} \rho(X,t_n) \right) - G_\rho(Y,t_n) \\
R_{m_1} (Y,t_n) &= E \left(D^{\partial_t}_{\text{BDF}} m_1(X,t_n) \right) - G_{m_1}(Y,t_n) \\
R_{m_2} (Y,t_n) &= E \left(D^{\partial_t}_{\text{BDF}} m_2(X,t_n) \right) - G_{m_2}(Y,t_n) \\
R_{\mathcal{E}} (Y,t_n)   &= E \left(D^{\partial_t}_{\text{BDF}} \mathcal{E}(X,t_n) \right) - G_\mathcal{E} (Y,t_n).
    \end{aligned}
\end{equation}
Here $G_\rho$, $G_{m_1}$, $G_{m_2}$ and $G_\mathcal{E}$ are defined in \eqref{eq:discretization:Euler_RHS_flux}, and 
the time derivative approximation $D^{\partial_t}_{\text{BDF}}$ is the same as defined in 
\eqref{eq:resviscosity:BDF} and \eqref{eq:resviscosity:BDF_detailed}. Each normalization $n_{*}$ in \eqref{eq:residualviscosity:euler_ep_rv} follows the same definition as 
in \eqref{eq:resviscosity:normalization}. For example, $n_\rho$ is defined as:
\begin{equation}
    \begin{aligned}
    n_\rho(x_i, t_n) &= \Big | \rho_{\text{loc}} - \| \rho(X,t_n) - \bar \rho(X,t_n) \|_\infty \Big |,\\
    \rho_{\text{loc}} &= \max_{y_j \in K_i}  \rho(y_j,t_n) - \min_{y_j \in K_i}  \rho(y_j,t_n), \quad i=1,..,N
    \end{aligned}
\end{equation}
where $\bar \rho(X,t_n)$ is the mean density at $t=t_n$ and $K_i$ is a Voronoi cell centered around $x_i \in X$.

\section{Numerical study I: linear advection problem}
\label{sec:experiments:linearadvection}
Throughout this section we compare collocation and oversampled RBF-FD methods, and investigate how do the hyperviscosity term 
\eqref{eq:discretization:hypervi} and the residual viscosity term \eqref{eq:resviscosity:P2} influence the numerical solution. 
Consider $\bm F$ is a linear flux with time-dependent velocity field given by:
\begin{equation}
    \label{eq:experiments_linear:velocityfield}
\bm F'(u) = \bm F' = \left(\cos(2 \pi\, t), \sin(2 \pi\,t) \right).
\end{equation}
The boundary of $\Omega$ is defined using polar coordinates $(r,\theta)$, where: $r = 1+\frac{1}{10}(\sin(7\theta)+ \sin(\theta))$ for $\theta \in [0, 2\pi]$.
We use the explicit classical Runge-Kutta 4 method to advance the solution in time. The constant time step is chosen as:
\begin{equation}
    \label{eq:experiments_linear:timestep}
\Delta t = \text{CFL}\, \min_{x_i \in X}\, \frac{h_{\text{loc}}(x_i)}{\sqrt{(F_1'(t))^2 + (F_2'(t))^2}} = \text{CFL}\, \min_{x_i \in X}\, h_{\text{loc}}(x_i),
\end{equation}
where we also used that $(F_1'(t))^2 + (F_2'(t))^2 = 1$, since $\bm F'$ is in \eqref{eq:experiments_linear:velocityfield} chosen as the rotational velocity field. The term $h_{\text{loc}}(x_i)$ is defined in the scope of \eqref{eq:resviscosity:epsilon_rv_uw}.
The constant $\text{CFL}$ is the Courant-Friedrichs-Lax number, specified in each subsection separately.
The relative approximation errors in $2$- and $1$-norm at final time $t=T$ are computed as:
\begin{equation}
    \label{eq:experiments_linear:errors}
\|e\|_2 = \frac{\|E(Y,X)\, u_h(X,T) - u(Y,T)\|_2}{\|u(Y,T)\|_2},\quad \|e\|_1 = \frac{\|E(Y,X)\, u_h(X,T) - u(Y,T)\|_1}{\|u(Y,T)\|_1},
\end{equation}
where $u(Y,T)$ is the exact solution at $t=T$, sampled at $Y$ points.

\subsection{Smooth initial condition: approximation error and stability properties}
Oversampled methods are normally associated with improved stability properties. 
In this section we investigate whether this is also true when solving hyperbolic problems, i.e., 
if there is a positive effect of oversampling on the eigenvalue spectrum of the discretized advection operator.

\begin{figure}
    \centering
\includegraphics[width=0.4\linewidth]{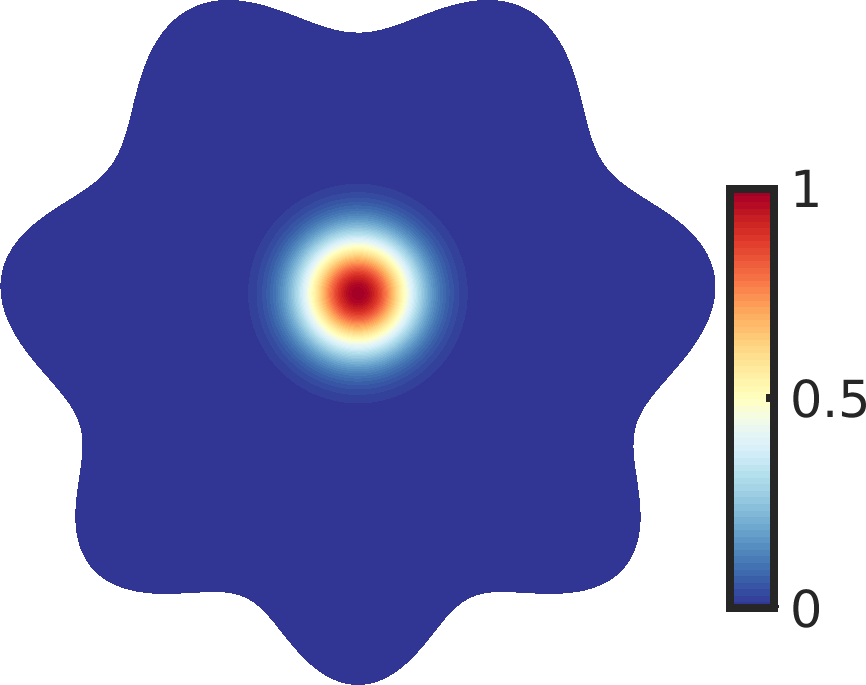}
\includegraphics[width=0.4\linewidth]{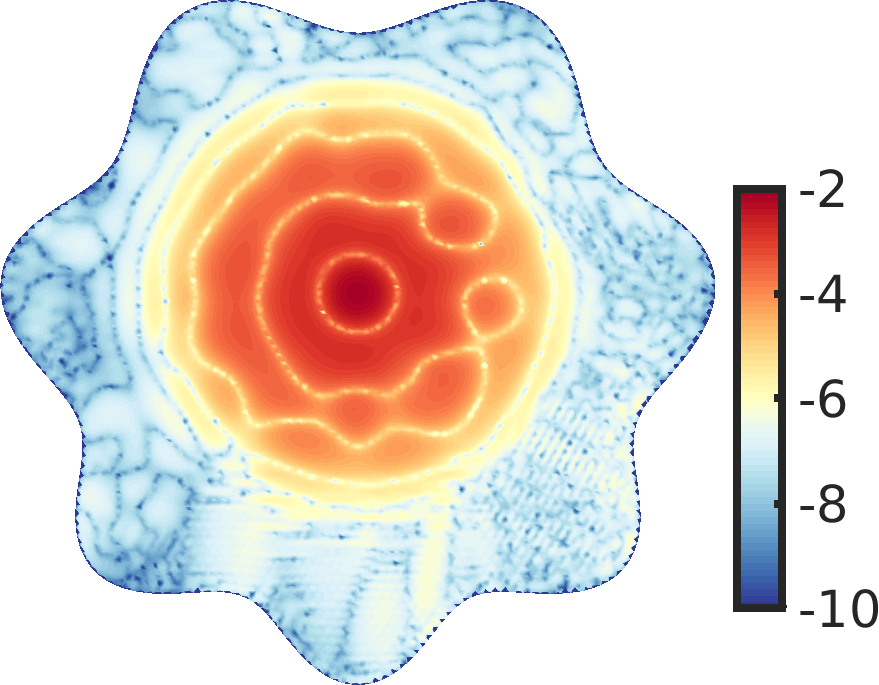}
\caption{Linear advection of a compactly supported $C^6$ function using a rotational velocity field. Initial condition is given in the left image. 
The spatial error distribution in the $\log_{10}$ scale, after the initial condition is rotated five times, and an oversampled RBF-FD method is in time stabilized using 
hyperviscosity, is given in the image on the right. Method parameters are set to $h=0.03$ and $p=4$ and $q=5$.}
\label{fig:linearadvection:solution_spatial_err}
\end{figure}

\begin{figure}
    \begin{tabular}{cccc}
        \hspace{0.2cm}\textbf{Collocation (q=1)} & \hspace{0.4cm} $\mathbf{q=5}$ & \hspace{1cm} $\mathbf{q=20}$ & \hspace{1cm}$\mathbf{q=30}$ \\
        \includegraphics[width=0.3\linewidth]{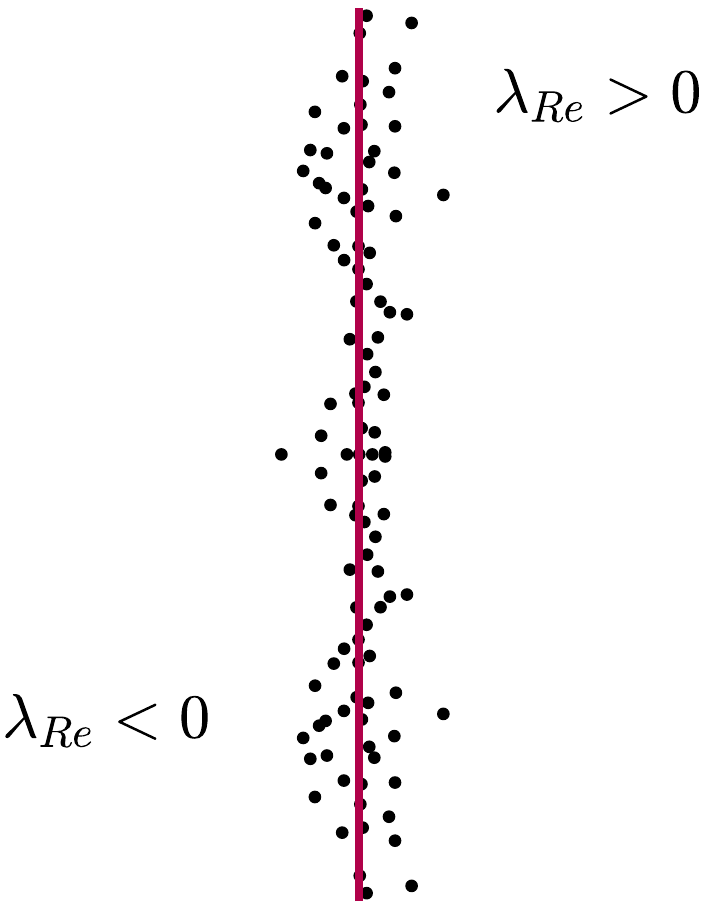} &
        \hspace{0.5cm} \includegraphics[width=0.067\linewidth]{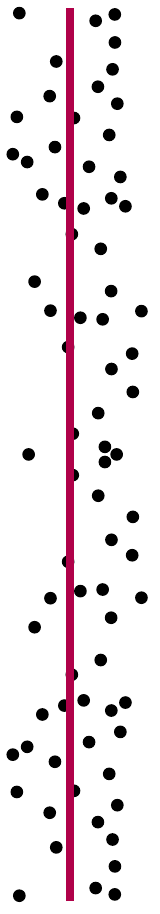} &
        \hspace{1cm} \includegraphics[width=0.115\linewidth]{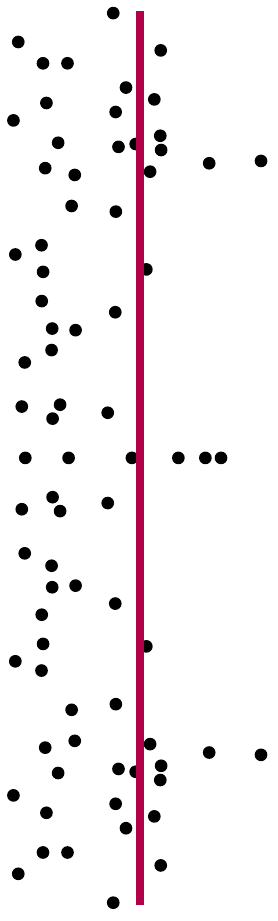}  &      
        \hspace{1cm} \includegraphics[width=0.075\linewidth]{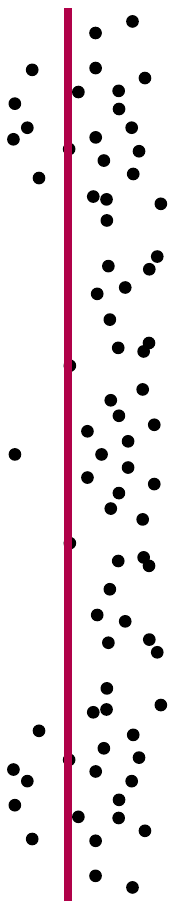}
    \end{tabular}
    \caption{Linear advection of a compactly supported $C^6$ function using a rotational velocity field: a close-up of the eigenvalue spectra for different 
    choices of the oversampling parameter $q$, when none of the cases is stabilized using hyperviscosity. Eigenvalues left of the red line have a negative 
    real part, while the eigenvalues on the right of the red line have a positive real part.}
    \label{fig:linearadvection:eigenvalues_notstabilized}    
\end{figure}

Provided that $d = d(y_1,y_2) = \sqrt{(y_1-0.2)^2 + (y_2-0.2)^2}$, the initial condition at $t=0$ is a $C^6$ compactly supported Wendland function 
centered in a point $(0.2, 0.2)$:
\begin{equation}
    \label{eq:experiments_linear:initCond_Wendland}
u(y,0) = 
\begin{cases}
(1-d)^6\, (35 d^2 + 18d + 3), & d \leq 0.4 \\
0, & \text{otherwise}.
\end{cases}
\end{equation}
The initial condition is plotted in Figure \ref{fig:linearadvection:solution_spatial_err}.
The discretized advection operator $D$ is given in \eqref{eq:discretization:advection_linear_final}. In this section \igor{we only want to 
observe the effect of oversampling on the time stability of the discretization, and compare the results to the collocation case. The test case does not involve any shocks}, 
therefore, we set $P_2=0$. When studying 
eigenvalues, we enforce zero Dirichlet boundary condition exactly in $D$, by removing all rows and columns corresponding to the boundary unknowns.
In Figure \ref{fig:linearadvection:eigenvalues_notstabilized} we display a close up of the eigenvalue spectrum around the imaginary axis, 
when the discretization is not stabilized, that is, when the parameter $\gamma$ in $P_1$ included into \eqref{eq:discretization:advection_linear_final} 
is set to $0$. We observe that the spectrum of the collocation RBF-FD method includes eigenvalues with a positive real part. These eigenvalues 
cause spurious growth when the solution is advanced in time. However, the situation is not improved as we introduce an oversampled discretization: 
for the oversampling parameters $q=5$, $q=20$ and $q=30$, the eigenvalue spectra still include eigenvalues with a positive real part. We do not observe that the positive real 
part is asymptotically shifted towards the negative imaginary half-plane, as $q$ is increased. Our conclusion is that oversampling by itself does not 
stabilize the eigenvalue spectrum of the semi-discrete PDE problem when the RBF-FD method is used.

\begin{figure}
    \centering
    \begin{tabular}{cccc}
        \hspace{0.5cm}\textbf{Collocation} & \hspace{0.7cm} $\mathbf{q=5}$ & \hspace{0.7cm} $\mathbf{q=20}$ & \hspace{0.7cm} $\mathbf{q=30}$ \\
        \includegraphics[width=0.21\linewidth]{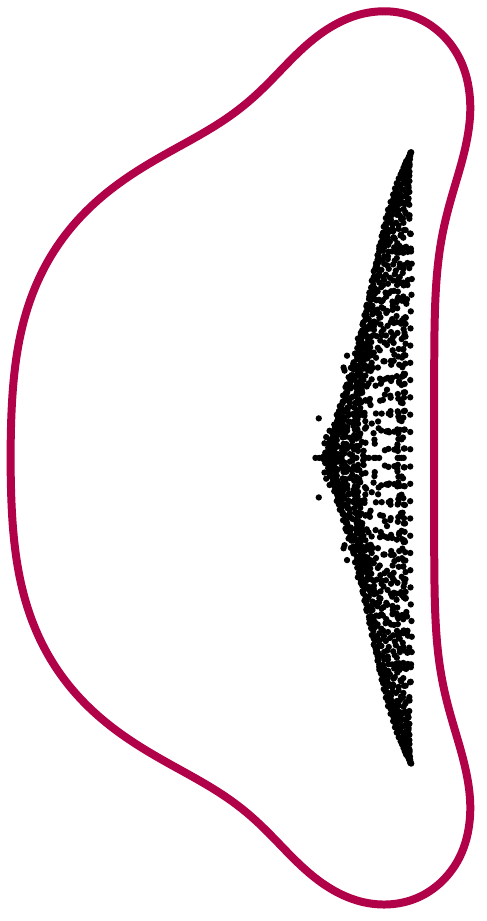} &
        \includegraphics[width=0.21\linewidth]{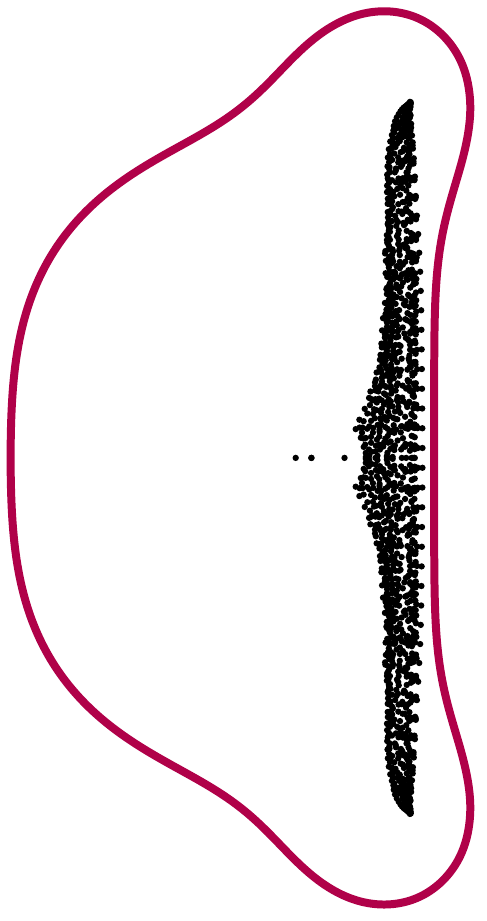} &
        \includegraphics[width=0.21\linewidth]{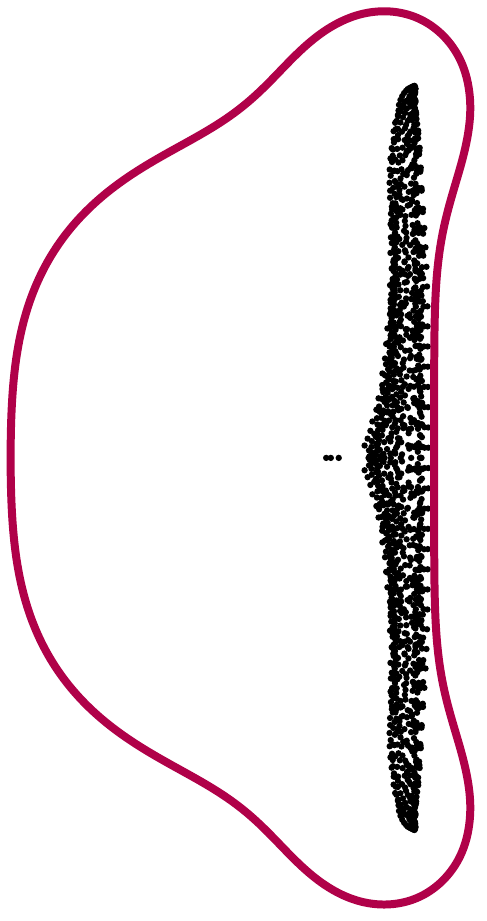} &
        \includegraphics[width=0.21\linewidth]{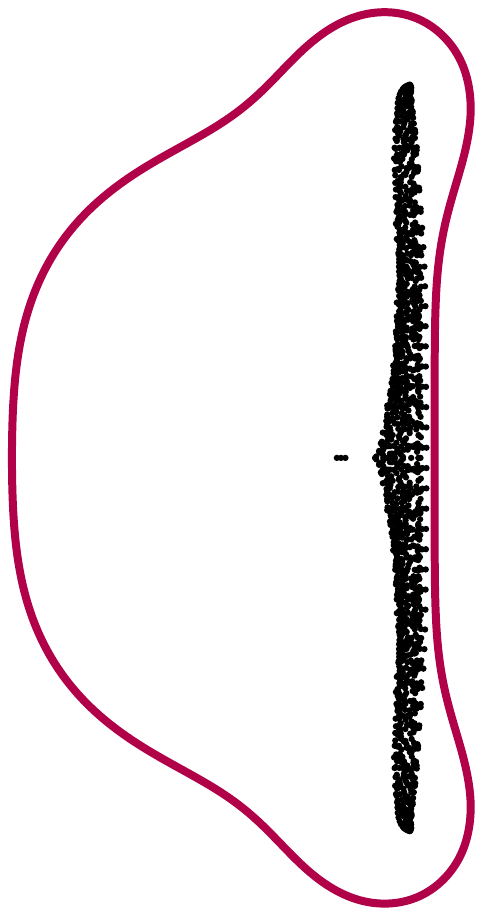} 
    \end{tabular}
    \caption{Linear advection of a compactly supported $C^6$ function using a rotational velocity field: eigenvalue spectra for different 
    choices of the oversampling parameter $q$, when the discretizations are stabilized in time using the hyperviscosity term. The closed red 
    line representsstability region of the Runge-Kutta 4 method for advancing the system of ODEs in time.}
    \label{fig:linearadvection:eigenvalues_hyperviscosity}    
\end{figure}

Next, we study the eigenvalue spectra when hyperviscosity is added to \eqref{eq:discretization:advection_linear_final}, that is, when 
the parameter $\gamma$ and the second-order hyperviscosity $P_1$ follow the definition in \eqref{eq:discretization:hypervi}. 
The results are given in Figure \ref{fig:linearadvection:eigenvalues_hyperviscosity}, where we observe that the eigenvalue spectra are stabilized. 
The eigenvalue distributions are similar across the different methods.

\begin{figure}
    \centering    
    \begin{tabular}{ccc}
        \multicolumn{3}{c}{\hspace{0.7cm} \textbf{Linear advection: collocation vs. oversampled RBF-FD methods}} \\
        \hspace{0.7cm} $\mathbf{p=2}$ & \hspace{0.7cm} $\mathbf{p=3}$ & \hspace{0.7cm} $\mathbf{p=4}$ \\
    \includegraphics[width=0.3\linewidth]{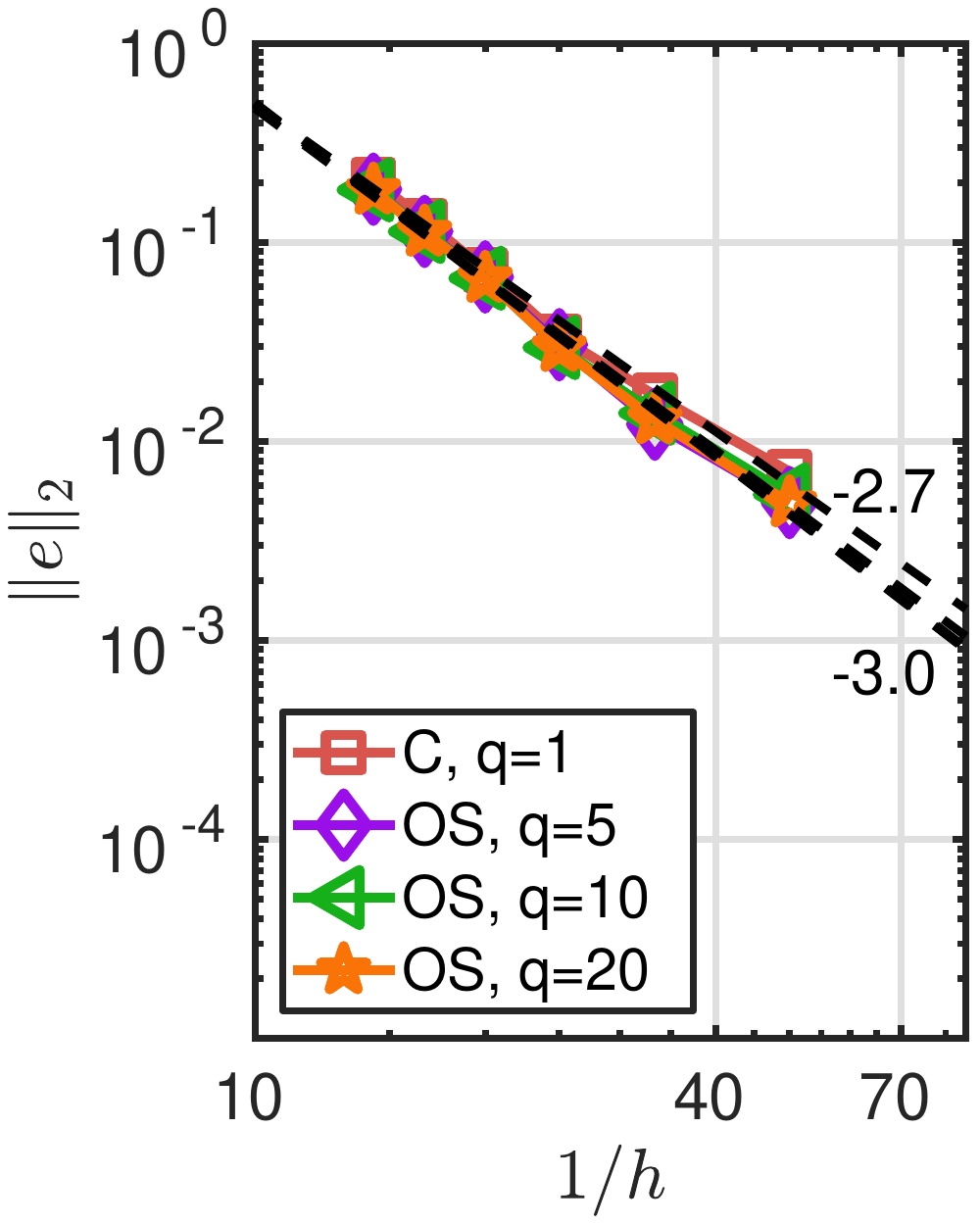} &
    \includegraphics[width=0.3\linewidth]{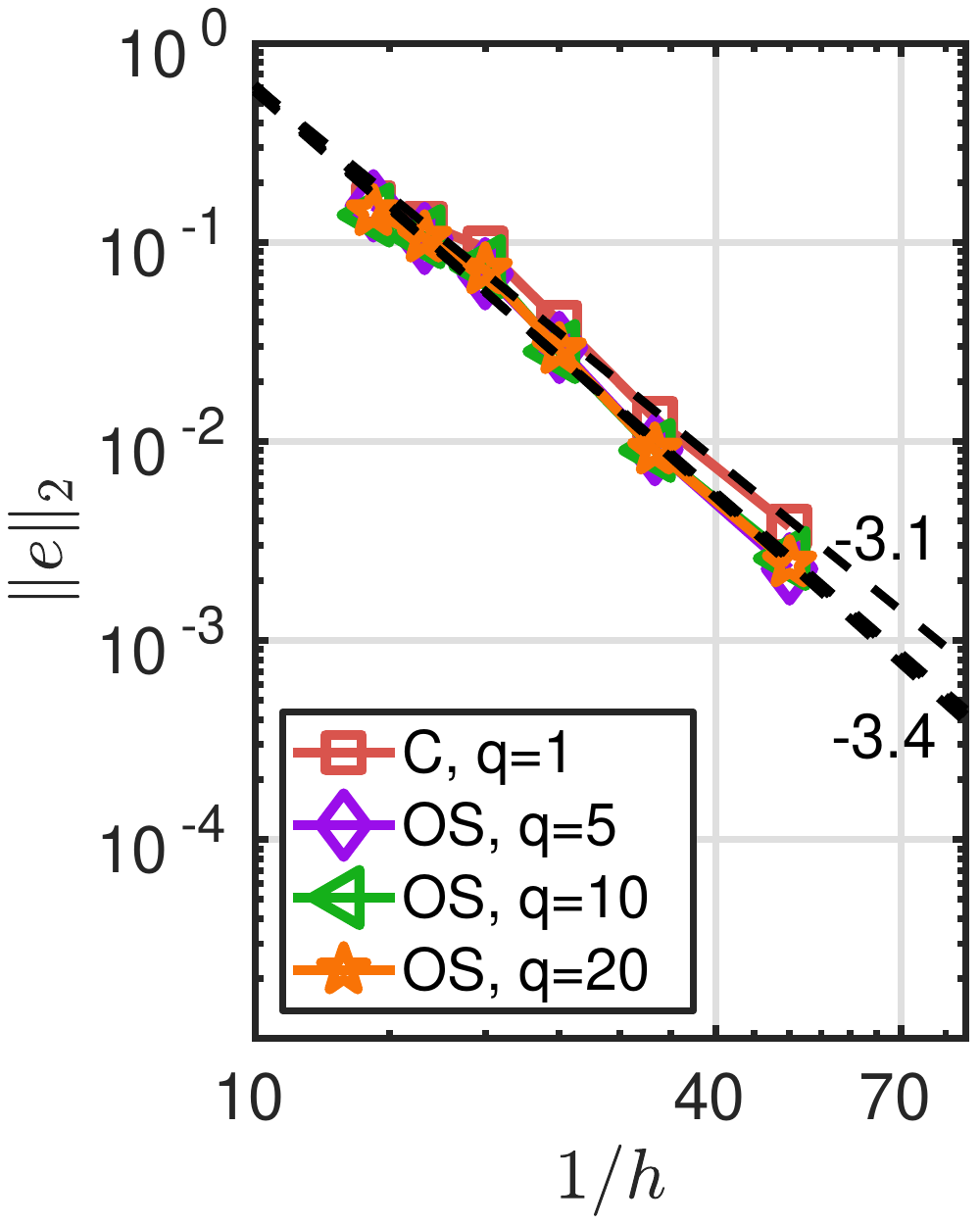} &
    \includegraphics[width=0.3\linewidth]{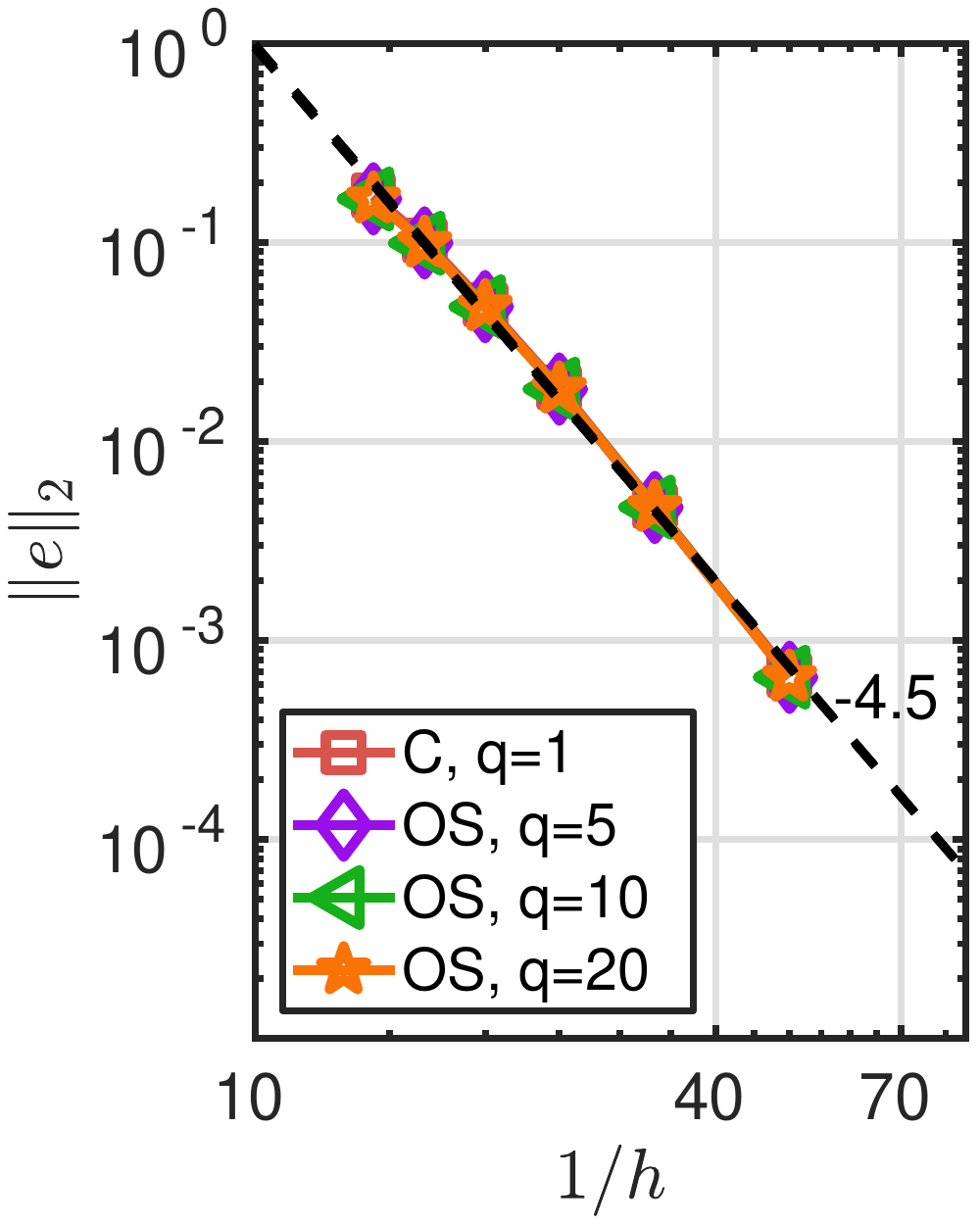}
    \end{tabular}
    \caption{Linear advection of a compactly supported $C^6$ function using a rotational velocity field: convergence as a function of the inverse internodal distance 
    $1/h$. Experiments in each plot correspond to using a different degree of monomial basis $p$ used to construct the RBF-FD trial space. 
    Within each plot, different oversampling parameters $q$ are used to oversample (OS) the semi-discrete system of ODEs. The case $q=1$ 
    represents the collocation RBF-FD method, and other choices of $q$ represent the oversampled RBF-FD method.}
    \label{fig:linearadvection:convergence}
\end{figure}

Now we compare the hyperviscosity stabilized collocation and oversampled RBF-FD methods in terms of the approximation error when the internodal distance $h$ is decreased ($1/h$ increased). 
The simulation is run until $t=5$. In this time the initial condition makes five full rotations around the point $(0,0)$. To 
determine the time step we use $\text{CFL}=0.2$.
Convergence results are collected in Figure \ref{fig:linearadvection:convergence}.
For all tested degrees of the appended monomial bases, we observe that 
the errors are very similar in magnitude and convergence trend, for both collocation and oversampled RBF-FD methods. 
An increase in the oversampling parameter 
does not lead to a better accuracy in this case. This conclusion is in line with the results on solving elliptic PDEs with pure Dirichlet boundary conditions, 
presented in \cite{ToLaHe21}. 
An example of a spatial error distribution after five revolutions of the initial condition is displayed in Figure \ref{fig:linearadvection:solution_spatial_err}.

\subsection{Smooth initial condition: a residual viscosity test}
\begin{figure}[h!]
    \centering
\begin{tabular}{cc}
    \multicolumn{2}{c}{\vspace{0.05cm} \hspace{0.3cm} \textbf{Linear advection: smooth initial condition}} \\
    \hspace{1cm}\textbf{No RV} & \hspace{1cm}\textbf{Using RV} \\
    \includegraphics[width=0.37\linewidth]{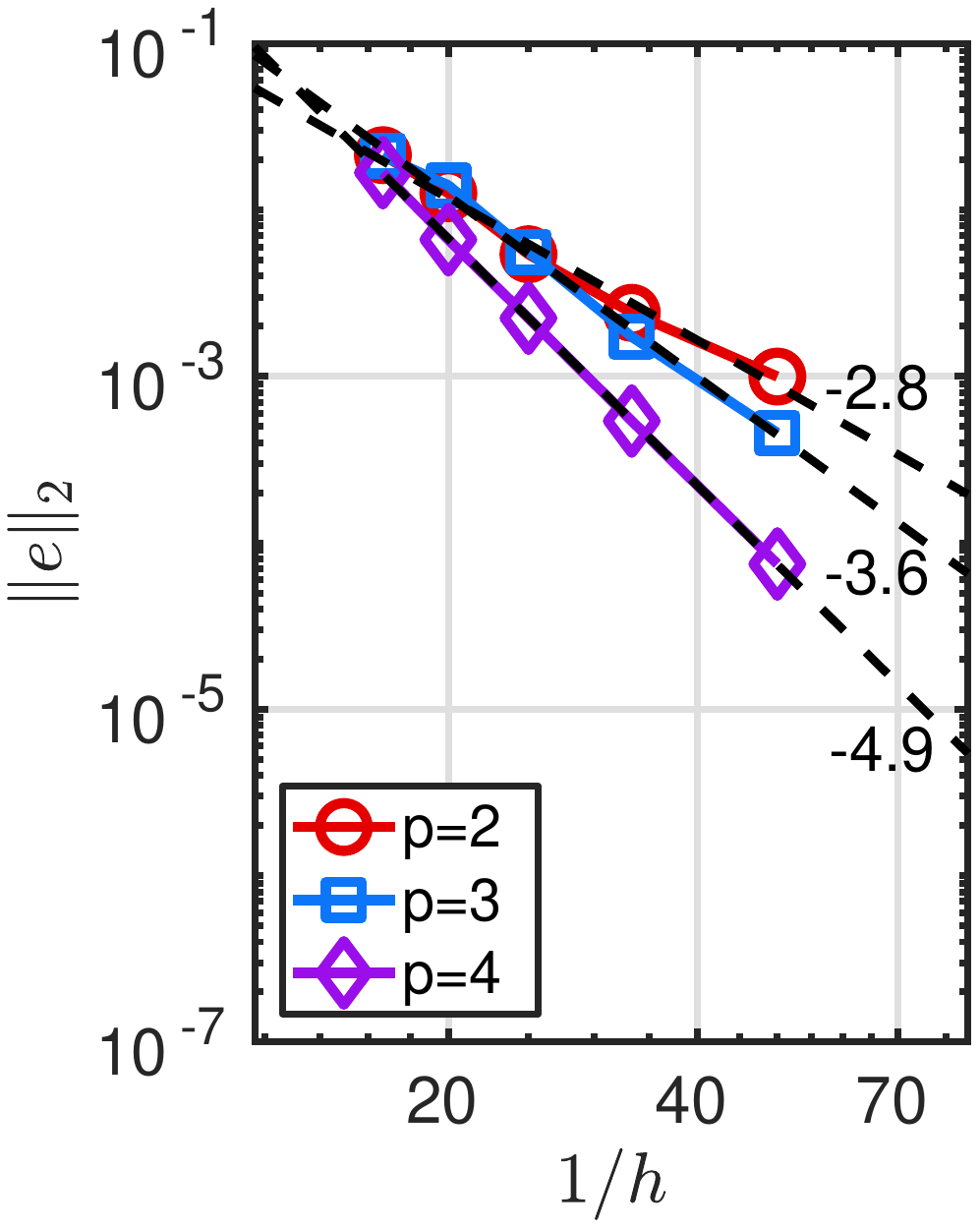} &
    \includegraphics[width=0.37\linewidth]{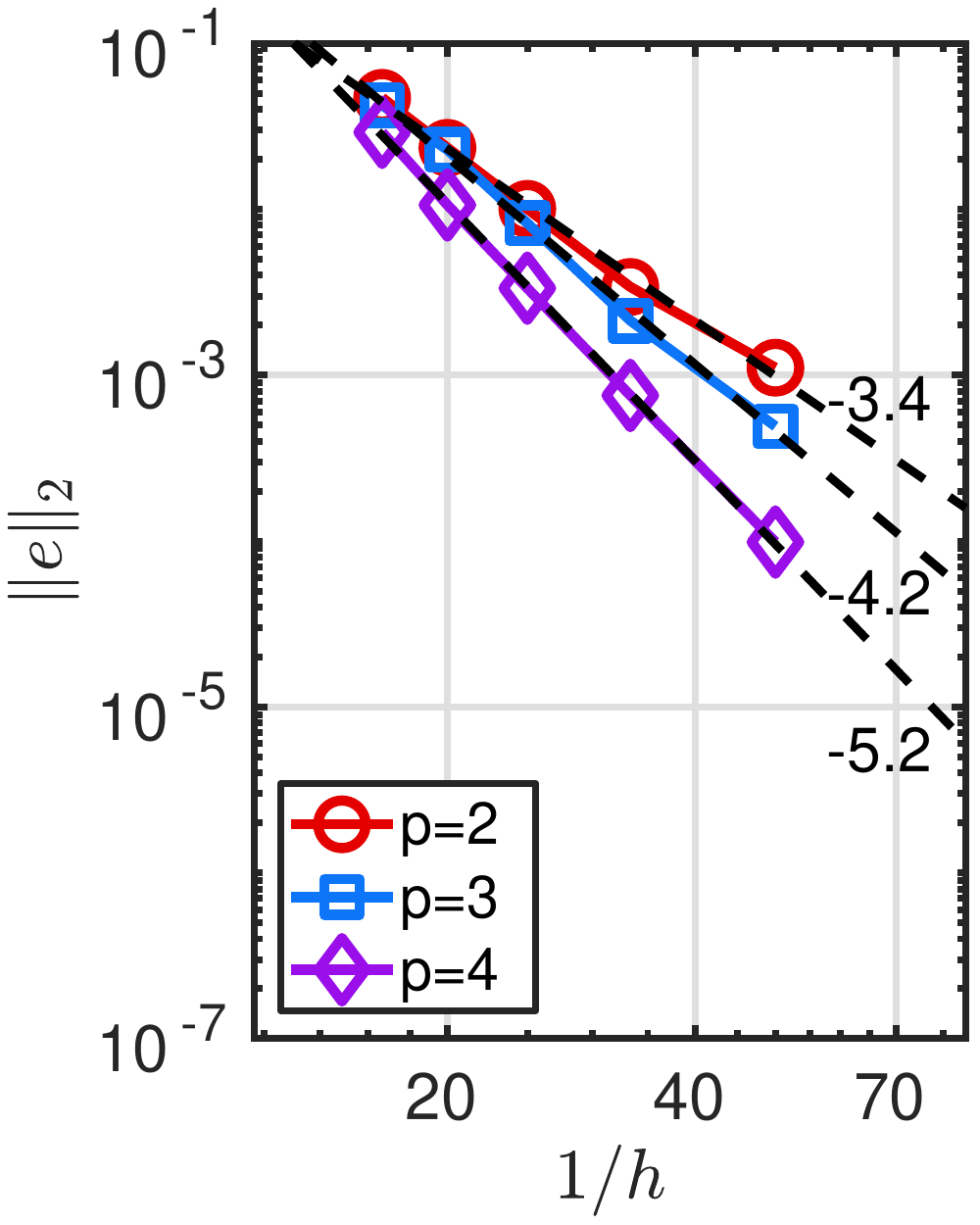}
\end{tabular}
    \caption{Linear advection of a compactly supported $C^6$ function using a rotational velocity field: convergence as a function of the inverse internodal distance $1/h$,
    for an oversampled RBF-FD method ($q=5$). The CFL number is set to $0.2$, the initial condition is rotated until time $t=1$.}
    \label{fig:linearadvection:smoothInitialCondition:RV:convergence}
\end{figure}

We again use \eqref{eq:experiments_linear:velocityfield} as the velocity field and a compactly supported $C^6$ Wendland function \eqref{eq:experiments_linear:initCond_Wendland} 
as the initial condition. The difference compared with the previous experiment is that we now activate the residual viscosity (RV) term $P_2$ 
in our discretization \eqref{eq:discretization:advection_linear_final}. The objective is to check if $P_2$ 
distorts high-order convergence when the solution is smooth. 
We use $C_{\text{RV}} = 1$ as the scaling used in \eqref{eq:resviscosity:epsilon_rv_uw}.
We use $\text{CFL} = 0.5$ and run our simulation until $t=1$. 
The results are collected in 
Figure \ref{fig:linearadvection:smoothInitialCondition:RV:convergence}. 
We observe that the approximation error is slightly larger throughout the refinement of $h$ when RV is used. 
However, the convergence trend when RV is used does not change significantly compared with the case when RV is not used. This is accounted to the fact that the residual is 
small when no discontinuity is present in the numerical solution. Since the coefficient $\varepsilon$ given in \eqref{eq:resviscosity:epsilon} 
that scales $P_2$ in \eqref{eq:resviscosity:P2} 
is residual dependent, the viscosity terms plays a negligible role in the numerical discretization given in \eqref{eq:discretization:advection_linear_final}. An important outcome of this subsection is that the residual viscosity can be 
used in cases when we do not know whether the solution is going to develop a shock, without sacrificing the order of convergence before the shock 
is developed.

\subsection{Discontinuous initial condition: a residual viscosity test}
\begin{figure}[h!]
    \centering
    \begin{tabular}{ccc}
    \includegraphics[width=0.4\linewidth]{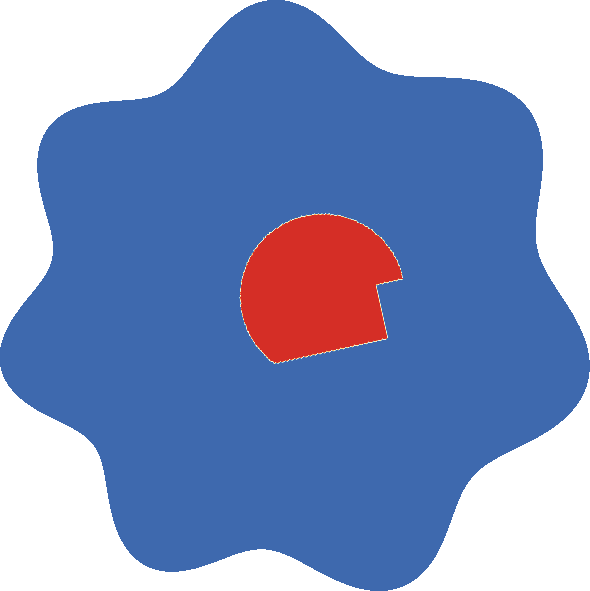} &
    \includegraphics[width=0.4\linewidth]{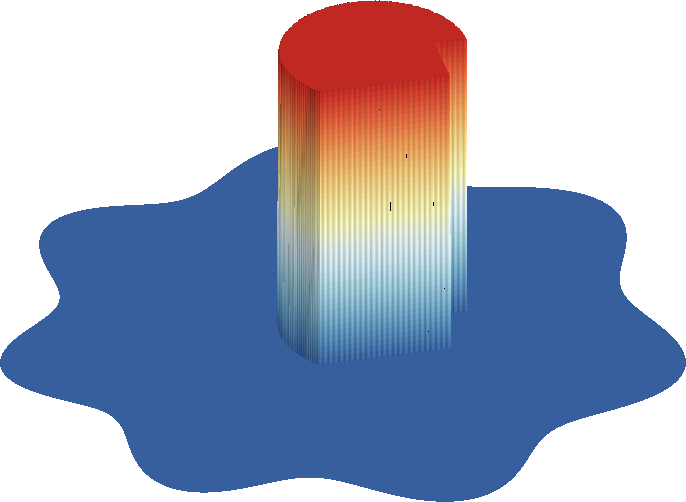} &
    \includegraphics[width=0.08\linewidth]{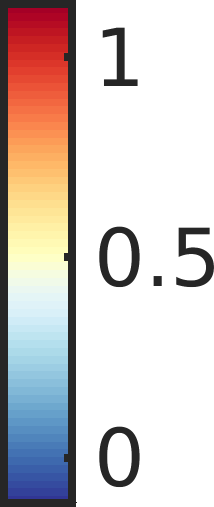}        
    \end{tabular}
    \caption{Linear advection of a discontinuous initial condition: exact solutions at $t=1$, drawn from two perspectives.
    }
    \label{fig:linearadvection:discontInitialCondition_exactsolution}
\end{figure}

\begin{figure}[h!]
    \centering
    \setlength\tabcolsep{1.5pt}
    \begin{tabular}{cccl}
        \multicolumn{4}{c}{\vspace{0.1cm} \textbf{Oversampled RBF-FD method}} \\
        \textbf{No stabilization} & \textbf{First-order viscosity} & \textbf{Residual viscosity} & \\
        \includegraphics[width=0.3\linewidth]{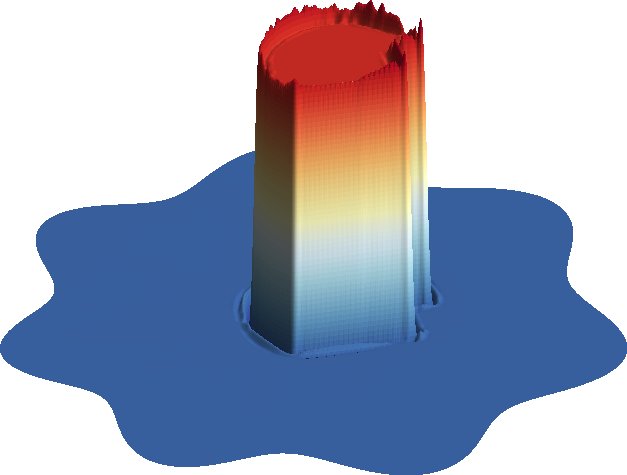} & 
        \includegraphics[width=0.3\linewidth]{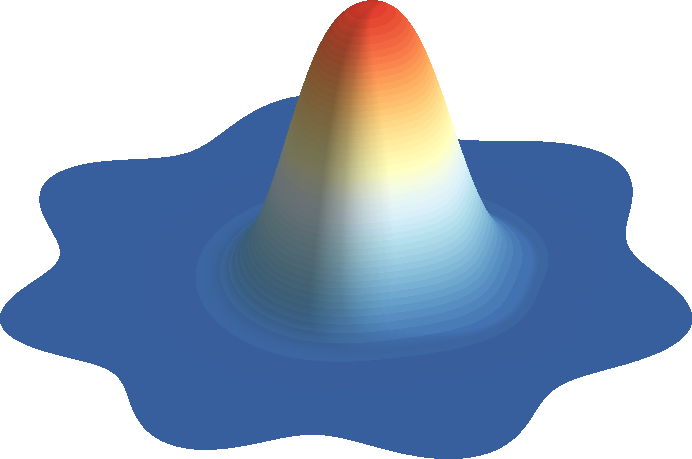} &
        \includegraphics[width=0.3\linewidth]{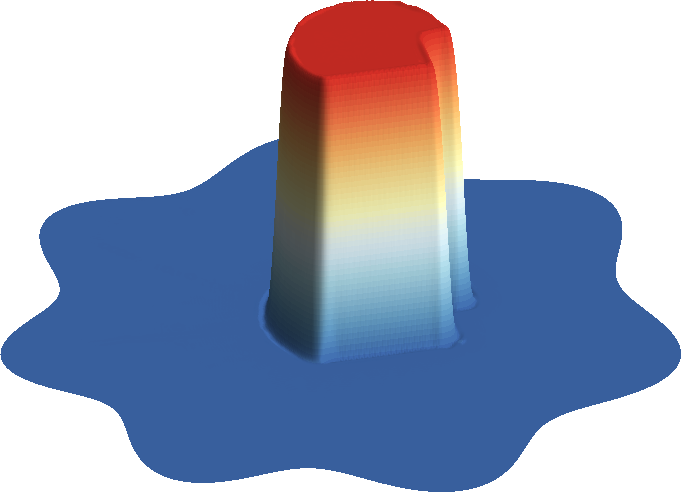} &
        \includegraphics[width=0.05\linewidth]{pics/Pics_linearadvection_discontinuous_solutions/colorbar.png}  \vspace{0.3cm}  \\
        
        \multicolumn{4}{c}{\vspace{0.1cm}\textbf{Collocation RBF-FD method}} \\
        \textbf{No stabilization} & \textbf{First-order viscosity} & \textbf{Residual viscosity} & \\
        \includegraphics[width=0.3\linewidth]{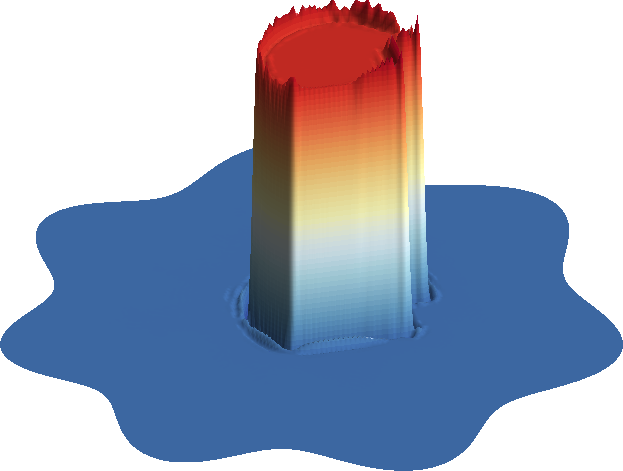} & 
        \includegraphics[width=0.3\linewidth]{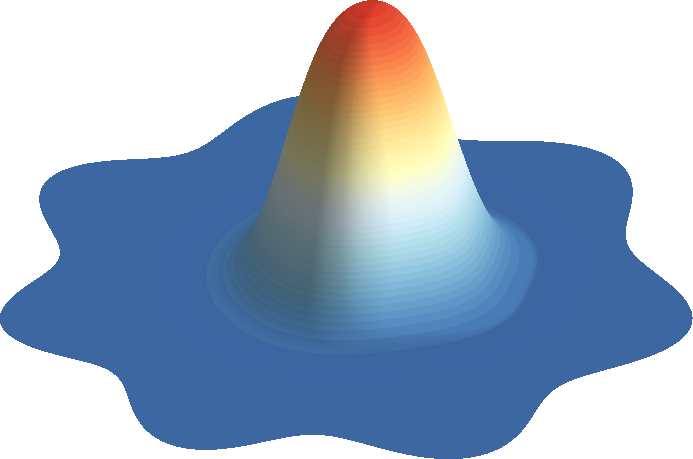} &
        \includegraphics[width=0.3\linewidth]{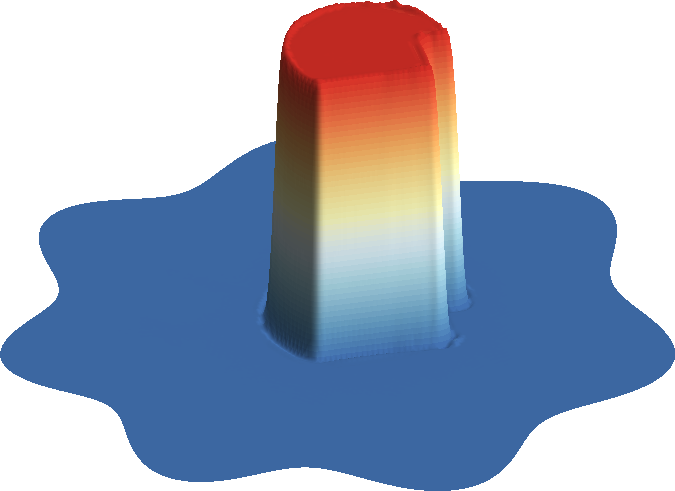} &
        \includegraphics[width=0.05\linewidth]{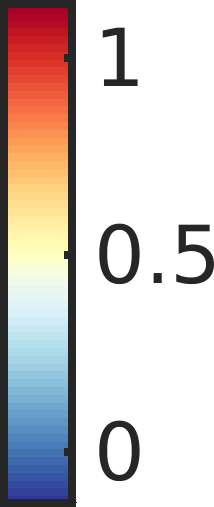}              
    \end{tabular}
    \caption{Linear advection of a discontinuous initial condition when different stabilization techniques are used, side view. All solutions are displayed after one rotation at $t=1$. 
    Parameters used were $h=0.01$, $p=4$, $q=5$. The CFL condition when first-order viscosity was set to $0.1$, and to $0.5$ in cases when residual viscosity method was used and when 
    no stabilization was used.
    }
    \label{fig:linearadvection:discontInitialCondition}
\end{figure}

Now we again consider a linear advection problem with \eqref{eq:experiments_linear:velocityfield} as the velocity field, but 
use a discontinuous initial condition; a cylinder with radius $0.3$, cut from two sides, given in Figure \ref{fig:linearadvection:discontInitialCondition_exactsolution}. 
We let the initial condition rotate around its axis one time (until $t=1$) when $h=0.01$, $p=4$, and compare three scenarios:
\begin{enumerate}
    \item Residual viscosity is not active: in \eqref{eq:resviscosity:P2} we set $\varepsilon = 0$.
\item  Only the first order viscosity is active: instead of following the definition of $\varepsilon$ in \eqref{eq:resviscosity:epsilon} we set $\varepsilon = \varepsilon_{\text{UW}}$ throughout $\Omega$.
\item  Residual viscosity is active: $\varepsilon$ is taken as defined in \eqref{eq:resviscosity:epsilon}. 
\end{enumerate}
Note that the hyperviscosity term $P_1$ is essential to stabilize the discretization in time  
and is active in all three scenarios. We use $\text{CFL}=0.1$ when using first-order viscosity and $\text{CFL}=0.5$ when using 
RV and when not stabilizing discontinuities. We use $C_{\text{RV}} = 1$ as the scaling used in \eqref{eq:resviscosity:epsilon_rv_uw}.
Results for the collocation RBF-FD method and an oversampled RBF-FD method are given in Figure \ref{fig:linearadvection:discontInitialCondition} (side view). 

From the two figures we observe that the numerical scheme without stabilization induces oscillations (Gibbs phenomenon) around the discontinuity. This is expected. 
When the numerical scheme is stabilized by the first-order viscosity term, we observe that the initial condition is highly smeared, 
to the extent where it is difficult 
to recognize the characteristic shape of the initial condition. The main observation is that after the oversampled RBF-FD method is stabilized using the residual viscosity term, the oscillations are significantly damped, while the initial condition 
kept its characteristic shape throughout the simulation.
We observe that the RV solution is not significantly different 
compared with the RV solution when an oversampled RBF-FD method is used. It is however possible to observe 
a few more oscillations on the outer edge of the slotted cylinder, when the collocation RBF-FD method is used.

\begin{figure}[h!]
    \centering
    \setlength\tabcolsep{1.5pt}    
    \begin{tabular}{ccc}
        \multicolumn{3}{c}{\vspace{0.1cm}\textbf{Oversampled RBF-FD method}} \\
        \hspace{-1.1cm}\textbf{Residual} & \hspace{-1.1cm}\textbf{RV coefficient} & \hspace{-1.1cm}\textbf{RV coefficient ($\mathbf{\log_{10}}$)} \\
        \includegraphics[width=0.32\linewidth]{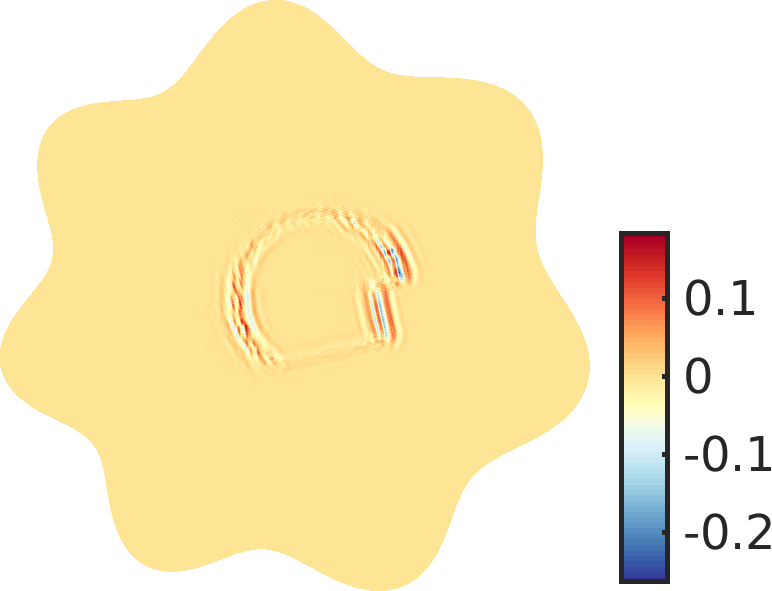} & 
        \includegraphics[width=0.32\linewidth]{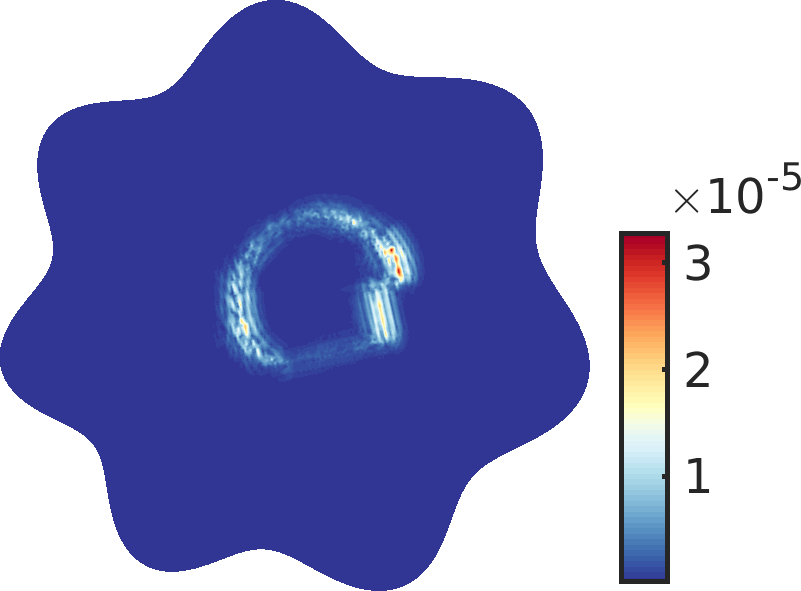} &
        \includegraphics[width=0.32\linewidth]{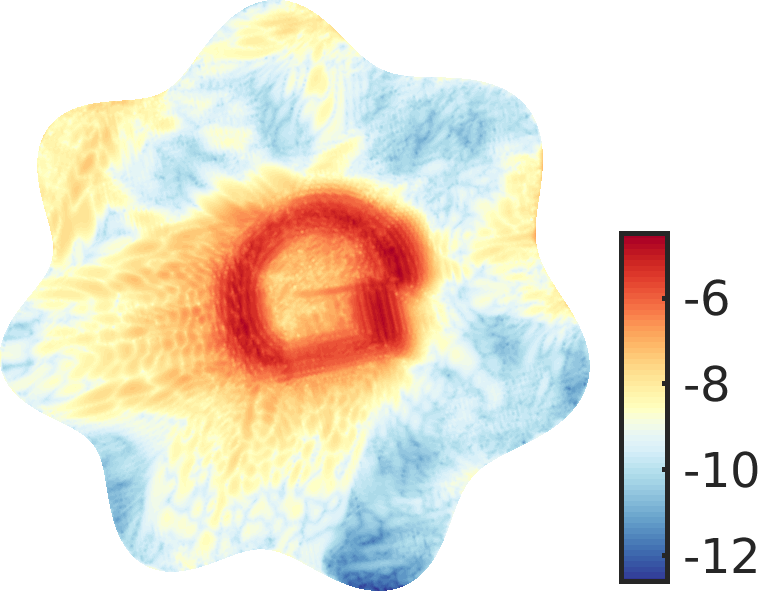} \vspace{0.3cm} \\
        \multicolumn{3}{c}{\vspace{0.1cm}\textbf{Collocation RBF-FD method}} \\
        \hspace{-1.1cm}\textbf{Residual} & \hspace{-1.1cm}\textbf{RV coefficient} & \hspace{-1.1cm}\textbf{RV coefficient ($\mathbf{\log_{10}}$)} \\
        \includegraphics[width=0.32\linewidth]{pics/Pics_linearadvection_discontinuous_residual_coefficient/rv_numerical_residual.png} & 
        \includegraphics[width=0.32\linewidth]{pics/Pics_linearadvection_discontinuous_residual_coefficient/rv_numerical_epsilon.png} &
        \includegraphics[width=0.32\linewidth]{pics/Pics_linearadvection_discontinuous_residual_coefficient/rv_numerical_epsilon_log10.png} 
        
    \end{tabular}
    \caption{Linear advection of a discontinuous initial condition when residual viscosity is used. The pictures show the residual and the residual-viscosity (RV) coefficients in 
    linear and logarithmic scale.}
    \label{fig:linearadvection:rv_coefficient}
\end{figure}

The effectiveness of the residual viscosity method when an oversampled RBF-FD method is used, is displayed in Figure \ref{fig:linearadvection:rv_coefficient}, 
where we show the spatial distribution of the residual \eqref{eq:resviscosity:residual} and the coefficient $\varepsilon$ given in \eqref{eq:resviscosity:epsilon}. 
We observe that 
the residual is giving the information about the position of the discontinuity (large oscillations) present in the solution. This makes the residual viscosity coefficient $\varepsilon$ 
(defined in \eqref{eq:resviscosity:epsilon} and \eqref{eq:resviscosity:epsilon_rv_uw}) large in the region of the shock, but small away from the shock. 
The results for the oversampled and the collocation RBF-FD methods are -- in the ''eyeball norm'' -- identical.
We conclude that for the cases considered in this subsection, the residual viscosity approach detects discontinuities in an accurate way, for both RBF-FD methods that we tested. 
In this paper we are mainly interested in the performance of the oversampled RBF-FD method, which we now focus on in the results below.

\begin{figure}[h!]
    \centering
    \setlength\tabcolsep{1.5pt}    
    \begin{tabular}{ccc}
        \multicolumn{3}{c}{\vspace{0.1cm}\textbf{Linear advection: discontinuous initial condition}} \\
        \hspace{0.8cm}\textbf{No stabilization} & \hspace{0.9cm}\textbf{First-order viscosity} & \hspace{0.7cm}\textbf{Residual viscosity} \\
        \includegraphics[width=0.31\linewidth]{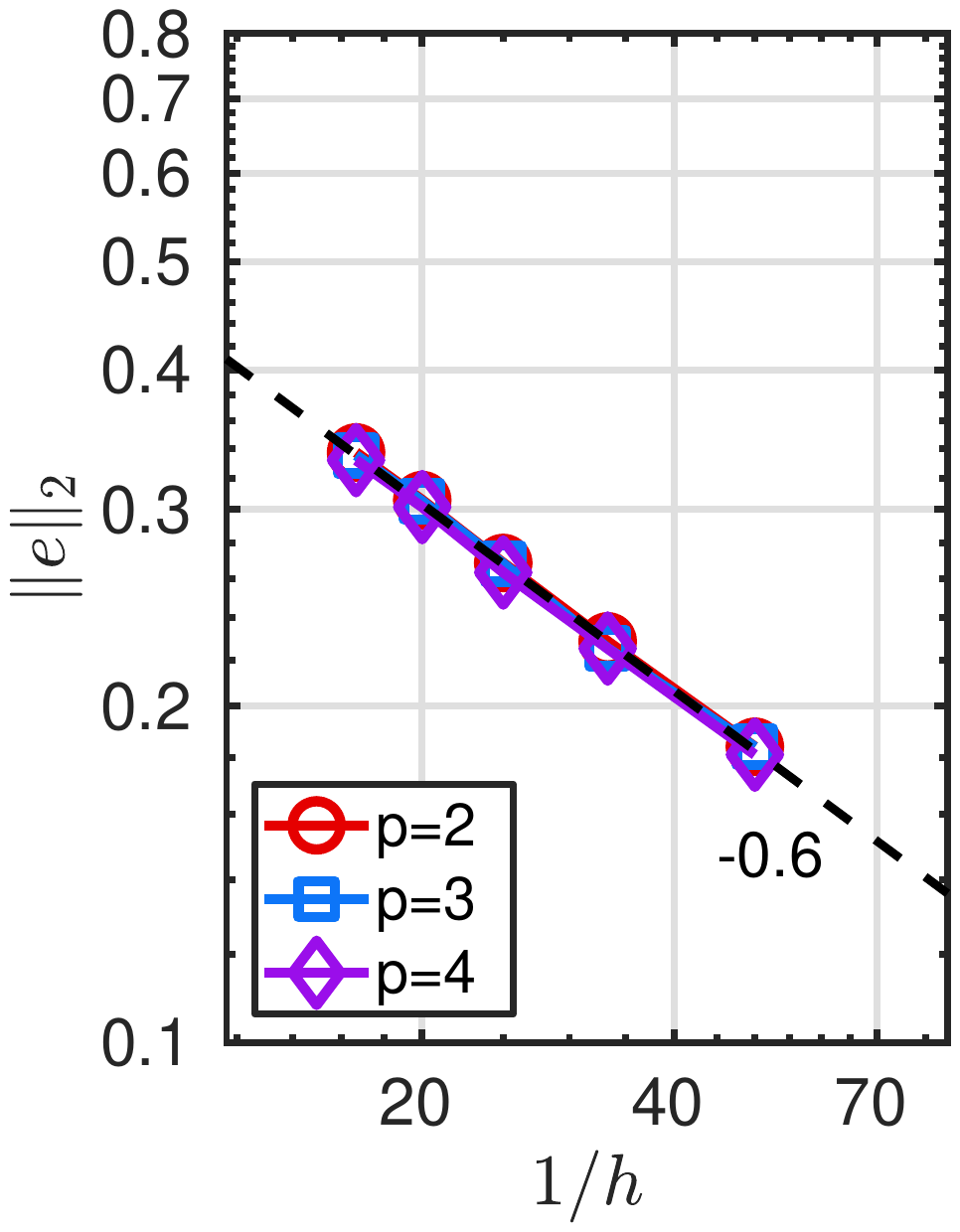} & 
        \includegraphics[width=0.31\linewidth]{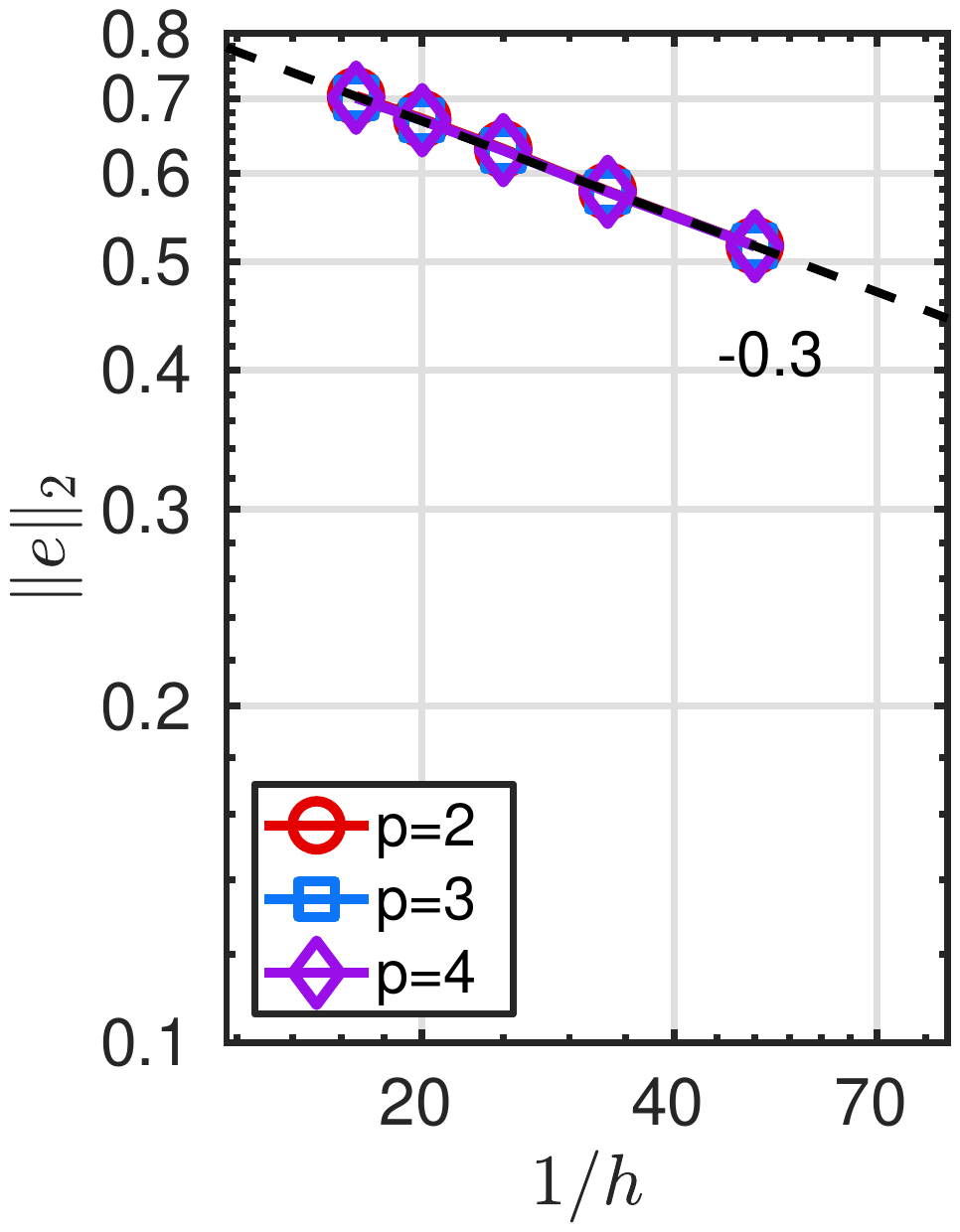} & 
        \includegraphics[width=0.31\linewidth]{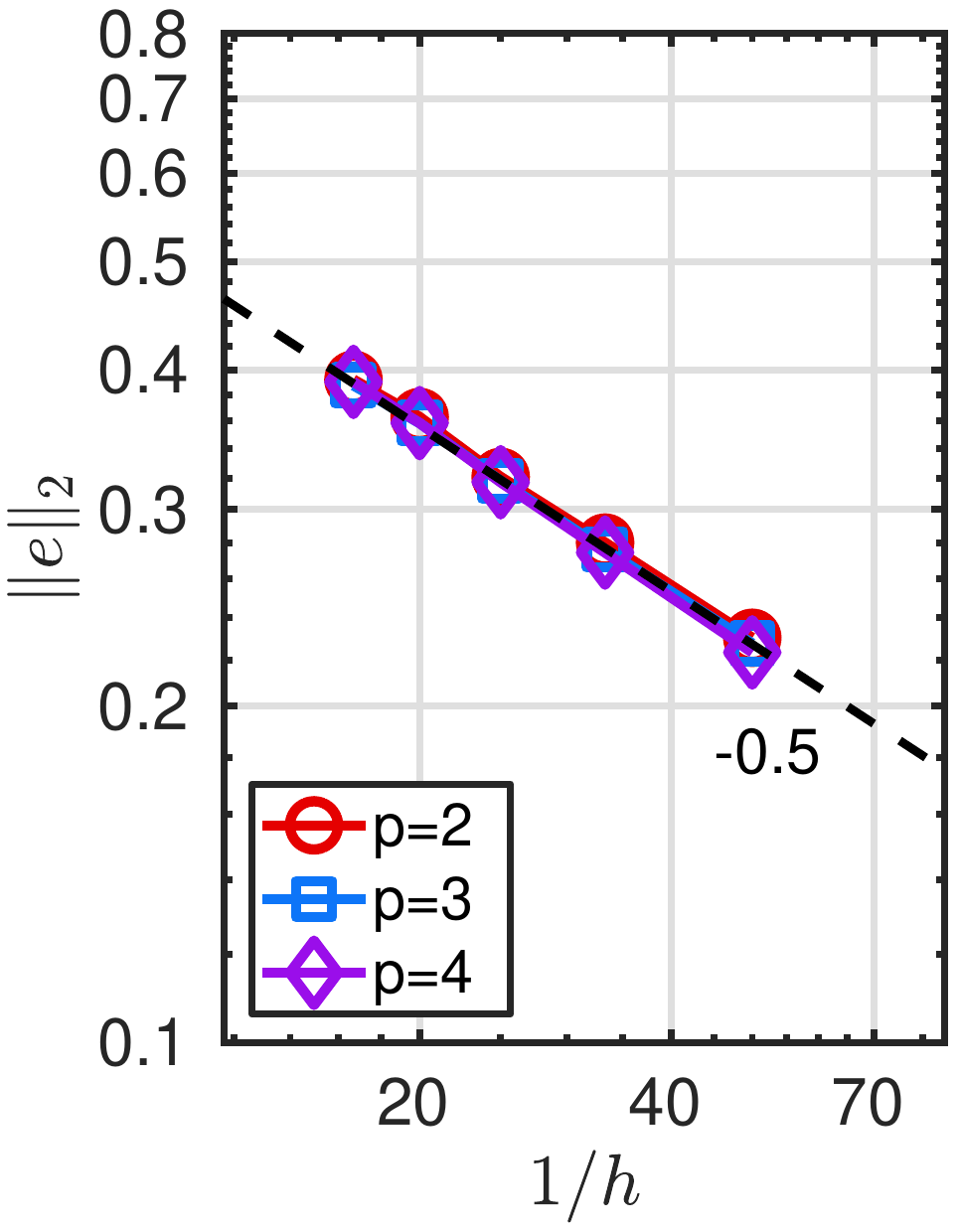} 
    \end{tabular}
    \caption{Linear advection of a discontinuous initial condition, convergence study of the approximation error under node refinement, for different degrees $p$ of the monomial basis used 
    to construct the interpolation problem over a stencil. The three plots show cases when numerical scheme includes: (i) no shock stabilization, (ii) first-order viscosity stabilization, (iii) 
    residual-viscosity stabilization.}
    \label{fig:linearadvection:discontinuous:convergence}
\end{figure}

In Figure \ref{fig:linearadvection:discontinuous:convergence} we provide convergence results for the oversampled RBF-FD 
method, for different choices of $p$, 
when the numerical scheme includes: (i) no stabilization, (ii) first-order viscosity, (iii) residual-viscosity. 
The CFL numbers used for the simulation were $0.2$ for 
cases (i) and (ii), and $0.1$ for case (iii). The observed convergence rates in all cases are what we expect, since the solution is discontinuous and the error is measured in $2$-norm.  
The approximation error is smallest in the case where no shock stabilization is added. This can be accounted to an observation made in Figure \ref{fig:linearadvection:discontInitialCondition}, 
where oscillations when no shock stabilization was added were not as severe as one expected. However, we do expect the oscillations to become significantly 
larger when the flux $\bm F(u)$ becomes nonlinear with respect to $u$ (studied in later sections). Additional experiments revealed that the approximation error in $\infty$-norm does not converge, 
which is expected.

\section{Numerical study II: Burger's equation}
\label{sec:experiments:burgers}
\begin{figure}[h!]
    \centering
    \setlength\tabcolsep{1.5pt}    
    \begin{tabular}{cccl}
        \multicolumn{4}{c}{\hspace{-1.3cm}\vspace{0.1cm}\textbf{Burger's equation: exact solution}} \\
        \hspace{0.1cm}\textbf{$\mathbf{t=0}$} & \hspace{0.1cm}\textbf{$\mathbf{t=0.5}$} & \hspace{0.1cm}\textbf{$\mathbf{t=0.5}$ (side)} & \\
        \includegraphics[width=0.28\linewidth]{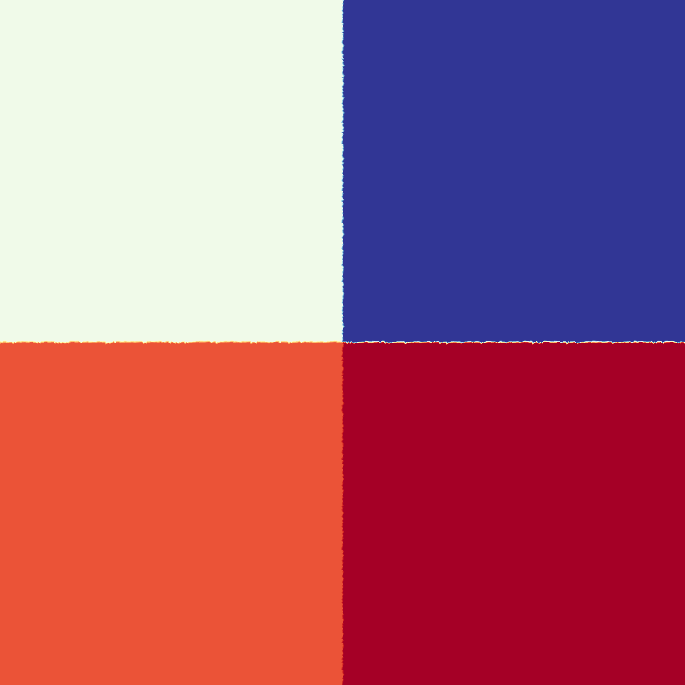} & 
        \includegraphics[width=0.28\linewidth]{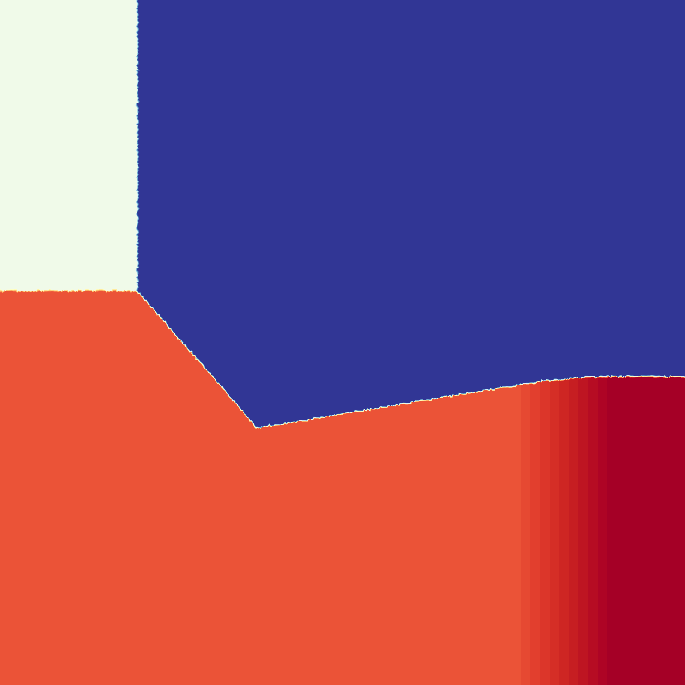} & 
        \includegraphics[width=0.28\linewidth]{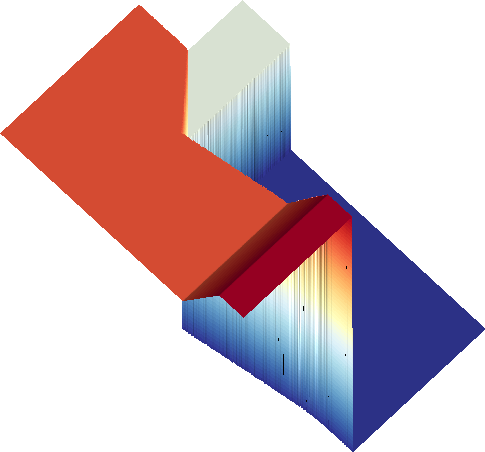} & 
        \includegraphics[width=0.07\linewidth]{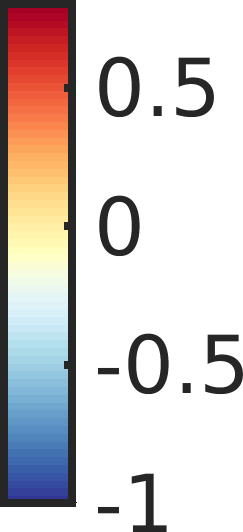}
    \end{tabular}
    \caption{Exact solution to the considered Burger's equation, at time $t=0$ (initial condition, and at time $t=0.5$ from two visual perspectives.}
    \label{fig:burgers:exactsolutions}
\end{figure}

In this section we solve Burger's equation in two dimensions. The flux is given by $\bm F(u) = (\frac{u^2}{2}, \frac{u^2}{2})$, 
and the velocity field by $\bm F'(u) = (u, u)$. 
The exact solution to the considered problem is:
  \begin{equation}
   u(y,t) = \left\{\!\!
      \begin{aligned}
        & \begin{aligned}
          -&0.2 \\
          &0.5
        \end{aligned}
        &\quad &\mbox{if } y_1<\tfrac12-\tfrac{3t}{5} \mbox{ and }
        \left\{
        \begin{aligned}
          &y_2 > \tfrac12 + \tfrac{3t}{20}, \\
          &\mbox{otherwise},
        \end{aligned}
        \right . \\
        & \begin{aligned}
          -&1 \\
          &0.5
        \end{aligned}
        &\quad &\mbox{if } \tfrac12-\tfrac{3t}{5}<y_1<\tfrac12-\tfrac{t}{4} \mbox{ and }
        \left\{
        \begin{aligned}
          &y_2 > -\tfrac{8y_1}{7} + \tfrac{15}{14} - \tfrac{15t}{28}, \\
          &\mbox{otherwise},
        \end{aligned}
        \right . \\
        & \begin{aligned}
          -&1 \\
          &0.5
        \end{aligned}
        &\quad &\mbox{if } \tfrac12-\tfrac{t}{4}<y_1<\tfrac12+\tfrac{t}{2} \mbox{ and }
        \left\{
        \begin{aligned}
          &y_2 > \tfrac{y_1}{6} + \tfrac{5}{12} - \tfrac{5t}{24}, \\
          &\mbox{otherwise},
        \end{aligned}
        \right . \\
        & \begin{aligned}
          -&1 \\
          &\tfrac{2y_1-1}{2t}
        \end{aligned}
        &\quad &\mbox{if } \tfrac12+\tfrac{t}{2}<y_1<\tfrac12+\tfrac{4t}{5} \mbox{ and }
        \left\{
        \begin{aligned}
          &y_2 > y_1 -\tfrac{5}{18t} \left(y_1+t-\tfrac12\right)^2, \\
          &\mbox{otherwise},
        \end{aligned}
        \right . \\
        & \begin{aligned}
          -&1 \\
          &0.8
        \end{aligned}
        &\quad &\mbox{if } y_1>\tfrac12+\tfrac{4t}{5} \mbox{ and }
        \left\{
        \begin{aligned}
          &y_2 > \tfrac12 - \tfrac{t}{10}, \\
          &\mbox{otherwise}.
        \end{aligned}
        \right . \\
      \end{aligned}
    \right.\hspace{-5mm}
  \end{equation}
The initial condition simplifies to:
\begin{equation}
    u(y,0) = \left\{
    \begin{aligned}
      -&0.2, &\quad \mbox{if } y_1<0.5 \mbox{ and } y_2>0.5, \\
      -&1, &\quad \mbox{if } y_1>0.5 \mbox{ and } y_2>0.5, \\
      &0.5, &\quad \mbox{if } y_1<0.5 \mbox{ and } y_2<0.5, \\
      &0.8, &\quad \mbox{if } y_1>0.5 \mbox{ and } y_2<0.5. \\
    \end{aligned}
    \right.
  \end{equation}
The initial condition and the exact solution to the considered problem are visualized in Figure \ref{fig:burgers:exactsolutions}.
A variable time step is chosen as:
\begin{equation}
  \label{eq:burgers:timestep}
\Delta t = \text{CFL}\, \min_{x_i\in X} \left [\frac{h_{\text{loc}}(x_i)}{\sqrt{(F_1'(u(x_i,t)))^2 + (F_2'(u(x_i,t)))^2}} \right ],
\end{equation}
where $h_{\text{loc}}(x_i)$ is defined in the scope of \eqref{eq:resviscosity:epsilon_rv_uw}.

We solve this problem on a unit square $\Omega = [0,1] \times [0,1]$ until $t=0.5$ and use Dirichlet boundary conditions with the data obtained from 
the exact solution. Our objective is to make observations on 
how the oversampled RBF-FD method behaves when the nonlinear flux $\bm F (u)$ is inducing shocks in the solution. We focus on cases 
when we stabilize the Gibbs phenomenon through $P_2$ defined in \eqref{eq:resviscosity:P2}) using: 
(i) first order viscosity (we override \eqref{eq:resviscosity:epsilon} by $\varepsilon = \varepsilon_{\text{UW}}$), 
(ii) residual viscosity ($\varepsilon$ as defined in \eqref{eq:resviscosity:epsilon}). The hyperviscosity term \eqref{eq:discretization:hypervi} 
is always added to the numerical scheme for the purposes of stabilization in time. 
Throughout this section we use an oversampling parameter $q=5$ and the CFL number 
$0.2$. The nodes in $X$ point set are generated using the algorithm presented in \cite{FBF15_nodes}. 

\begin{figure}[h!]
    \centering
    \setlength\tabcolsep{1.5pt}    
    \begin{tabular}{cccl}
        \multicolumn{4}{c}{\hspace{-1.6cm} \vspace{0.1cm}\textbf{Burger's equation: numerical solutions}} \\
        \hspace{0.1cm}\textbf{First-order viscosity} & \hspace{0.1cm}\textbf{RV ($\mathbf{c_{\text{RV}}=1}$)} & \hspace{0.1cm}\textbf{RV ($\mathbf{c_{\text{RV}}=4}$) } & \\
        \includegraphics[width=0.28\linewidth]{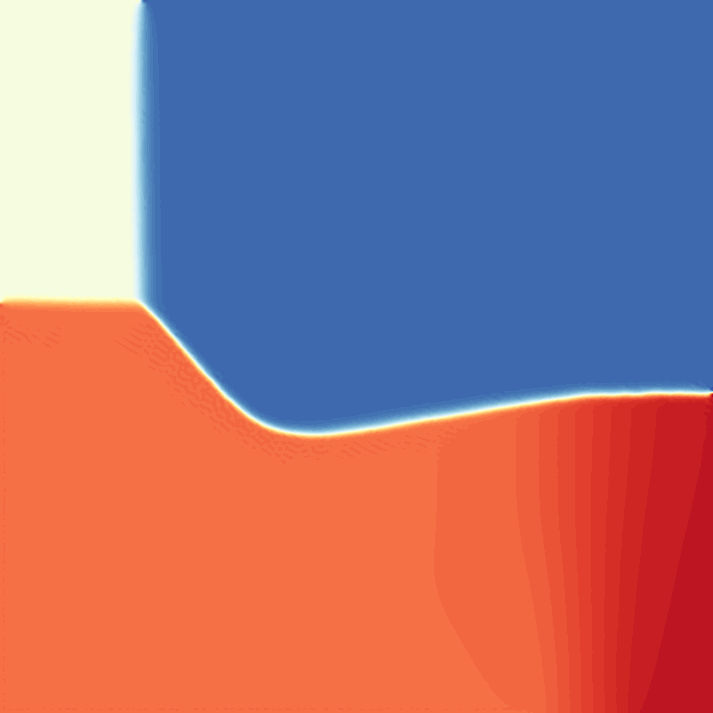} & 
        \includegraphics[width=0.28\linewidth]{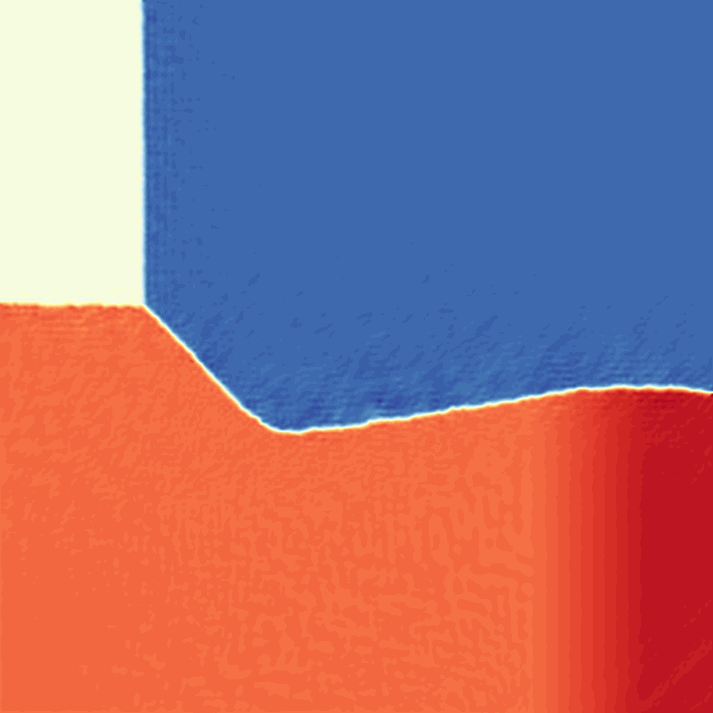} & 
        \includegraphics[width=0.28\linewidth]{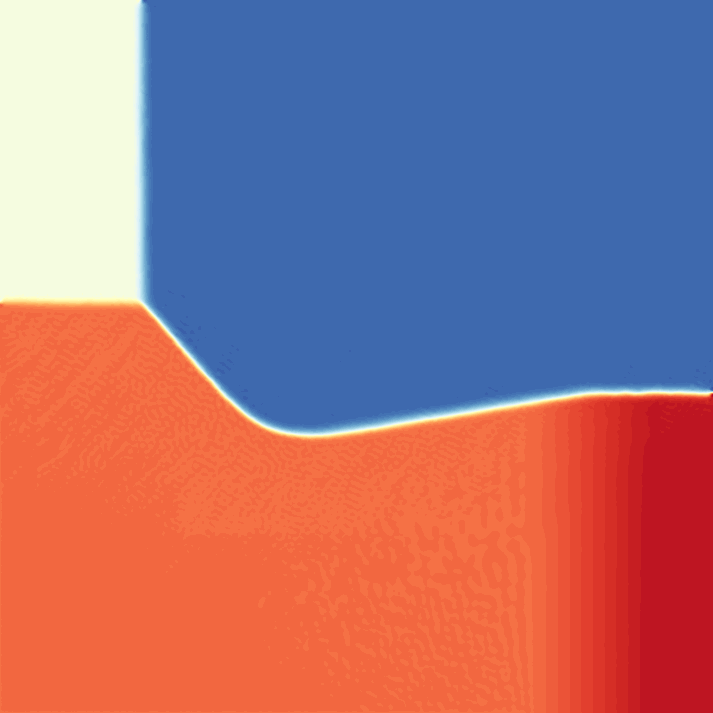} & 
        \includegraphics[width=0.07\linewidth]{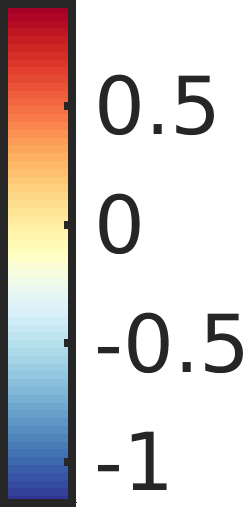}
    \end{tabular}
    \caption{Numerical solutions to the Burger's equations, when the oversampled RBF-FD discretization is shock-stabilized using 
    first-order viscosity, and residual viscosity for two choices of constants $C_{\text{RV}}$ in \eqref{eq:resviscosity:epsilon_rv_uw}. The internodal 
    distance is set to $h=0.005$ (corresponds to $N=39504$ unknowns), polynomial degree $p=3$ is used to construct the stencil-based interpolation problems.}
    \label{fig:burgers:solutions_normalview}
\end{figure}

\begin{figure}[h!]
    \centering
    \setlength\tabcolsep{1.5pt}    
    \begin{tabular}{cccl}
        \multicolumn{4}{c}{\hspace{-1.6cm}\vspace{0.1cm}\textbf{Burger's equation: numerical solutions (side view)}} \\
        \hspace{0.1cm}\textbf{First-order viscosity} & \hspace{0.1cm}\textbf{RV ($\mathbf{c_{\text{RV}}=1}$)} & \hspace{0.1cm}\textbf{RV ($\mathbf{c_{\text{RV}}=4}$) } & \\
        \includegraphics[width=0.28\linewidth]{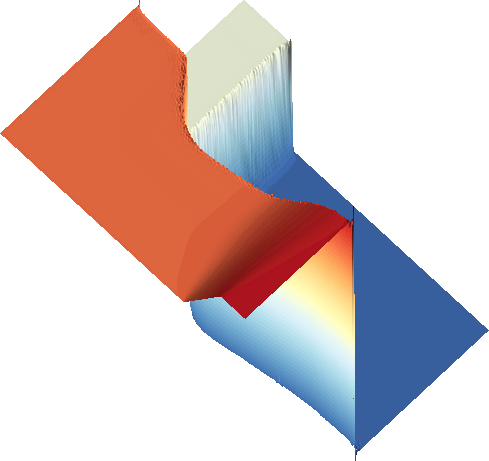} & 
        \includegraphics[width=0.28\linewidth]{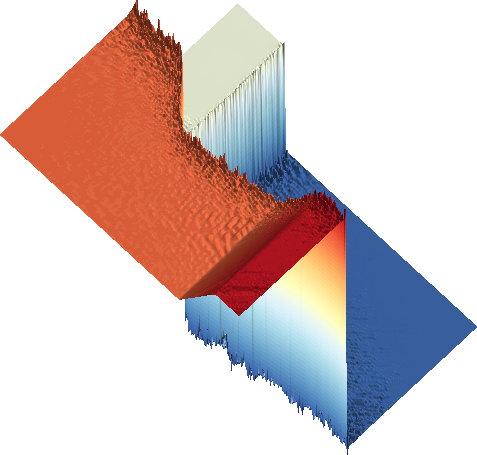} & 
        \includegraphics[width=0.28\linewidth]{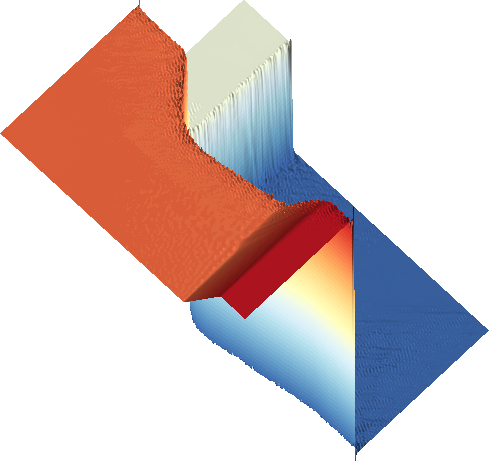} & 
        \includegraphics[width=0.07\linewidth]{pics/Pics_Burgers/numerical_solution_colorbar.png}
    \end{tabular}
    \caption{Numerical solutions to the Burger's equations from a rotated perspective, when the oversampled RBF-FD discretization is shock-stabilized using 
    first-order viscosity, and residual viscosity for two choices of constants $c_{\text{RV}}$ in \eqref{eq:resviscosity:epsilon_rv_uw}. 
    The internodal 
    distance is set to $h=0.005$ (corresponds to $N=39504$ unknowns), polynomial degree $p=3$ is used to construct the stencil-based interpolation problems.}
    \label{fig:burgers:solutions}
\end{figure}

First, we examine the numerical solution from a visual perspective in the three cases: when the RBF-FD method is not shock-stabilized, 
when the method is stabilized with the first-order viscosity and when the 
method is stabilized using the residual viscosity method. 
The parameters that we use are $h=0.005$ ($N=40000$), $p=3$ and $q=5$. 
The numerical solution is a subject to rapid exponential growth when the Gibbs phenomenon is not stabilized with any of the approaches, 
and we therefore can not provide any visual results.
The solution when the numerical scheme is stabilized using first order viscosity and residual viscosity, is given in 
Figure \ref{fig:burgers:solutions_normalview} and in Figure \ref{fig:burgers:solutions}.

\begin{figure}[h!]
    \centering
    \setlength\tabcolsep{1.5pt}    
    \begin{tabular}{cccl}
        \multicolumn{4}{c}{\hspace{-1.6cm} \vspace{0.1cm}\textbf{Burger's equation: viscosity coefficients}} \\
        \hspace{0.1cm}\textbf{First-order viscosity} & \hspace{0.1cm}\textbf{RV ($\mathbf{c_{\text{RV}}=1}$)} & \hspace{0.1cm}\textbf{RV ($\mathbf{c_{\text{RV}}=4}$) } & \\
        \includegraphics[width=0.28\linewidth]{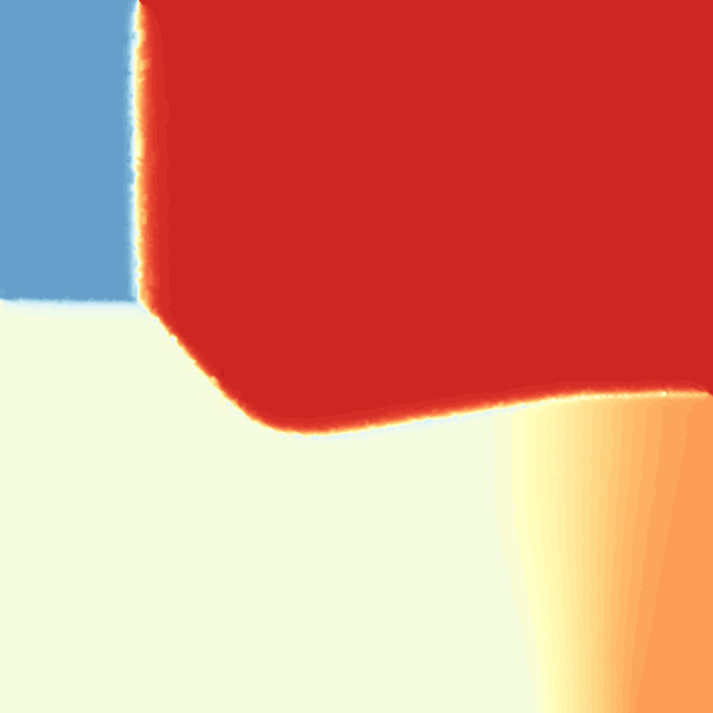} & 
        \includegraphics[width=0.28\linewidth]{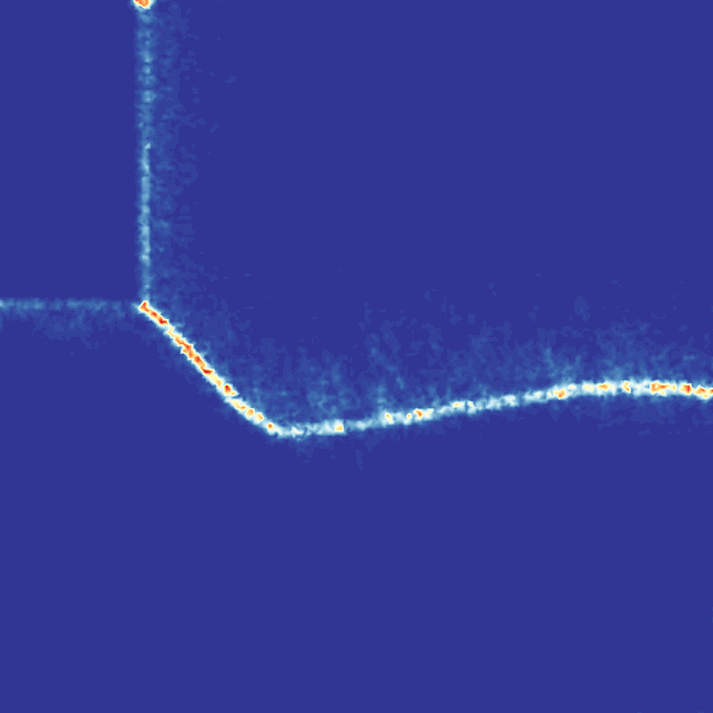} & 
        \includegraphics[width=0.28\linewidth]{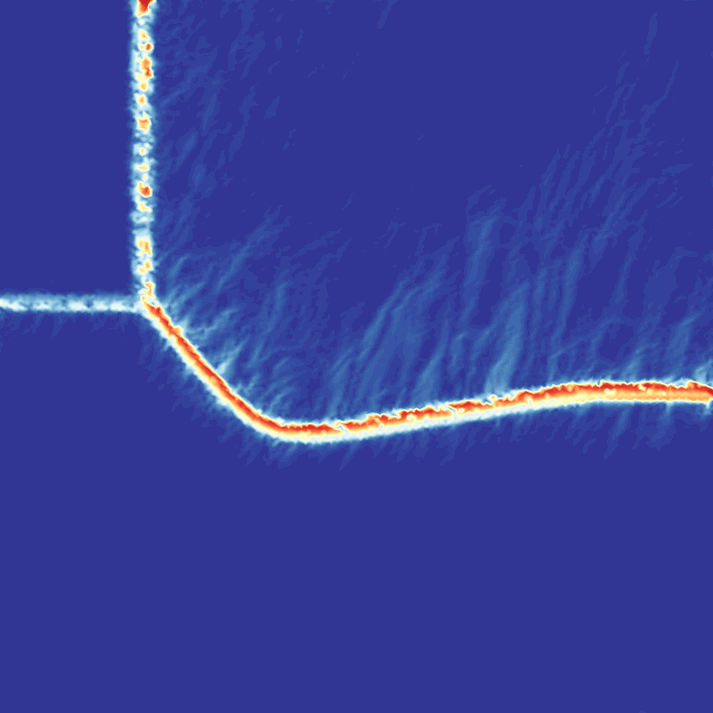} & 
        \includegraphics[width=0.09\linewidth]{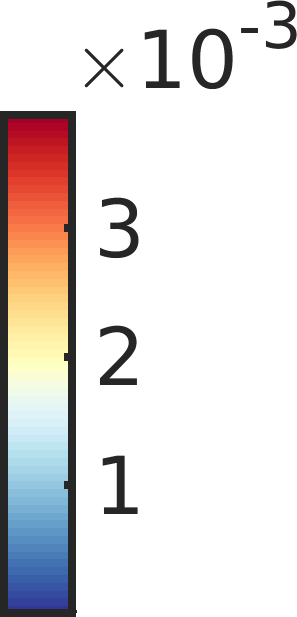}
    \end{tabular}
    \caption{Spatial variation of the viscosity coefficients, when the oversampled RBF-FD discretization is shock-stabilized using 
    first-order viscosity, and residual viscosity for two choices of constants $c_{\text{RV}}$ in \eqref{eq:resviscosity:epsilon_rv_uw}. 
    The internodal 
    distance is set to $h=0.005$ (corresponds to $N=39504$ unknowns), polynomial degree $p=3$ is used to construct the stencil-based interpolation problems.}
    \label{fig:burgers:ep_coefficients}
\end{figure}

In all cases where we attempted to stabilize the solution, the numerical solution was stable at the end. 
The smoothest solution is expectedly the one where the first-order viscosity is active throughout $\Omega$. 
When RV is used, the details around the discontinuities are preserved in a better way compared with the first-order viscosity solution. 
The RV solution when $c_{\text{RV}}$ in \eqref{eq:resviscosity:epsilon_rv_uw} is set to $1$ is more oscillatory compared with the 
RV solution when $c_{\text{RV}}$ is set to 4. Despite the oscillations, the important outcome is that the solution is stable in both cases.
Here we point out that these solutions are computed on scattered points. We made the same tests on Cartesian points, where 
we noticed that the oscillations were smaller. 
In Figure \ref{fig:burgers:ep_coefficients} we display the spatial variation of the viscosity coefficient $\varepsilon$ defined in 
\eqref{eq:resviscosity:epsilon}, for cases when the first-order viscosity is used and when RV is used. We see that the viscosity coefficient is large only 
in the vicinity of the shock, when RV is used, thus, the shocks are well captured. 

\begin{figure}[h!]
    \centering
    \setlength\tabcolsep{1.5pt}    
    \begin{tabular}{ccc}
        \multicolumn{3}{c}{\vspace{0.1cm}\textbf{Burger's equation: convergence in $2$-norm}} \\
        \hspace{0.8cm}\textbf{First-order viscosity} & \hspace{0.9cm} \textbf{RV ($c_{\text{RV}}=1$)} & \hspace{0.7cm} \textbf{RV ($c_{\text{RV}}=4$)} \\
        \includegraphics[width=0.31\linewidth]{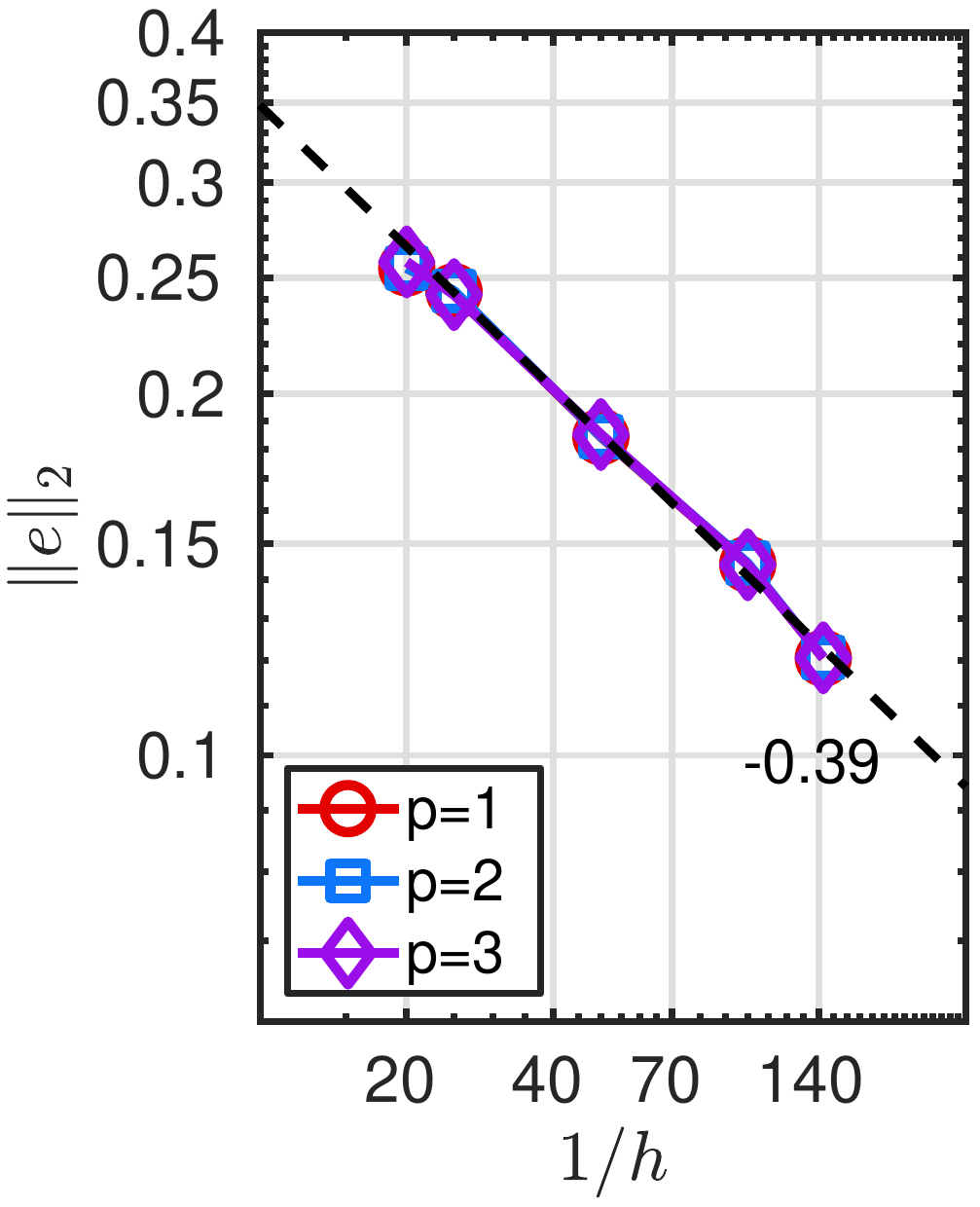} & 
        \includegraphics[width=0.31\linewidth]{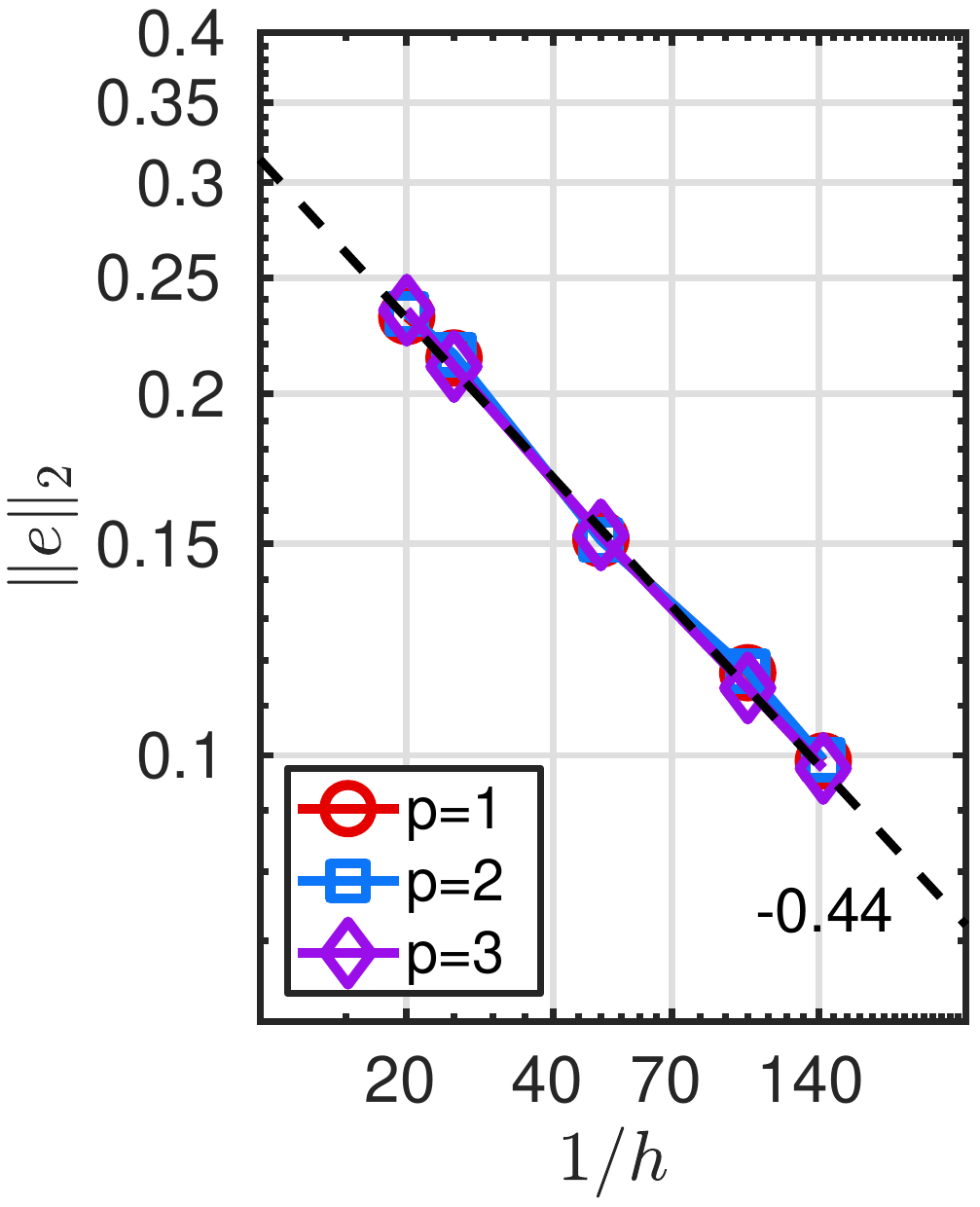} & 
        \includegraphics[width=0.31\linewidth]{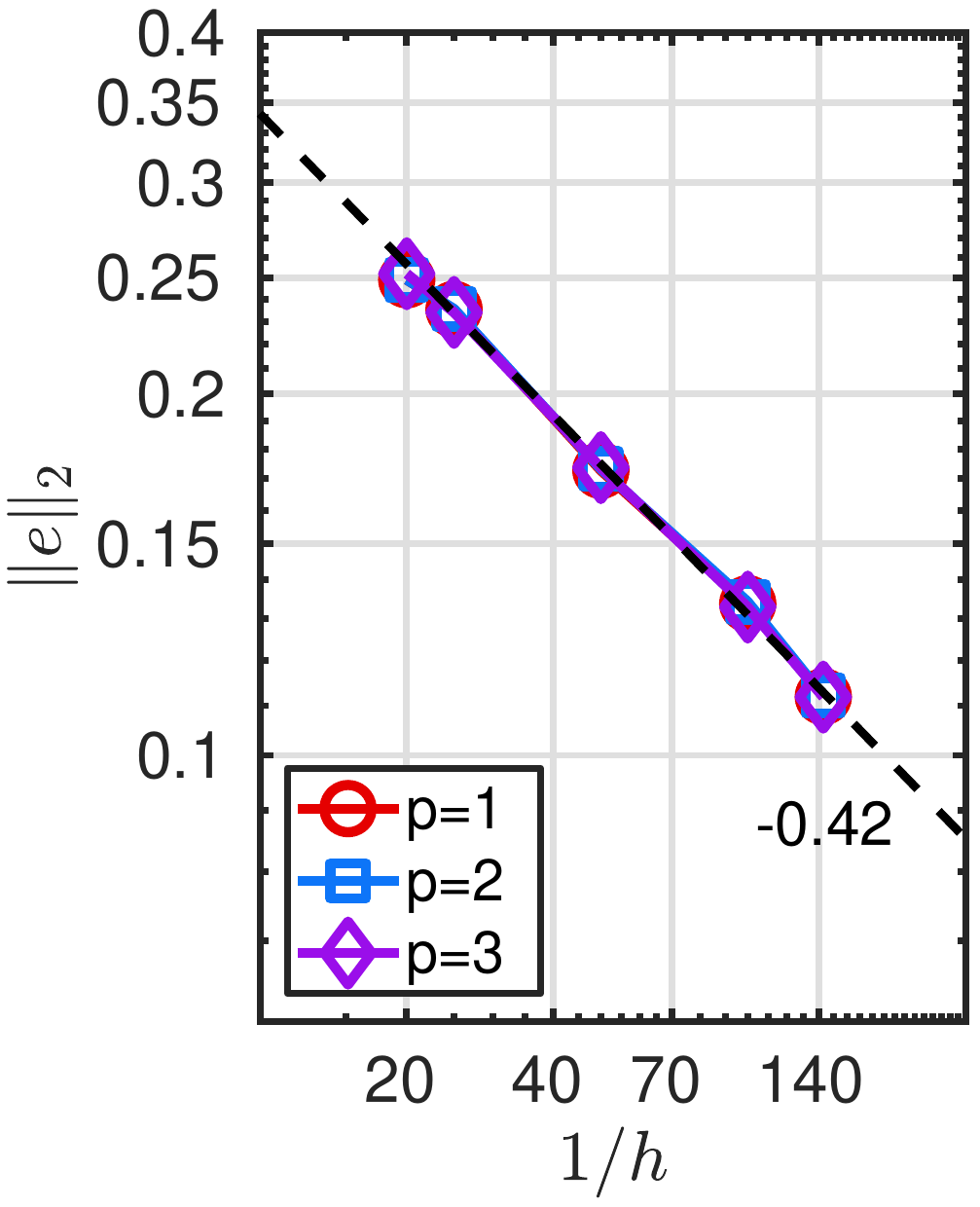} 
    \end{tabular}
    \caption{Convergence in $2$-norm, under node refinement for different choices of the monomial basis degree $p$, used to 
    construct the stencil-based approximants. The convergence plots are displayed for cases when the oversampled RBF-FD 
    is shock-stabilized using first-order viscosity and the residual viscosity.}
    \label{fig:burgers:error_convergence_l2}
\end{figure}

\begin{figure}[h!]
    \centering
    \setlength\tabcolsep{1.5pt}    
    \begin{tabular}{ccc}
        \multicolumn{3}{c}{\vspace{0.1cm}\textbf{Burger's equation: convergence in $1$-norm}} \\
        \hspace{0.8cm}\textbf{First-order viscosity} & \hspace{0.9cm} \textbf{RV ($c_{\text{RV}}=1$)} & \hspace{0.7cm} \textbf{RV ($c_{\text{RV}}=4$)} \\
        \includegraphics[width=0.31\linewidth]{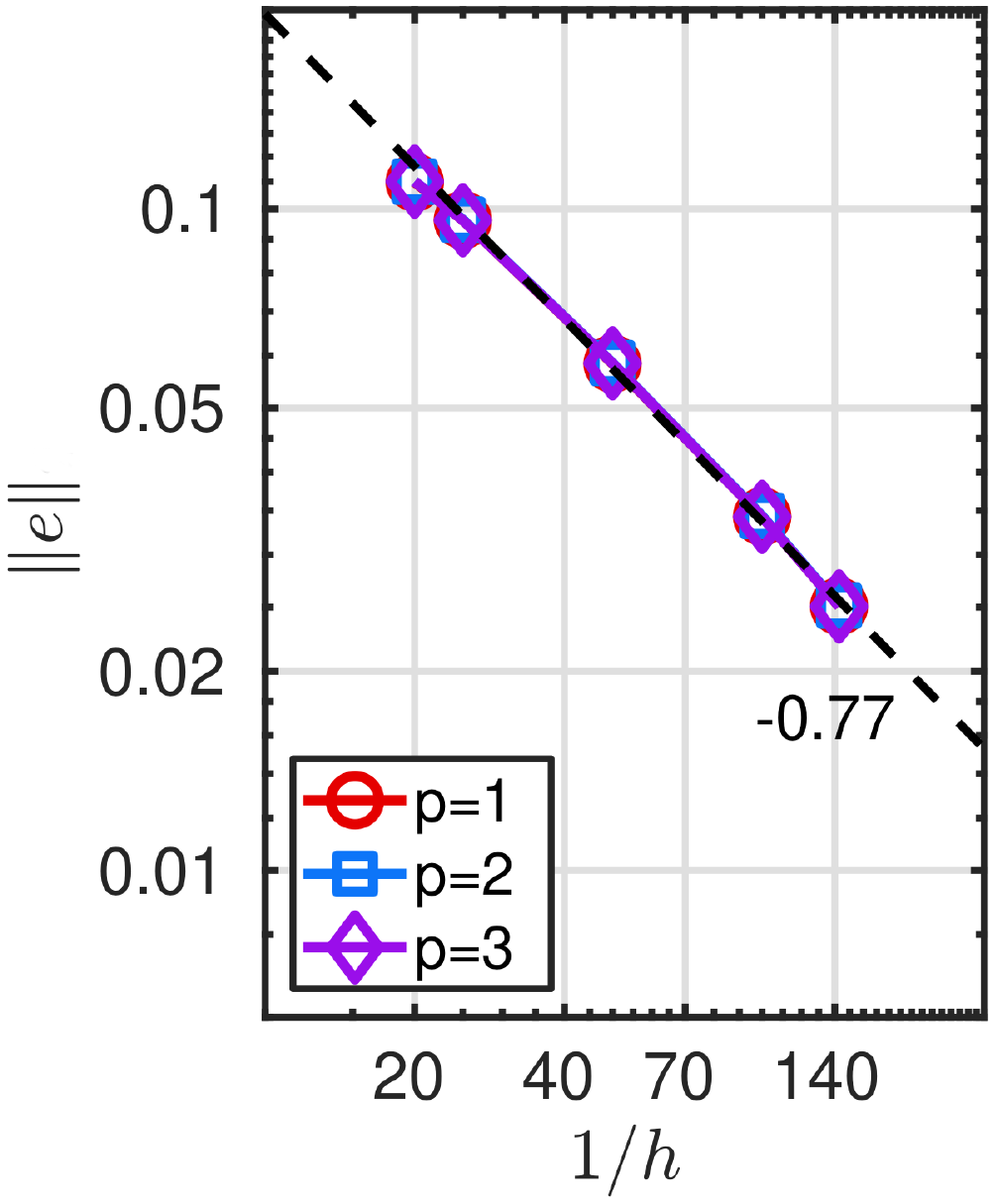} & 
        \includegraphics[width=0.31\linewidth]{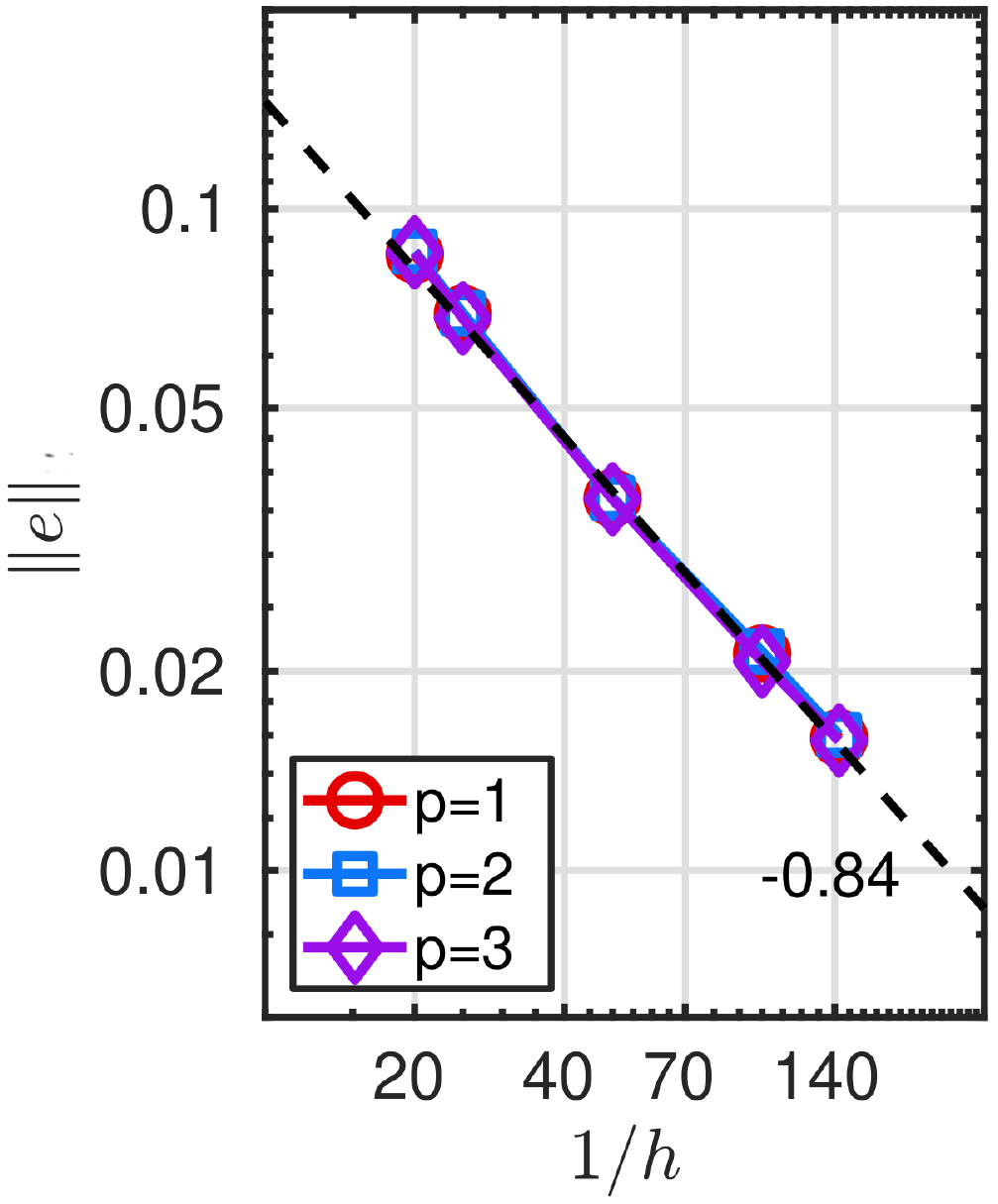} & 
        \includegraphics[width=0.31\linewidth]{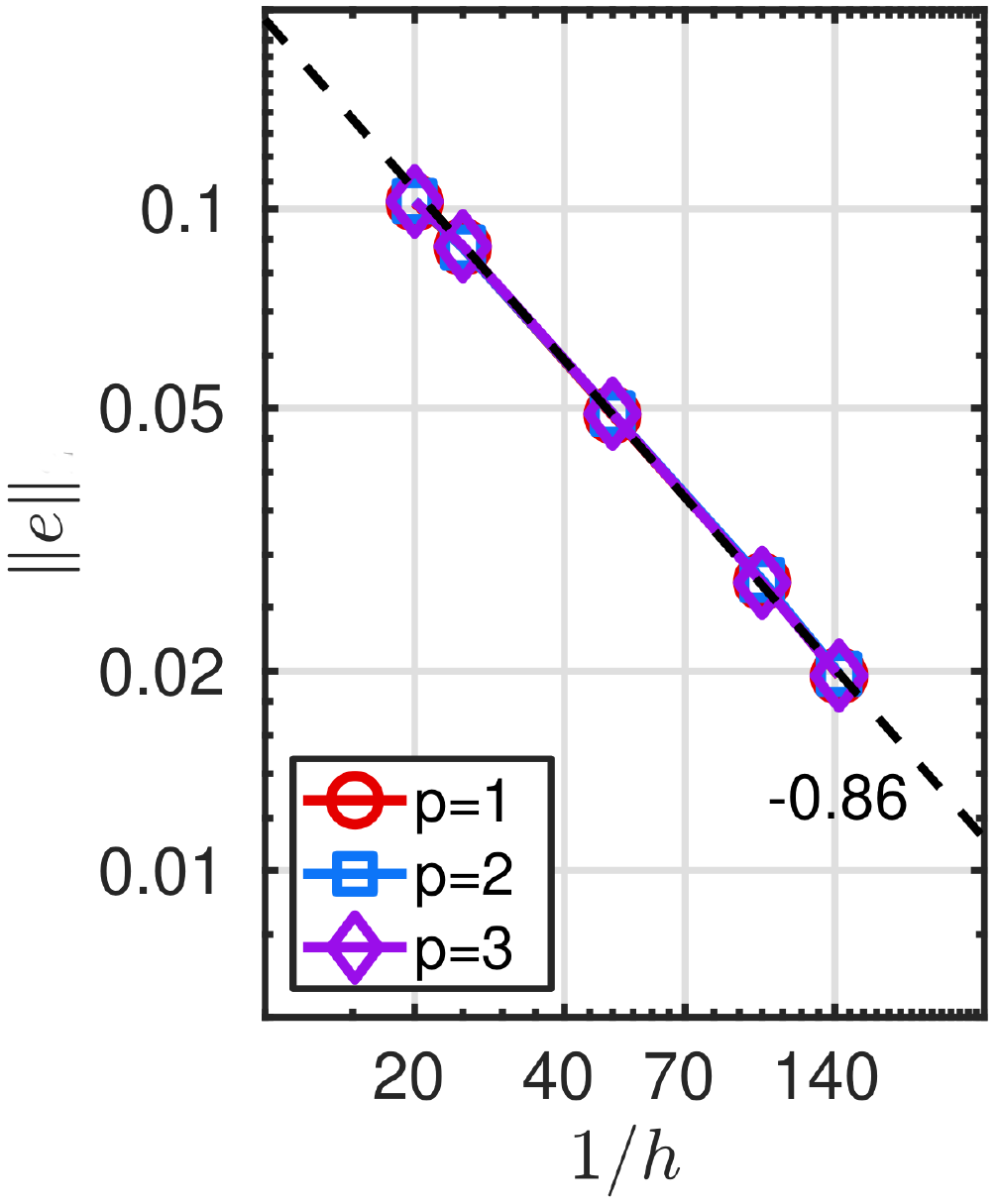} 
    \end{tabular}
    \caption{Convergence in $1$-norm, under node refinement for different choices of the monomial basis degree $p$, used to 
    construct the stencil-based approximants. The convergence plots are displayed for cases when the oversampled RBF-FD 
    is shock-stabilized using first-order viscosity and the residual viscosity.}
    \label{fig:burgers:error_convergence_l1}
\end{figure}

In Figure \ref{fig:burgers:error_convergence_l2} and Figure \ref{fig:burgers:error_convergence_l1}, we display the convergence under node refinement, in 
$2$-norm and $1$-norm respectively, where the norms are defined in \eqref{eq:experiments_linear:errors}. The optimal convergence rates are expected to be $0.5$ and $1$. 
The convergence trends in all stabilized cases are close to these optimal rates. When stabilizing the scheme using the first-order viscosity, the error is constant-wise larger 
compared to when stabilizing the scheme using RV. We do not observe that an increase of $p$ would -- in the present case -- lead to a significant improvement of the error by means of 
a constant.

\section{Numerical study III: The Kurganov-Petrova-Popov rotating wave problem}
\label{sec:experiments:kpp}

We solve another scalar conservation law with a nonlinear flux: the Kurganov-Petrova-Popov (KPP) problem, initially introduced in \cite{KPP}. 
The flux is given by $\bm F(u) = (\sin u, \cos u)$. The computational domain is a square $\Omega = [-2, 1.5] \times [2, 2.5]$. 
The initial condition and the Dirichlet boundary condition are in the same order defined as:
\begin{equation}
  u(y,0) = 
  \begin{cases}
    \frac{14 \pi}{4} & \text{if } \sqrt{y_1^2 + y_2^2} \leq 1 \\
    \frac{\pi}{4} & \text{otherwise}.
  \end{cases},
  \qquad u(y,t)\big |_{\partial\Omega} = \frac{\pi}{4}.
\end{equation}
We use the monomial basis degree $p=3$ to construct the stencil based approximations, 
on a scattered point set $X$ with a mean internodal distance 
$h=0.02$ (corresponding to $39498$ nodes). The point set is obtained using the algorithm introduced in \cite{FBF15_nodes}. 
The oversampling parameter is set to $q=5$ and the RV constant to $C_{\text{RV}} = 3$. 
We run the simulation until time $t=1$ with the CFL number $0.2$. 
Since the velocity field is for this benchmark rotational we use a fixed time step $\Delta t$, computed according to \eqref{eq:experiments_linear:timestep}.

\begin{figure}[h!]
    \centering
    \setlength\tabcolsep{1.5pt}    
    \begin{tabular}{ccl}
        \multicolumn{3}{c}{\hspace{-0.7cm}\vspace{0.1cm}\textbf{The KPP problem: unstabilized numerical solutions}} \\
        \includegraphics[width=0.37\linewidth]{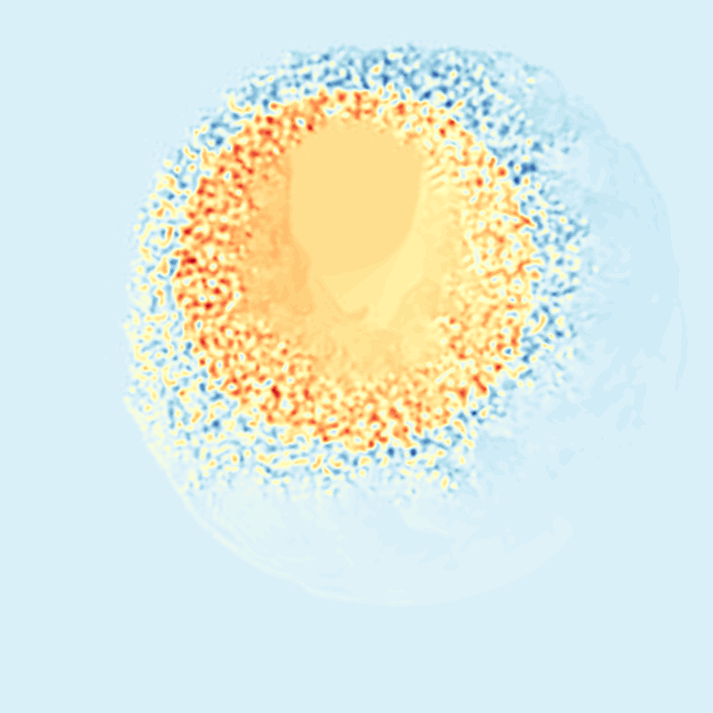} &
        \hspace{1cm}
        \includegraphics[width=0.12\linewidth]{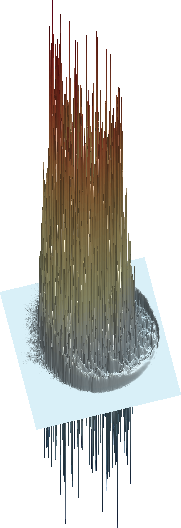} &
        \hspace{0.2cm}
        \includegraphics[width=0.07\linewidth]{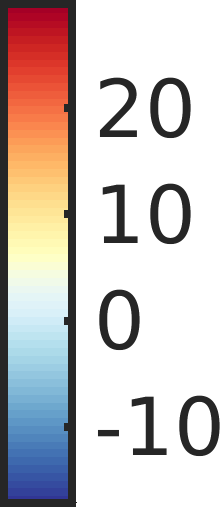}
    \end{tabular}
    \caption{Solution to the Kurganov-Petrova-Popov problem, when no shock stabilization is added to the oversampled RBF-FD discretization. 
    The hyperviscosity term in \eqref{eq:discretization:advection_final} is active in order to stabilize the numerical scheme in time.}
    \label{fig:kpp:solutions_unstabilized}
\end{figure}

\begin{figure}[h!]
    \centering
    \setlength\tabcolsep{1.5pt}    
    \begin{tabular}{cclcl}
        \multicolumn{3}{c}{\hspace{-1cm}\vspace{0.1cm}\textbf{The KPP problem: stabilized numerical solutions}} \\
        \vspace{0.05cm}\hspace{0cm}\textbf{First-order viscosity} & \hspace{-0.3cm}\textbf{RV} & \\
        \includegraphics[width=0.37\linewidth]{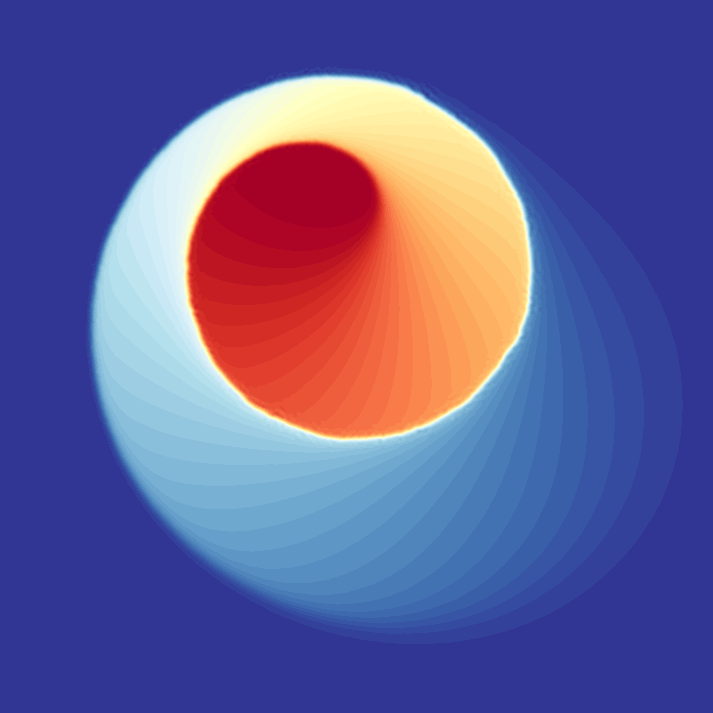} &
        \includegraphics[width=0.37\linewidth]{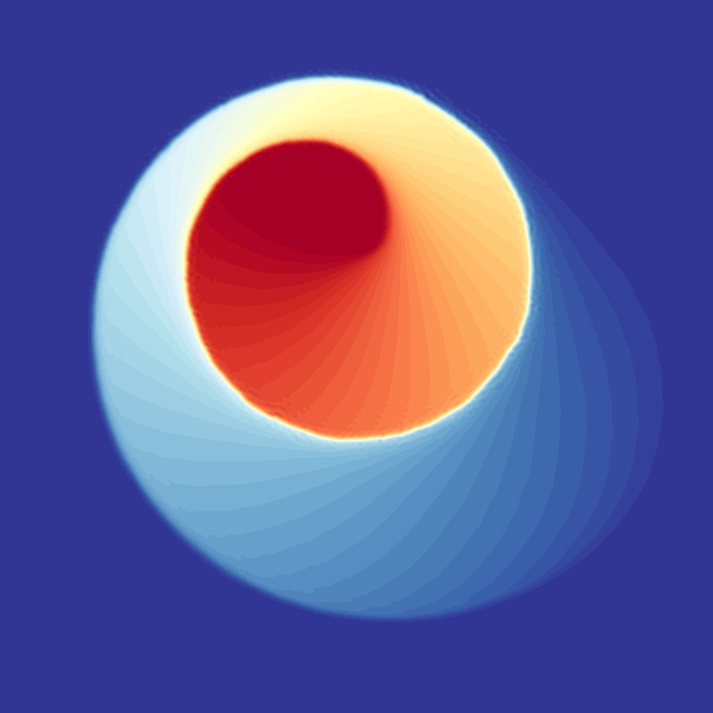} & 
        \includegraphics[width=0.07\linewidth]{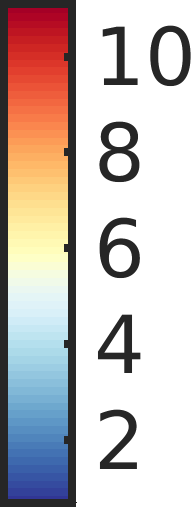}  \\
        \includegraphics[width=0.37\linewidth]{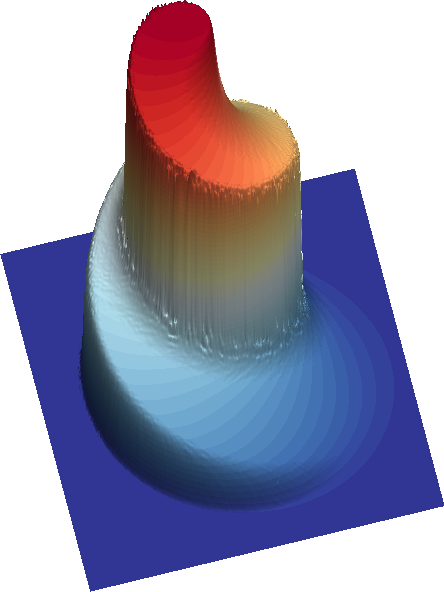} &
        \includegraphics[width=0.37\linewidth]{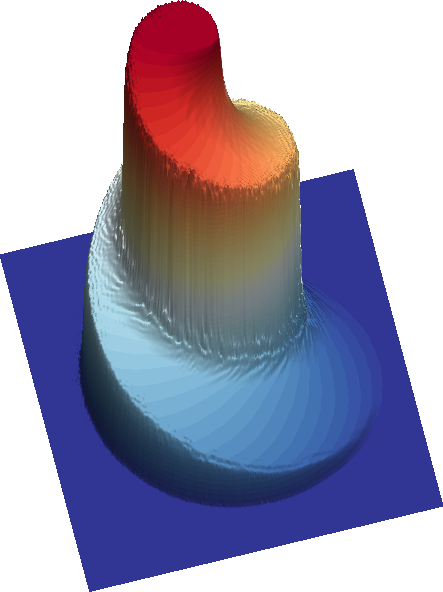} & 
        \includegraphics[width=0.07\linewidth]{pics/Pics_KPP/numerical_solution_colorbar_RV.png}        
    \end{tabular}
    \caption{A numerical solution to the Kurganov-Petrova-Popov problem from a top perspective (first row) and a side perspective (second row), 
    when two different shock stabilizations are added to the oversampled RBF-FD discretization. 
    The hyperviscosity term in \eqref{eq:discretization:advection_final} is active in addition, in order to stabilize the numerical scheme in time.}
    \label{fig:kpp:solutions}
\end{figure}

\begin{figure}[h!]
    \centering
    \setlength\tabcolsep{1.5pt}    
    \begin{tabular}{ccl}
        \multicolumn{3}{c}{\hspace{-1.5cm}\vspace{0.1cm}\textbf{The KPP problem: viscosity coefficients}} \\
        \vspace{0.2cm}\hspace{-0.3cm}\textbf{First-order viscosity} & \hspace{-0.3cm}\textbf{RV} & \\
        \includegraphics[width=0.37\linewidth]{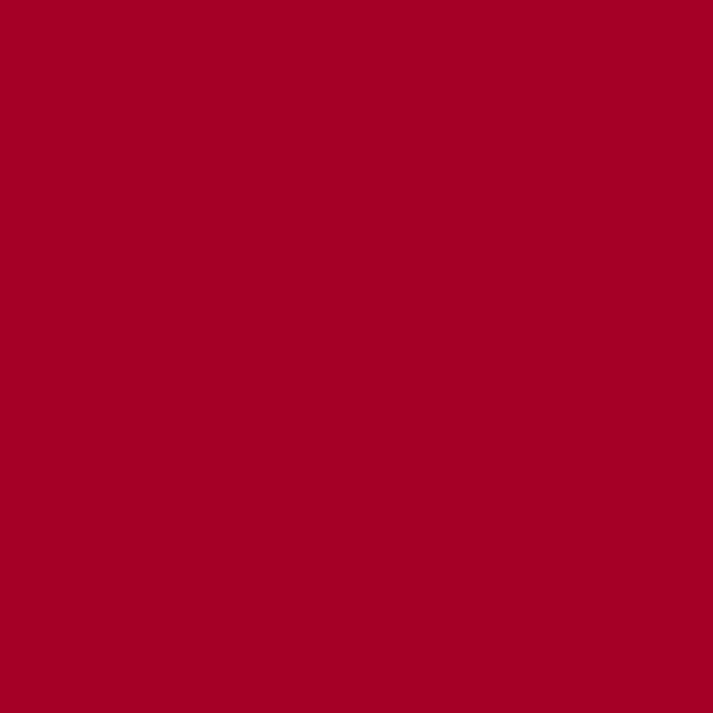} &
        \includegraphics[width=0.37\linewidth]{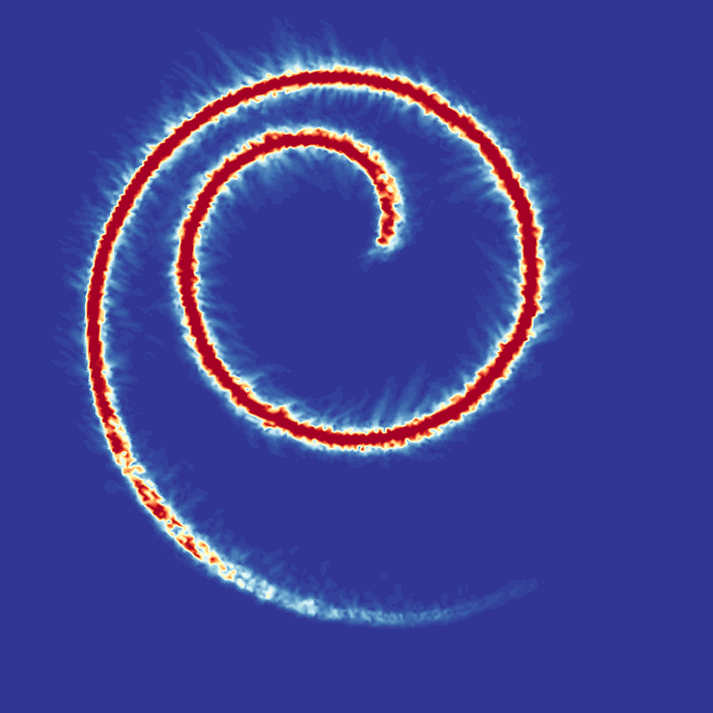} & 
        \includegraphics[width=0.11\linewidth]{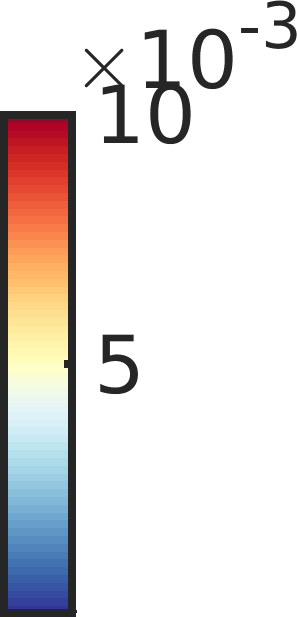}
    \end{tabular}
    \caption{The images display the spatial variation of the viscosity coefficient when using the first-order viscosity and the 
    residual viscosity.}
    \label{fig:kpp:coefficients}
\end{figure}

The solution at $t=1$ when the shocks are not stabilized is displayed in Figure \ref{fig:kpp:solutions_unstabilized}, 
where we see oscillations so large, that the solution can not be seen as physical.
Stabilized solutions are displayed in Figure \ref{fig:kpp:solutions}. The solutions stabilized with the first-order viscosity and the residual viscosity 
are in an overall sense similar. A closer 
look, however, reveals that the residual viscosity solution includes sharper details along the edge of the spiral. The RV solution is less diffused 
compared to the first-order viscosity solution.
The viscosity coefficients are given in Figure \ref{fig:kpp:coefficients}, where we see that the residual viscosity is, 
roughly speaking, only active in the region of the discontinuity.

\section{Numerical study IV: compressible Euler equations}
\label{sec:experiments:Euler}
In this section we solve compressible Euler equations discretized in \eqref{eq:discretization:Euler_final}. 
The Euler system can be seen as a limiting case of the Navier-Stokes 
system of equations, when the physical viscosity term tends to $0$. 
We solve several benchmarks which model a fluid in the state of a compressed gas. For all experiments 
we use the adiabatic constant $C_{\text{ad}} = 1.4$.
The variable time steps for all considered cases are computed according to the formula:
\begin{equation}
    \Delta t = \text{CFL}\, \min_{x_i \in X} \left [\frac{h_{\text{loc}}(x_i)}{\sqrt{v_1^2(x_i, t) + v_2^2(x_i, t)} + C_{\text{ad}}\, p(x_i, t)/\rho(x_i, t)} \right ],
\end{equation}
where $h_{\text{loc}}$ is defined within the context of \eqref{eq:resviscosity:epsilon_rv_uw}.
The approximation errors for each unknown function are computed analogously to \eqref{eq:experiments_linear:errors}.

\subsection{Sod's shock tube problem}
The computational domain is a rectangle $\Omega = [0, 1] \times [0,0.4]$, in which two fluids with different physical properties 
are separated by a membrane positioned at 
$y_1 = 0.5$. The initial condition is given by:
\begin{equation}
    \begin{aligned}
y_1 < 0.5:\quad & \rho = 1,\, & \bm v = (0,0), & \quad p = 1, \\
y_1 > 0.5:\quad & \rho = 0.125,\, & \bm v = (0,0), &\quad p = 0.1.
    \end{aligned}
\end{equation}
The simulation of the two fluids is started by removing the membrane at $t=0$, the final state is observed at $t=0.25$.
The boundary condition is \emph{slip} ($\bm v \cdot \bm n = 0$) along all boundaries of the rectangle, where $\bm n$ is the outward pointing normal.

We use scattered nodes generated using an algorithm 
from \cite{FBF15_nodes} to discretize $\Omega$.  Throughout the section, the CFL condition is set to $0.2$, the oversampling parameter is set to $q=5$ and the 
monomial basis degree to $p=3$, if not stated otherwise.

\begin{figure}[h!]
    \centering
    \begin{tabular}{ccc}
        \multicolumn{3}{c}{\textbf{Euler, Sod's shock tube: numerical solution}} \\
    \textbf{Top view} & \textbf{Side view} & \\
    \includegraphics[width=0.4\linewidth]{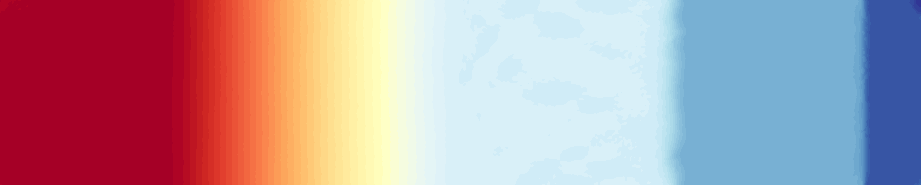} &
    \includegraphics[width=0.4\linewidth]{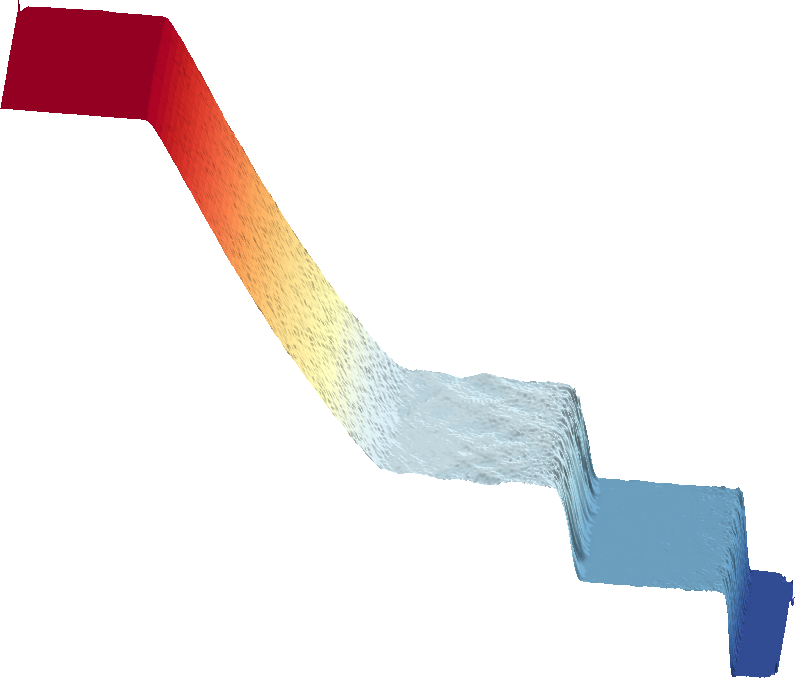}  &
    \includegraphics[width=0.06\linewidth]{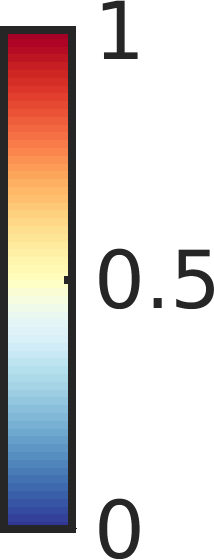} 
    \end{tabular}
    \caption{The numerical solution (density $\rho$) of Sod's shock tube problem at $t=0.25$, computed using the oversampled RBF-FD method. 
    The internodal distance is set to $h=0.003$ (corresponds to $N=22244$ unknowns). 
    The monomial basis degree used to construct the stencil-based approximants is $p=3$. The oversampling parameter is set to $q=5$.}
    \label{fig:euler:shocktube:solution_2d}
\end{figure}

The numerical solution (density), obtained using the RV stabilized oversampled RBF-FD method, is displayed in Figure \ref{fig:euler:shocktube:solution_2d}. 
We observe that the shock discontinuity and the contact discontinuity are well captured. 
The mean distance between the computational nodes is set to $h=0.003$ ($N=22244$ unknowns).

\begin{figure}[h!]
    \centering
    \begin{tabular}{cc}
        \multicolumn{2}{c}{\textbf{Euler, Sod's shock tube: solutions over a horizontal cross-section}} \\
        \hspace{0.5cm} \textbf{Density} $\rho$ & \textbf{Temperature} $T$  \\
    \includegraphics[width=0.45\linewidth]{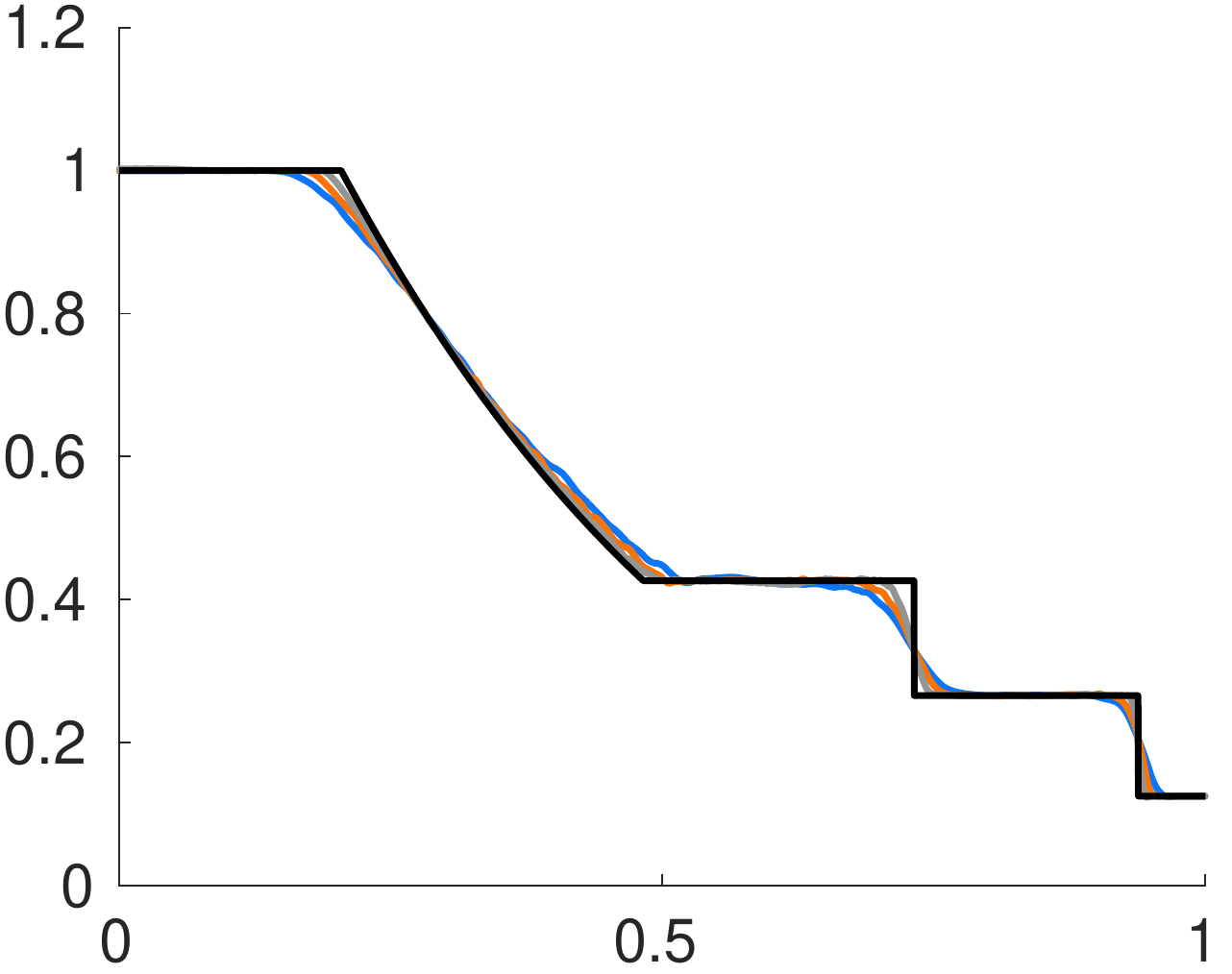} &
    \includegraphics[width=0.45\linewidth]{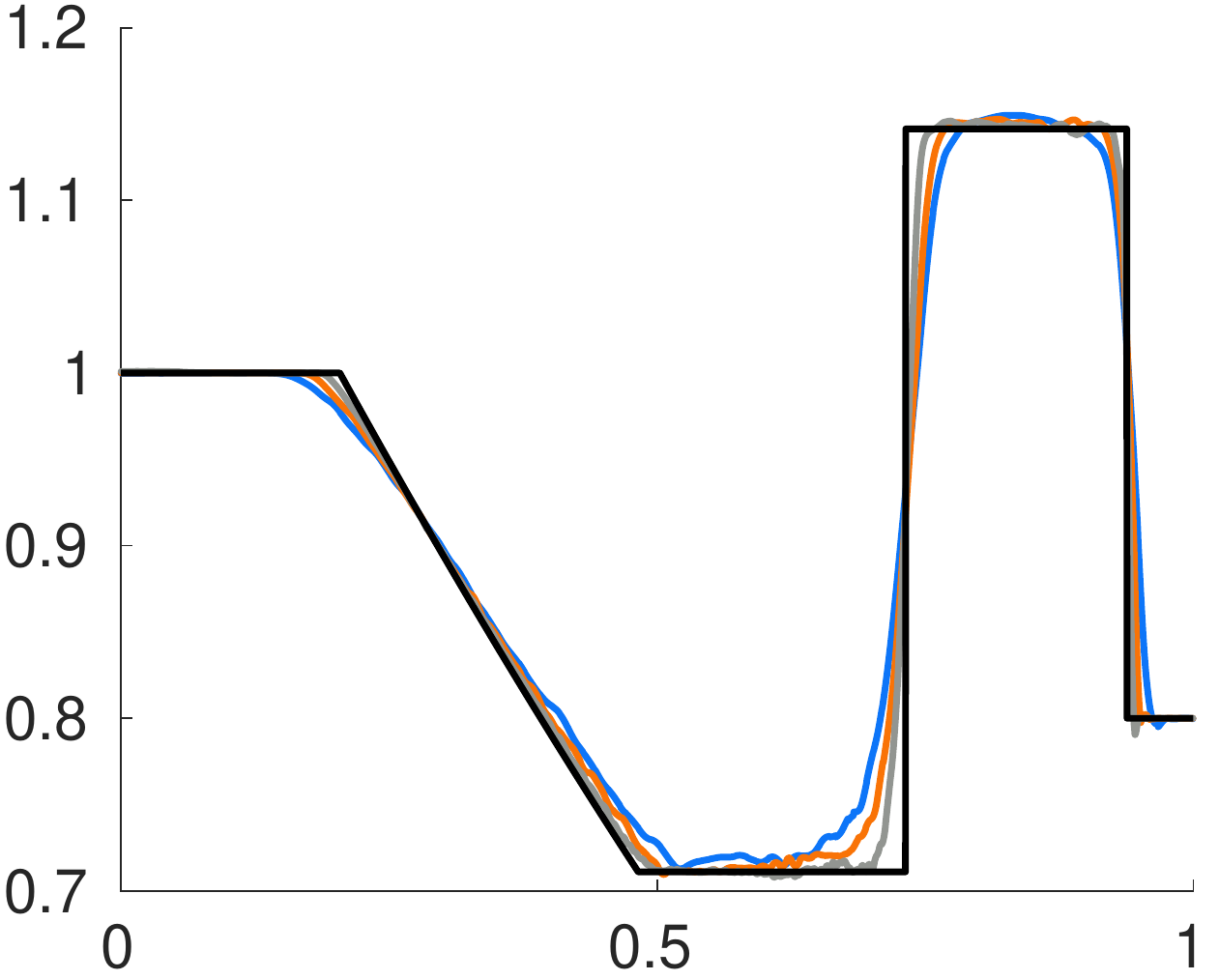}\\
    \multicolumn{2}{c}{\includegraphics[width=0.5\linewidth]{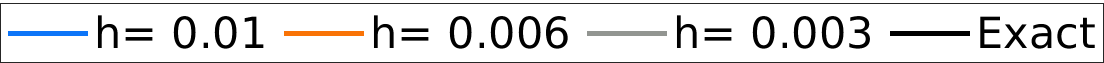}}
    \end{tabular}
    \caption{A slice of the solution of Sod's shock tube problem at $t=0.25$, computed using the oversampled RBF-FD method. The numerical solutions for different $h$ 
    are compared with 
    the exact solution. 
    The monomial basis degree used to construct the stencil-based approximants is $p=3$. The oversampling parameter is set to $q=5$.}
    \label{fig:euler:shocktube:solutions}
\end{figure}

\begin{figure}[h!]
    \centering
    \begin{tabular}{ccc}
        \multicolumn{3}{c}{\textbf{Euler, Sod's shock tube: convergence in $2$-norm}\vspace{0.05cm}} \\
        \hspace{1cm} \textbf{Density} & \hspace{1cm}\textbf{Momentum} & \hspace{1cm}\textbf{Total energy} \\ 
    \includegraphics[width=0.3\linewidth]{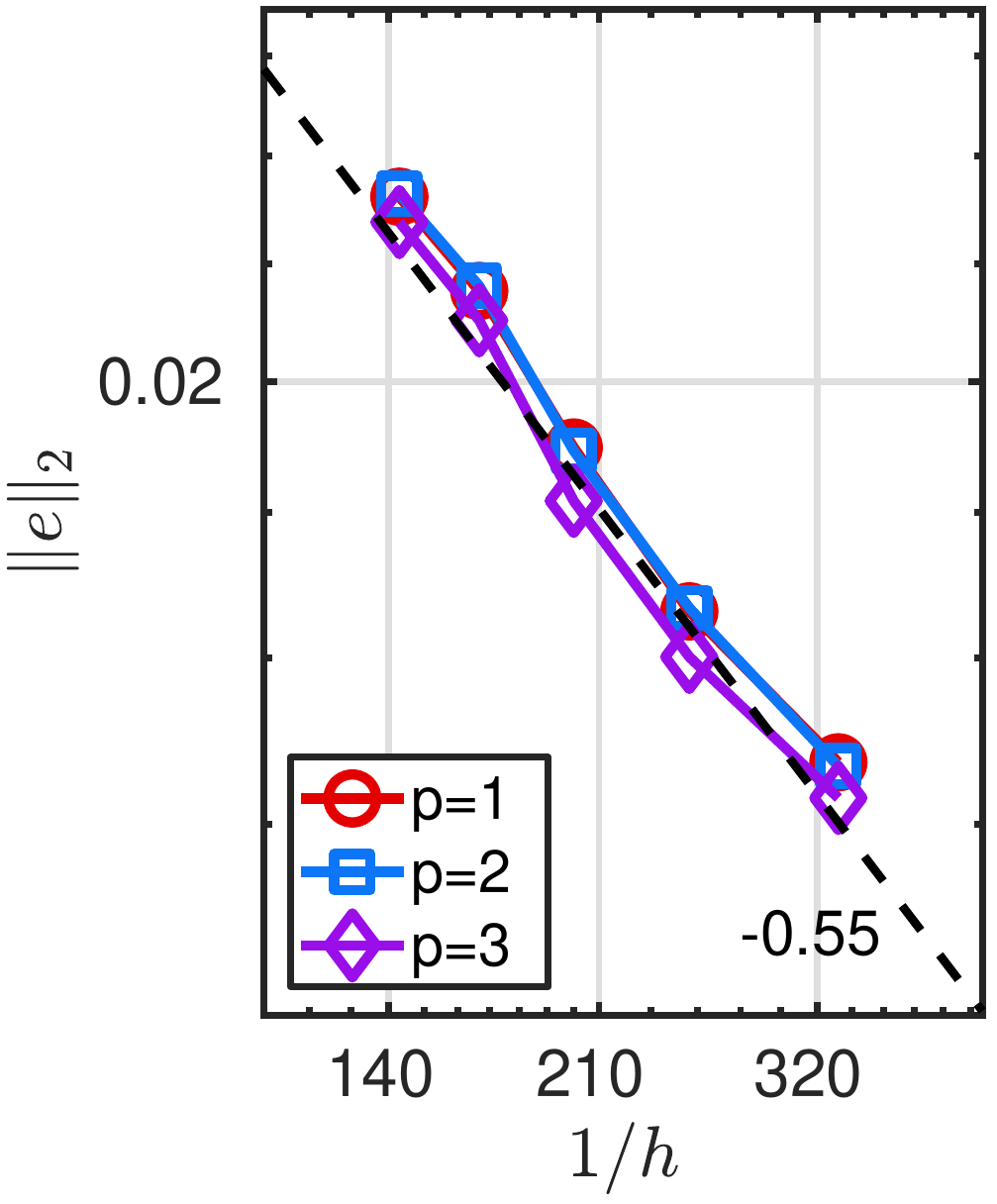} &
    \includegraphics[width=0.3\linewidth]{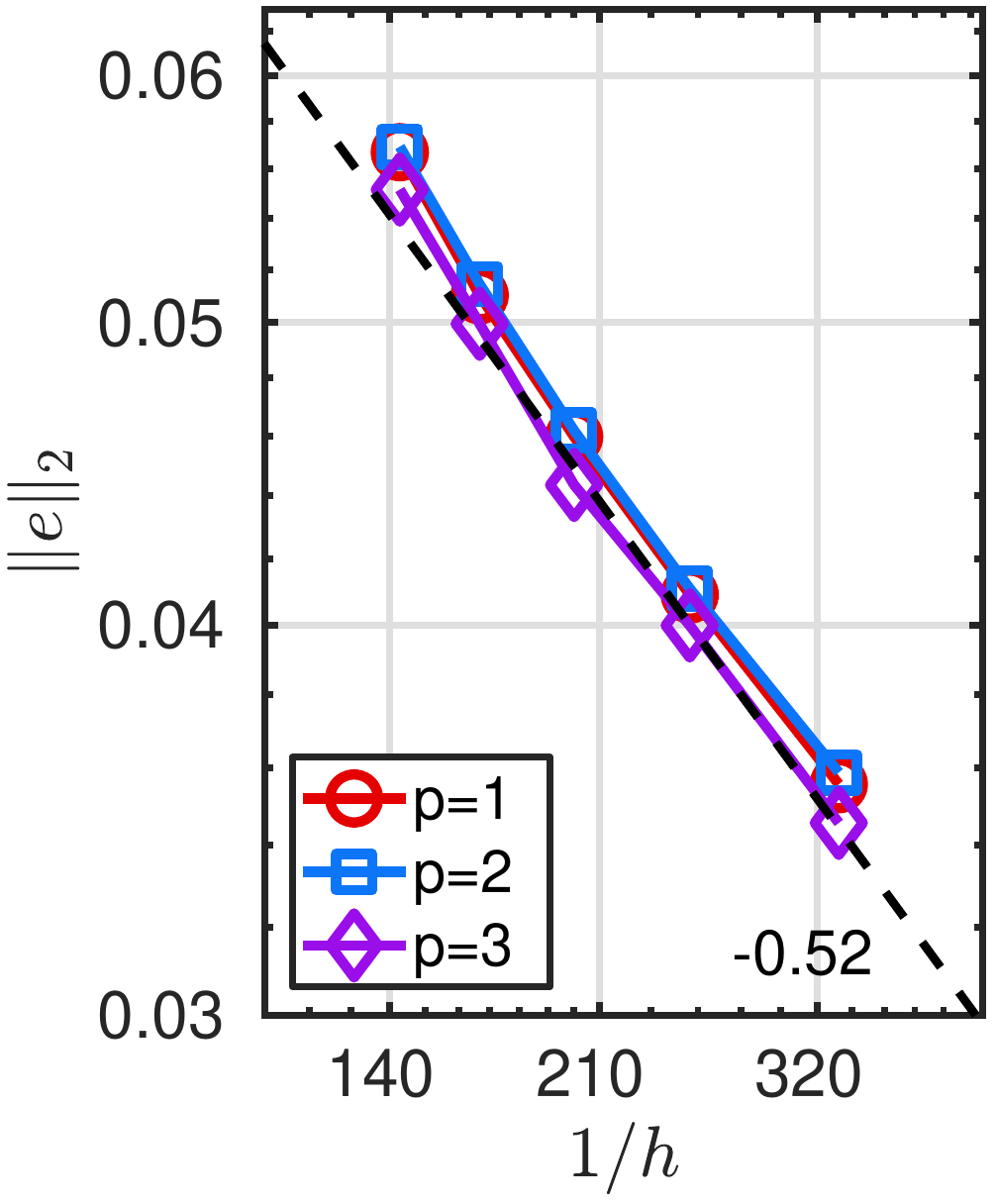} &
    \includegraphics[width=0.3\linewidth]{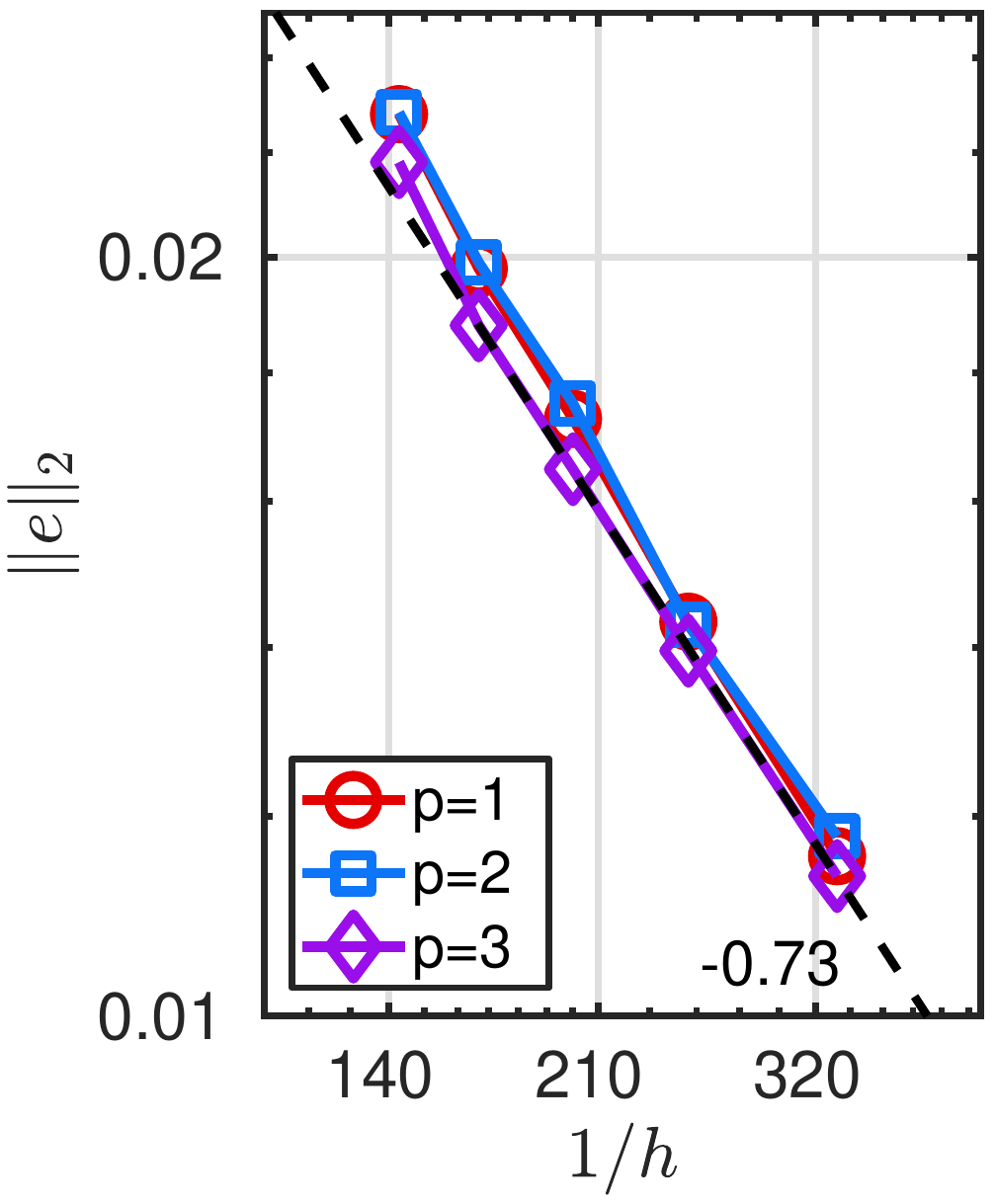}
    \end{tabular}
    \caption{Convergence in $2$-norm under node refinement, for different choices of the monomial basis degree $p$, used to 
    construct the stencil-based approximants.}
    \label{fig:euler:shocktube:error_l2}
\end{figure}

\begin{figure}[h!]
    \centering
    \begin{tabular}{ccc}
        \multicolumn{3}{c}{\textbf{Euler, Sod's shock tube: convergence in $1$-norm}\vspace{0.05cm}} \\
        \vspace{-0.08cm}\hspace{1cm} \textbf{Density} & \hspace{1cm} \textbf{Momentum} & \hspace{1cm} \textbf{Total energy} \\ 
    \includegraphics[width=0.3\linewidth]{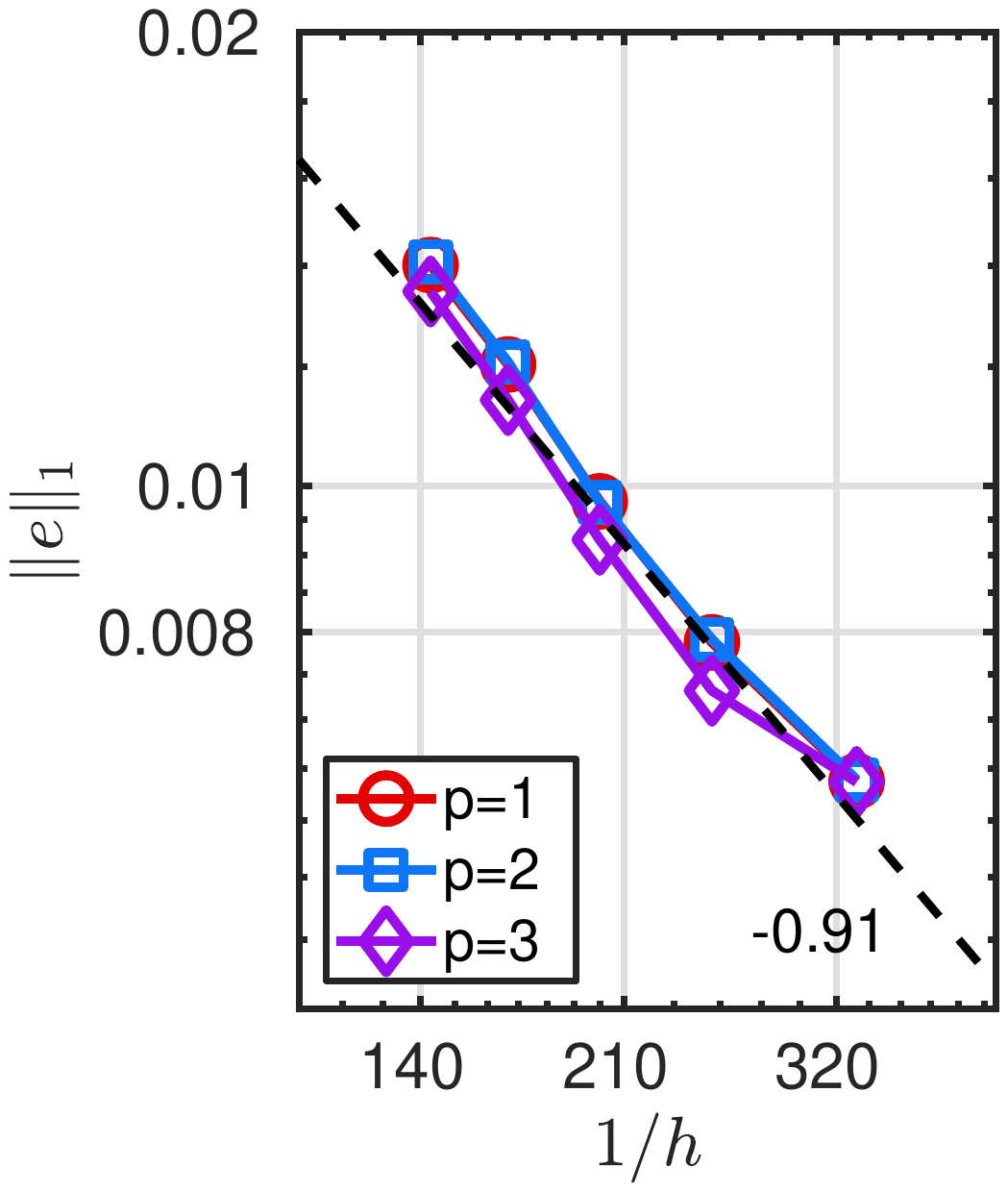} &
    \includegraphics[width=0.3\linewidth]{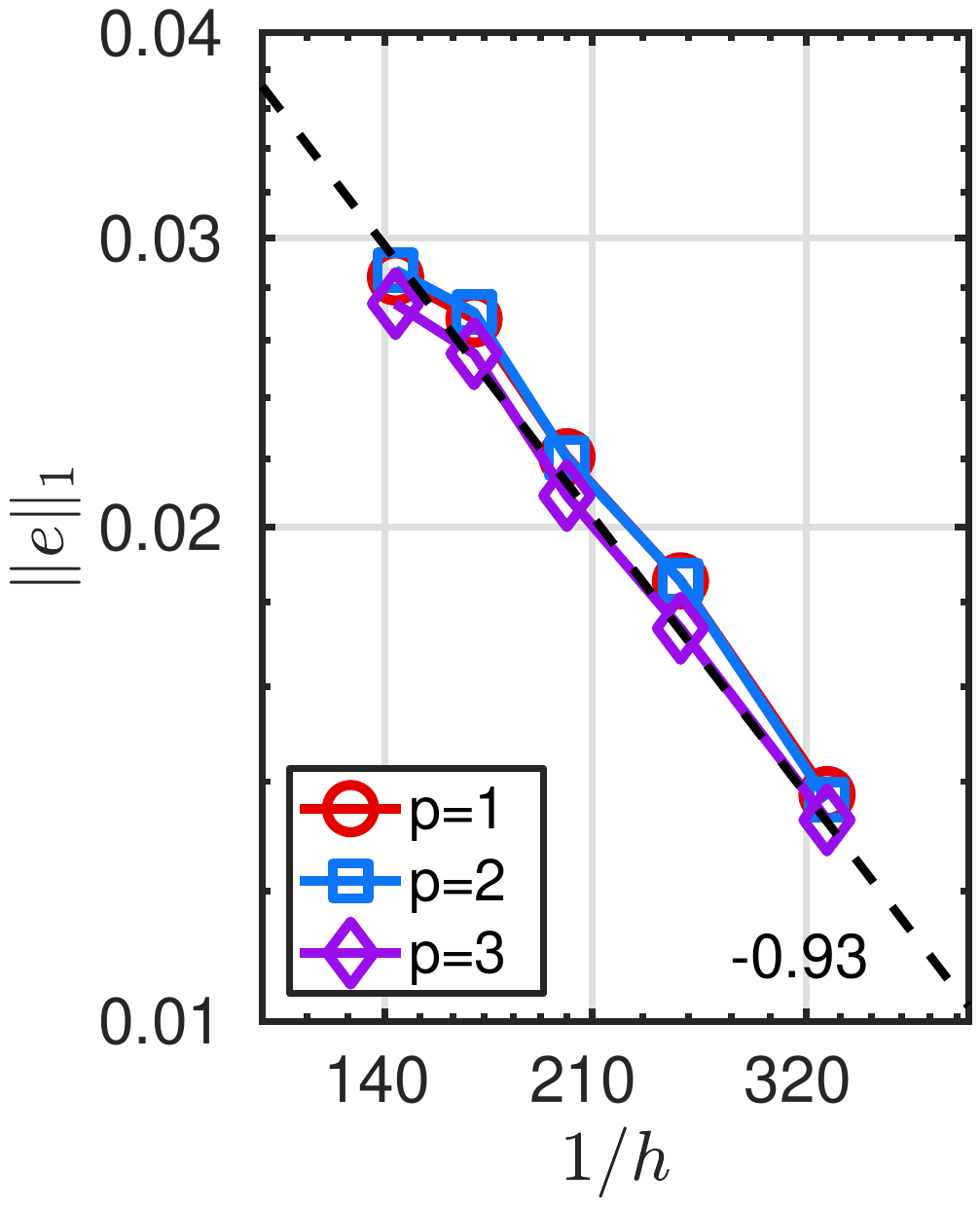} &
    \includegraphics[width=0.3\linewidth]{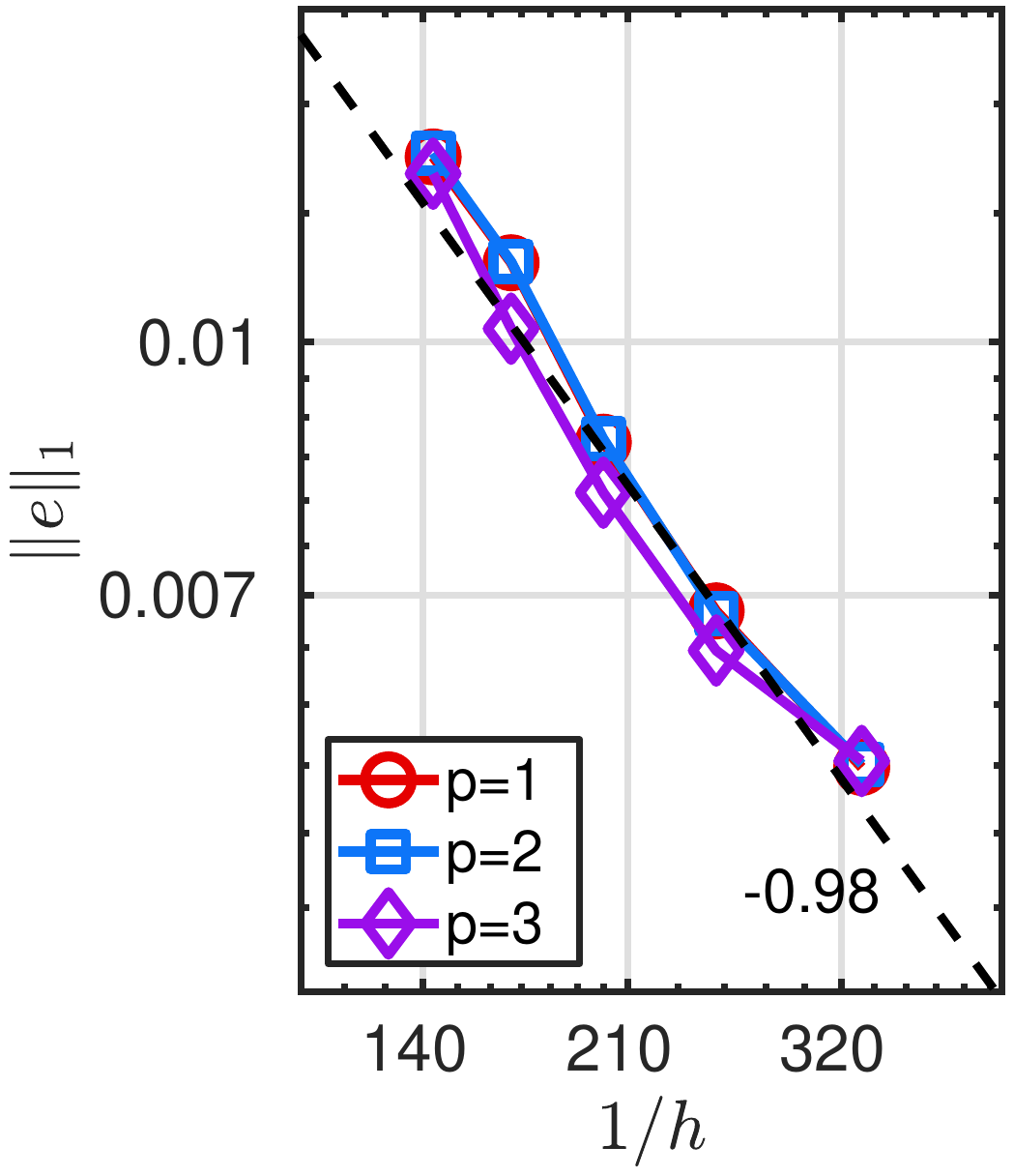}
    \end{tabular}
    \caption{Convergence in $1$-norm under node refinement, for different choices of the monomial basis degree $p$, used to 
    construct the stencil-based approximants.}
    \label{fig:euler:shocktube:error_l1}
\end{figure}

Now we compare the RV stabilized RBF-FD solution with the exact solution. The exact solution is obtained for an equivalent problem in 
one dimension. We interpolate the 2-dimensional RBF-FD solution along the line $y_2 = 0.2$ in order to obtain an approximate 1-dimensional solution. 
This can be done since the variation of the 2-dimensional solution in the $y_2$ direction is very small. 
The first result is displayed in Figure \ref{fig:euler:shocktube:solutions}, where we show exact and numerical profiles of the density and the temperature. 
We observe that the discontinuous exact solution is well approximated by the numerical solution, 
and that the approximation gets increasingly better with a decreasing $h$.
In Figure \ref{fig:euler:shocktube:error_l2} we plot the $2$-norm convergence of the numerical solution along the line $y_2 = 0.2$, 
where the exact 1D solution is taken as a reference. The optimal convergence rates for this norm are expected to be $0.5$. For the density and the momentum, 
our convergence rates are optimal, while the converence rate for the total energy is slightly higher than optimal. We also observe that for the density, 
the approximation error in $p=3$ case 
is slightly smaller compared to $p=2$ and $p=1$.
In Figure \ref{fig:euler:shocktube:error_l1} we show an analogous convergence plot, but in $1$-norm. 
Here the optimal convergence rate is $1$. The convergence plots in all cases are close to the optimum. An increase in $p$ 
gives a slight improvement of the approximation error by means of a constant. 

In \cite{NazarovLarcher17}, the authors used the RV stabilized finite element method to solve Sod's shock tube problem using the same set 
of physical parameters. Comparing Figure \ref{fig:euler:shocktube:error_l2} and Figure \ref{fig:euler:shocktube:error_l1} from the present paper, 
with Figure 6 from \cite{NazarovLarcher17}, we can observe that the errors in $2$-norm and $1$-norm are comparable in magnitude and the slopes. 

\subsection{A 2D Riemann problem}
The domain is a square $\Omega = [0,1] \times [0,1]$. The initial condition is defined such that it 
takes different values in different quadrants of the square:
\begin{equation}
\begin{aligned}
\rho &= 4/5, &  p&=1, & \bm v &= (0, 0), & \text{ in } & 0< y_1 <0.5, & 0 < y_2 < 0.5, \\
\rho &= 1,  & p&=1, & \bm v &= (0.7276, 0), & \text{ in } & 0< y_1 <0.5, & 0.5 < y_2 < 1, \\
\rho &= 1, & p&=1, & \bm v &= (0, 0.7276), & \text{ in } & 0.5 < y_1 < 1, & 0 < y_2 < 0.5, \\
\rho &= 17/32, &  p&=2/5, & \bm v &= (0, 0), & \text{ in } & 0.5 < y_1 < 1, & 0.5 < y_2 < 1.
\end{aligned}
\end{equation}
The boundary conditions are \emph{slip} on the top side and the right side of the square domain. 
On the bottom and the left side of the square domain, we use the inflow Dirichlet boundary conditions. These are for all $t$ 
set to the corresponding values of the initial condition. The scattered nodes over $\Omega$ are placed using the algorithm provided in 
\cite{FBF15_nodes}. We set the monomial basis degree to $p=3$, the oversampling parameter to $q=5$, the internodal distance to $h=0.0025$ (corresponds to $N=156981$ unknowns), 
the CFL number to $0.2$ and the RV constant to $C_{\text{RV}} = 3$, and run the simulation until $t=0.25$. 

\begin{figure}[h!]
    \centering
    \begin{tabular}{cc}
        \multicolumn{2}{c}{\textbf{Euler, Riemann problem}} \\
        \hspace{-0.95cm}\textbf{Numerical solution} & \hspace{-1.2cm}\textbf{RV coefficient} \\
\includegraphics[width=0.45\linewidth]{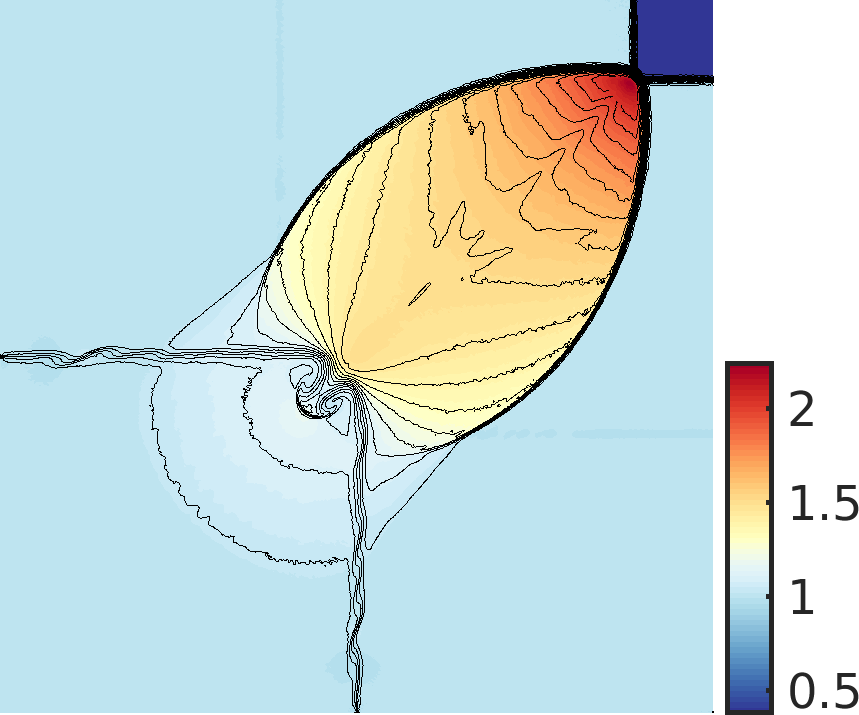} &
\includegraphics[width=0.475\linewidth]{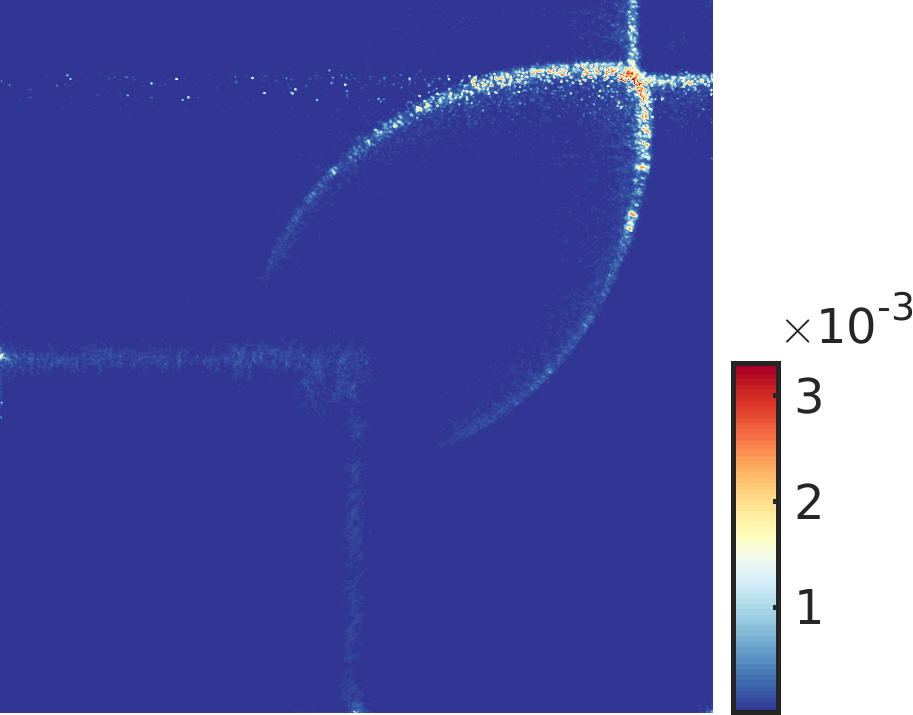}
    \end{tabular}
    \caption{The numerical solution of the Riemann problem, together with the spatial distribution of the residual viscosity coefficient. 
    The colormap of the numerical solution shows the pressure in the interval $[0.37, 2.23]$, 
    while the contour lines show the density in the interval $[0.54. 1.7]$.}
    \label{fig:euler:riemann:solutions}
\end{figure}

The result is given in Figure \ref{fig:euler:riemann:solutions}, where we display the RBF-FD solution stabilized using RV 
and the corresponding spatial distribution of the viscosity coefficient. We observe that the solution is sharp around the shock in the upper right 
corner of the domain, and that the contact discontinuity in the lower left quadrant is well defined. Furthermore, the magnitude of the viscosity coefficient 
is largest in the shock region and slightly smaller in the region of the contact discontinuity. The solution looks very similar to the one 
obtained using the RV stabilized finite element method in Figure 7 of \cite{NazarovLarcher17}, and also to all of the solutions obtained using six different 
shock-capturing schemes in Figure 4.2 of \cite{LiskaWendroff03}. 

\subsection{A channel flow over a forward facing step}
We study a Mach 3 supersonic flow over a forward facing step in a rectangular channel 
$\Omega = [0,3] \times [0,1]$. 
The sharp inward corner of the forward facing step is located at coordinate $(0.6,0.2)$.
The initial condition is: $\rho = 1.4$, $\bm m = (4.2, 0)$, $\mathcal{E} = 8.8$. We use 
slip boundary conditions on the top and the bottom boundaries of the channel. On the left boundary of the channel we use an 
inflow (Dirichlet) boundary condition, where for all $t$ we prescribe the values of the initial condition. On the right 
side of the channel we mimic the outflow by not imposing any boundary conditions. We run the simulation until $t=4$ 
using $N=117432$ nodes obtained using Gmsh \cite{Gmsh}. Other parameters are $C_{\text{RV}} = 5$, $\text{CFL} = 0.3$. 
We also consider two choices of monomial basis degrees: $p=1$ and $p=3$.

The flow entering the channel hits the step and then creates a bow shock that propagates towards the upper boundary, 
and then keeps reflecting between upper and lower boundaries. Close to the first reflection where the shocks meet, 
there is a so-called triple point. From the direction of the triple point towards the outflow, 
the physical solution has a contact discontinuity, which is challenging to capture in a numerical sense.

\begin{figure}[h!]
    \centering
    \begin{tabular}{ccc}
        \multicolumn{3}{c}{\textbf{Euler, Flow over a forward facing step} \vspace{0.1cm}} \\
        & \hspace{-0.5cm} \textbf{Numerical solution} & \hspace{-1cm} \textbf{RV coefficient} \\
        \raisebox{2\normalbaselineskip}[0pt][0pt]{\rotatebox[origin=c]{90}{$\mathbf{p=1}$}} & 
        \includegraphics[width=0.42\linewidth]{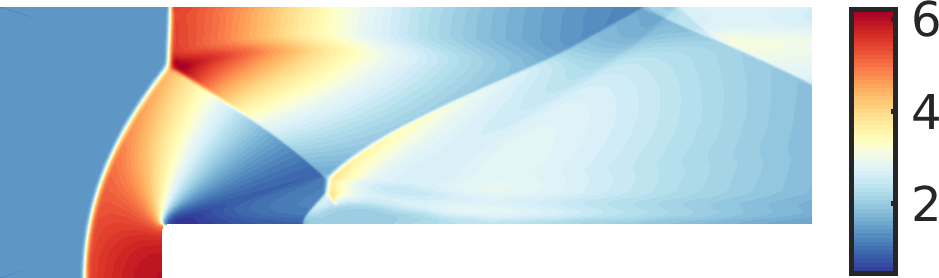} &
        \includegraphics[width=0.46\linewidth]{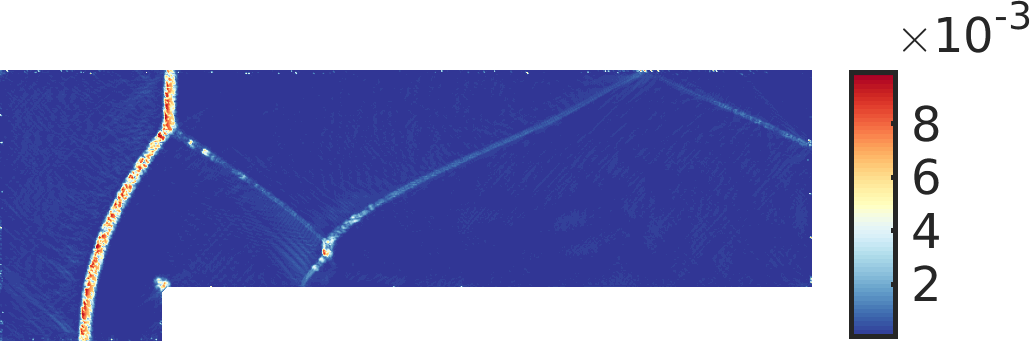} \\
         \raisebox{2\normalbaselineskip}[0pt][0pt]{\rotatebox[origin=c]{90}{$\mathbf{p=3}$}} &
        \includegraphics[width=0.42\linewidth]{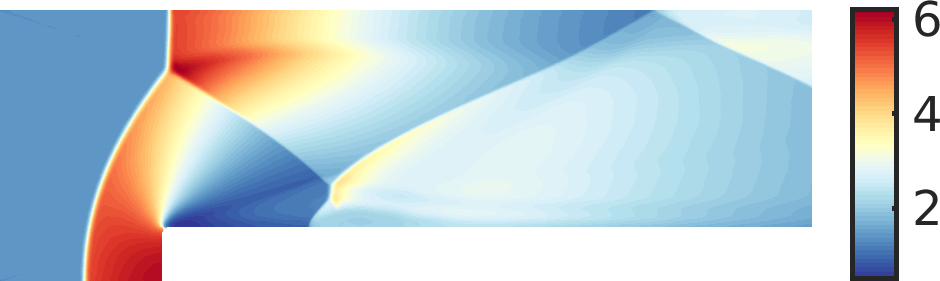} &
        \includegraphics[width=0.46\linewidth]{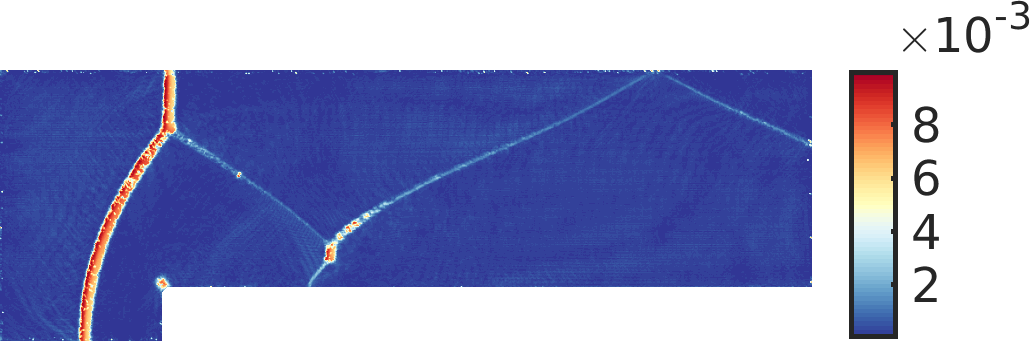}
    \end{tabular}
    \caption{The results are made using $N=117432$ nodes. The numerical solution (left column) is given in terms of density $\rho$ at $t=4$, 
    for two choices of a monomial basis degree $p$ to construct the stencil based approximation. The plots in the right column are residual viscosity 
    coefficients for the two choices of $p$.}
    \label{fig:euler:forwardStep}
\end{figure}

The numerical results are collected in Figure \ref{fig:euler:forwardStep}. The numerical solutions for $p=1$ and $p=3$ 
are similar to the solutions displayed in Figure 9 of \cite{NazarovLarcher17}. A prevailing difference 
is that the shock reflection close to the inward facing corner has a different location. This is ascribed to the fact that 
in this paper we do not smooth out the inward corner as the authors in \cite{NazarovLarcher17} do, 
but keep it sharp. In our case, the contact discontinuity is 
captured for both choices of $p$. When $p=1$ we notice that the solution around the reflection closest to the outflow 
of the channel, has a large error by means of the phase, compared to when using $p=3$. This shows a benefit of using a
high-order discretization, even if the underlying solution does not have a sufficient regularity to support the high-order 
convergence.






\subsection{Explosion in a domain with cylinders}
This problem was initially introduced in \cite{NazarovLarcher17}. 
The computational domain is the disc with radius $2$. In addition, the domain has 
eight inner circular boundaries with radius $0.3$, which are placed at distance $1.4$ from the origin. 
The distance between the neighboring inner circular boundaries are equal. The initial condition is a discontinuous function, where 
a compressed gas with $\rho = 1$, $p=1$ is put inside the disc with radius $\sqrt{0.4}$, located at the origin. 
Outside of that disc we have $\rho = 0.125$ and $p=1$. 
The velocity is $\bm v = (0,0)$ throughout the domain. Slip boundary conditions are used for the exterior boundary and all of
the interior boundaries. The simulation is run until $t=4.25$.

\begin{figure}[h!]
    \centering
    \begin{tabular}{cc}
        \multicolumn{2}{c}{\textbf{Euler, Explosion in a domain with cylinders}\vspace{0.05cm}} \\
        \multicolumn{2}{c}{\textbf{Numerical solution}\vspace{0.2cm}} \\
        \includegraphics[width=0.33\linewidth]{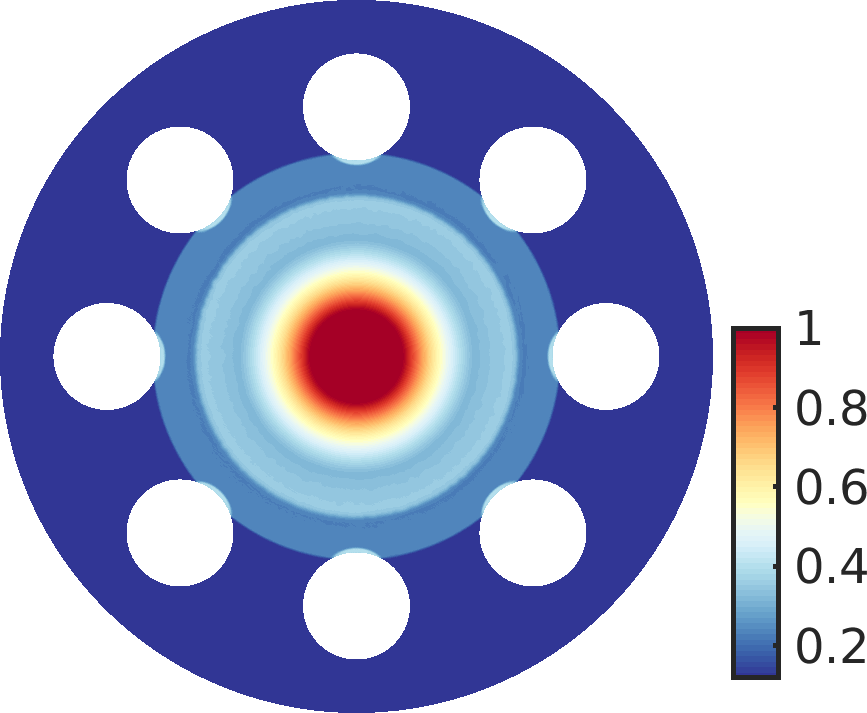} & 
        \includegraphics[width=0.33\linewidth]{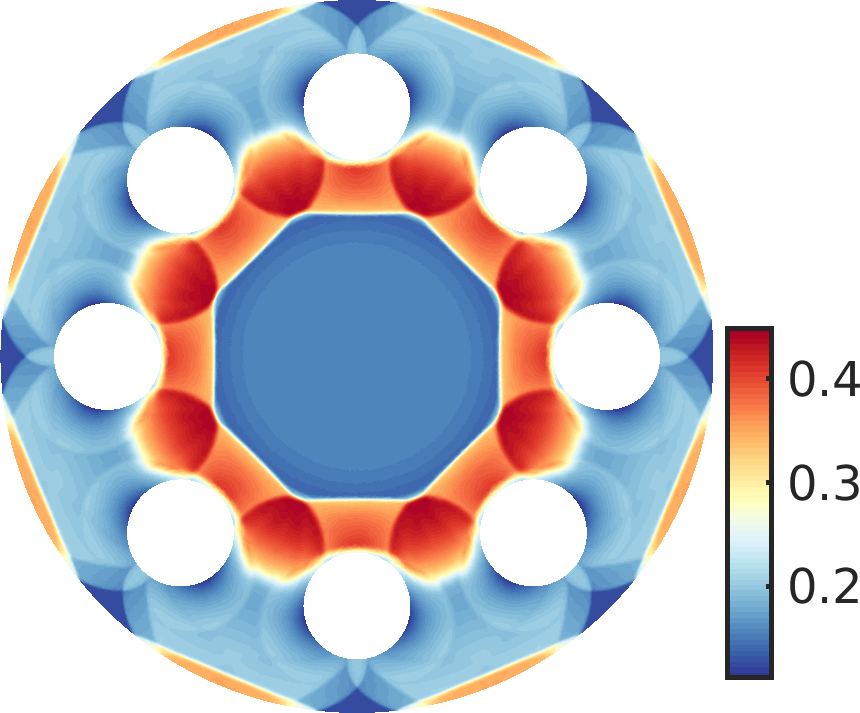} \\
        \hspace{-0.8cm}$\mathbf{t=0.3}$ \vspace{0.2cm} & \hspace{-0.8cm}$\mathbf{t=0.9}$ \\
        \includegraphics[width=0.33\linewidth]{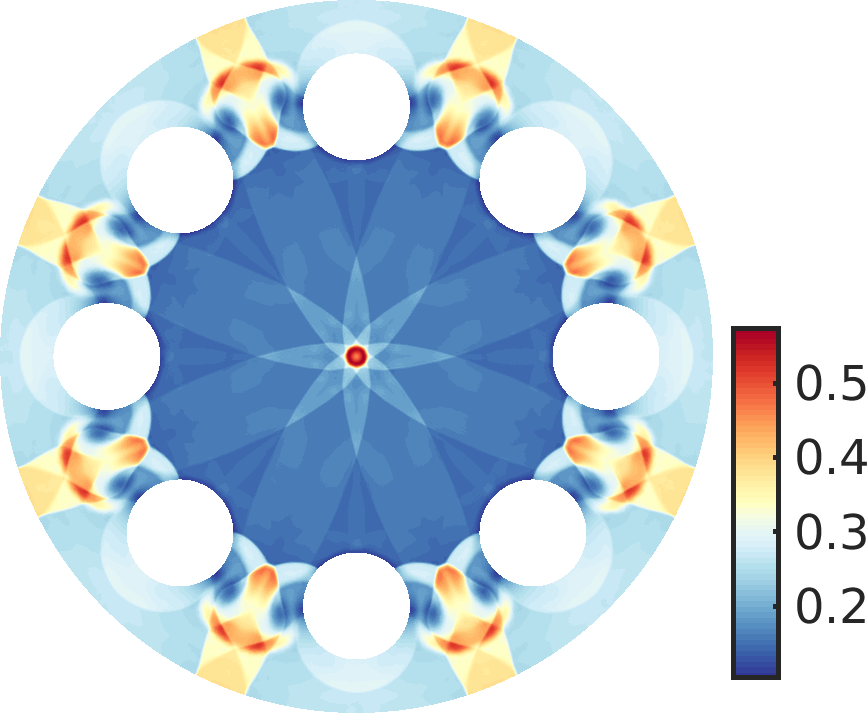} &
        \includegraphics[width=0.33\linewidth]{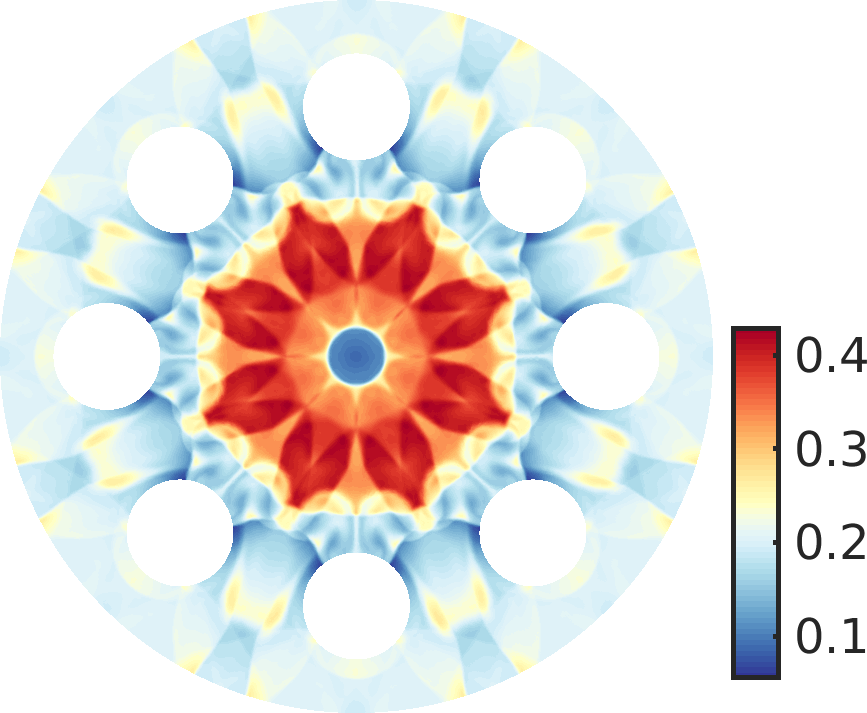} \\
        \hspace{-0.8cm}$\mathbf{t=1.8}$ \vspace{0.2cm} & \hspace{-0.8cm}$\mathbf{t=2.71}$ \\
        \includegraphics[width=0.33\linewidth]{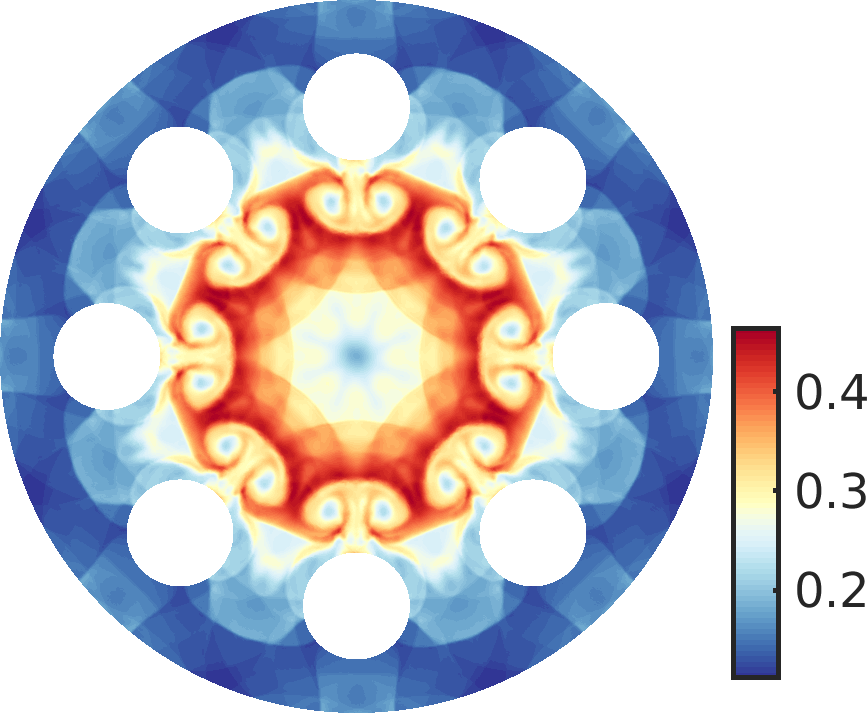} & 
        \includegraphics[width=0.33\linewidth]{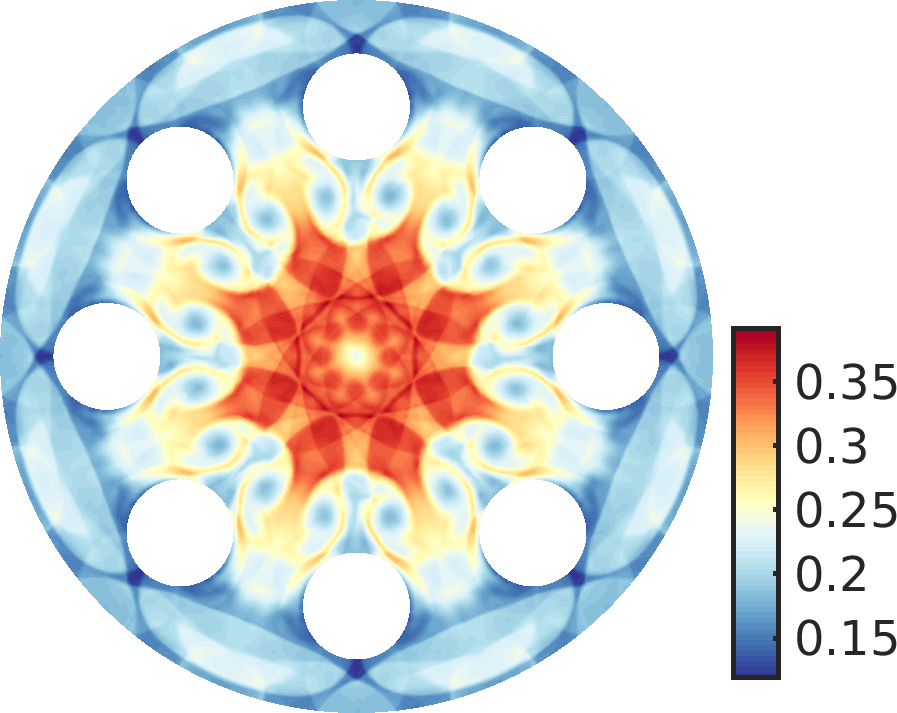} \\
        \hspace{-0.8cm}$\mathbf{t=3.61}$ \vspace{0.2cm} & \hspace{-0.8cm}$\mathbf{t=4.25}$ \\
    \end{tabular}
    \caption{Numerical solution of the benchmark Explosion in a domain with cylinders, at different points in time $t$. 
    The solution is obtained using $N=337356$ nodes. The oversampling parameter is $q=5$. The polynomial degree used to construct 
    the stencil-based approximation is $p=3$.}
    \label{fig:euler:explosion}
\end{figure}

\begin{figure}[h!]
    \centering
    \begin{tabular}{cc}
        \multicolumn{2}{c}{\textbf{Euler, Explosion in a domain with cylinders}\vspace{0.05cm}} \\
        \multicolumn{2}{c}{\textbf{Residual viscosity coefficient}\vspace{0.2cm}} \\
        \includegraphics[width=0.33\linewidth]{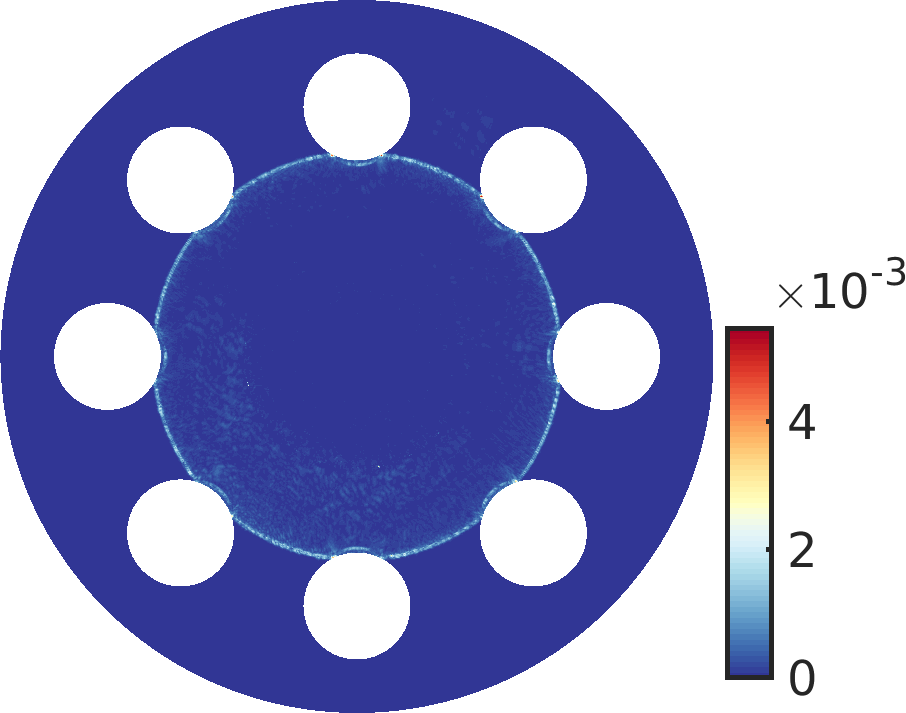} & 
        \includegraphics[width=0.33\linewidth]{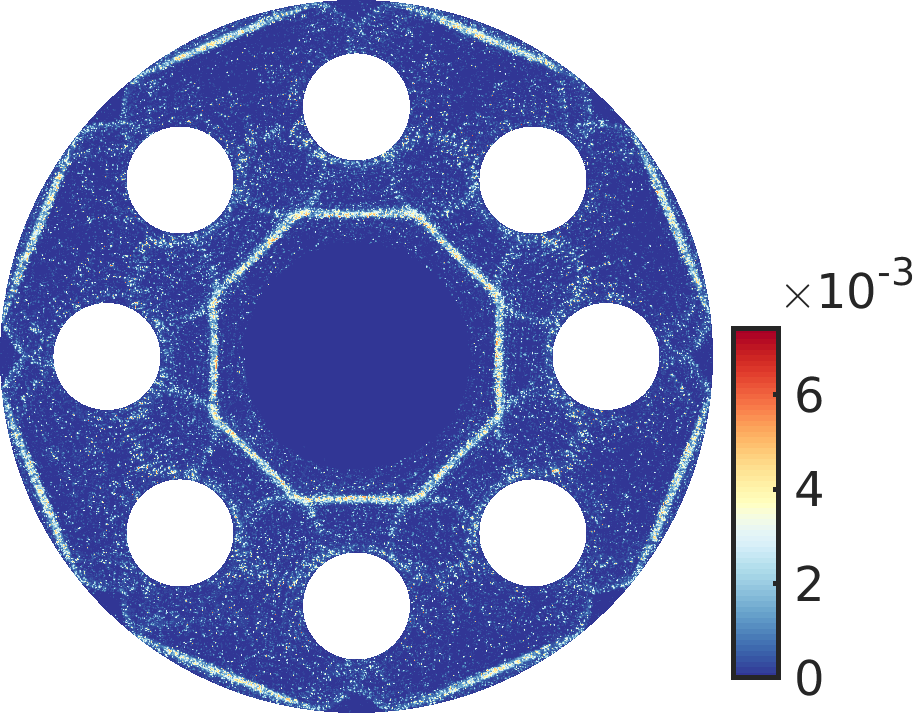} \\
        \hspace{-0.8cm}$\mathbf{t=0.3}$ \vspace{0.2cm} & \hspace{-0.8cm}$\mathbf{t=0.9}$ \\
        \includegraphics[width=0.33\linewidth]{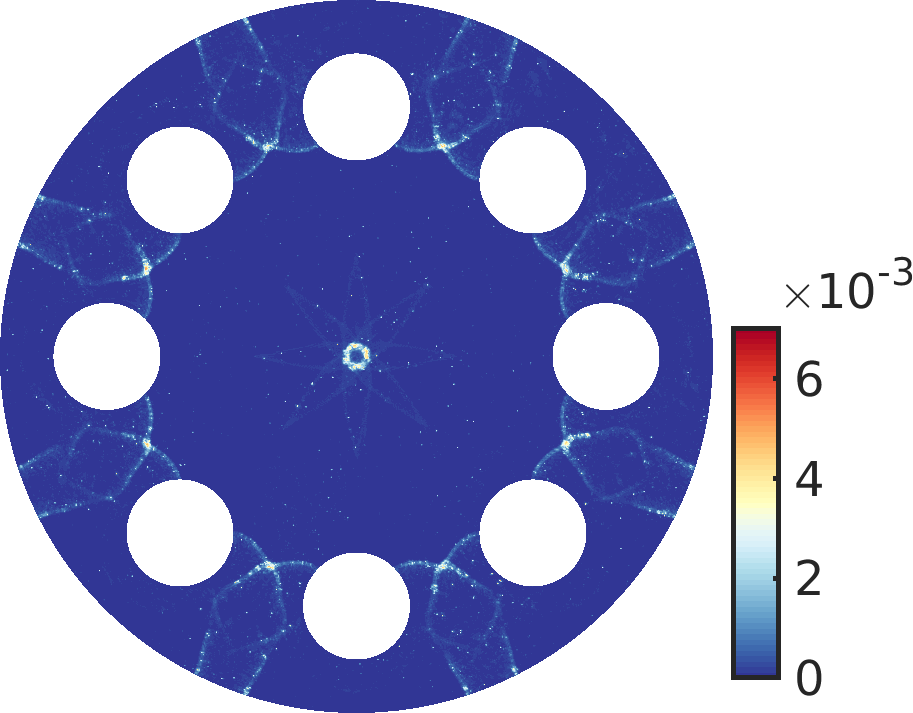} &
        \includegraphics[width=0.33\linewidth]{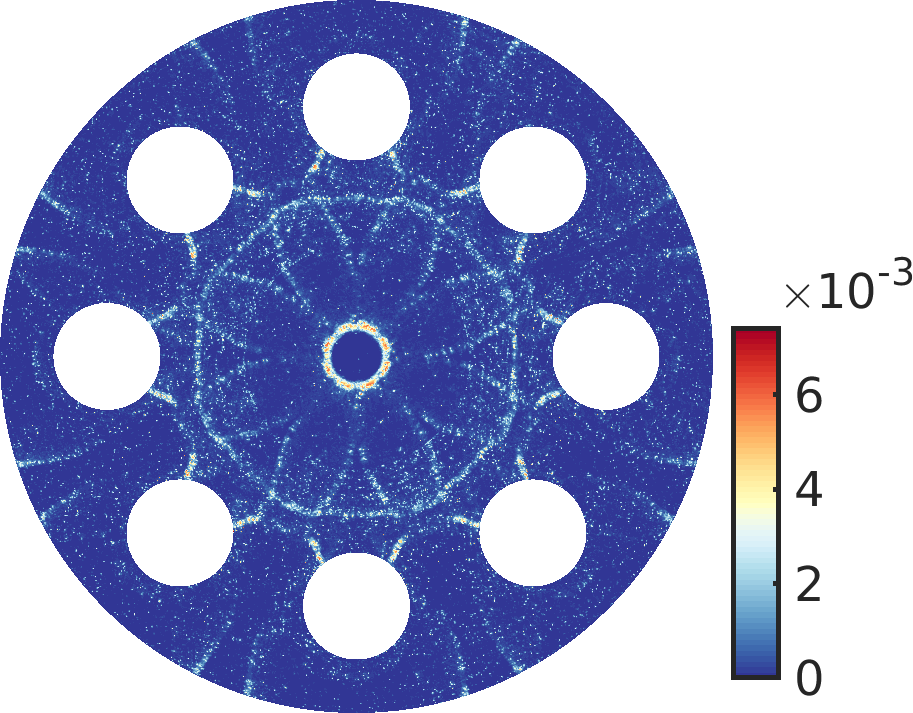} \\
        \hspace{-0.8cm}$\mathbf{t=1.8}$ \vspace{0.2cm} & \hspace{-0.8cm}$\mathbf{t=2.71}$ \\
        \includegraphics[width=0.33\linewidth]{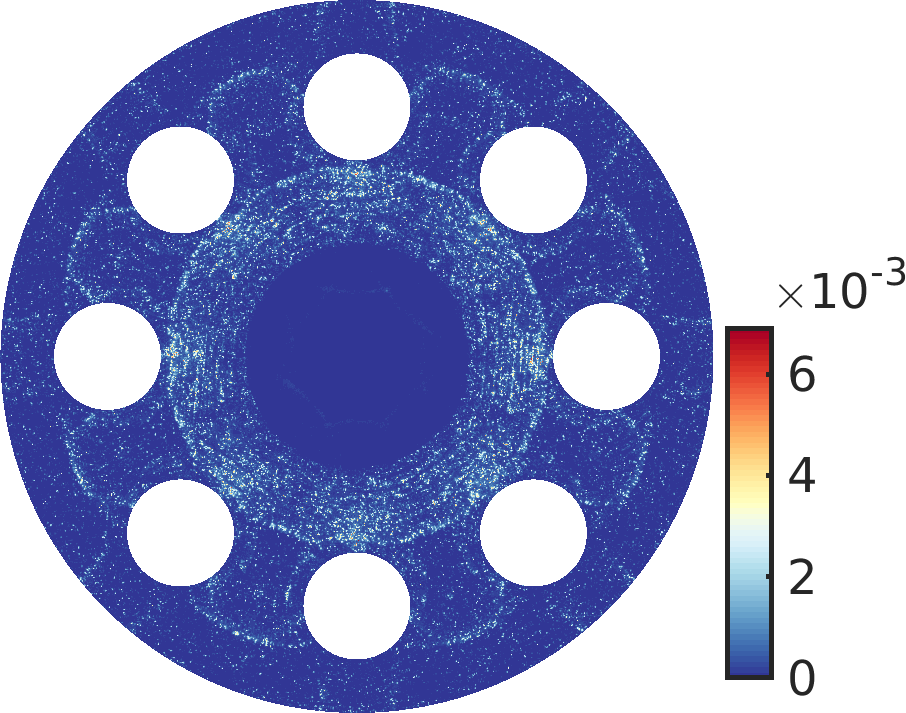} & 
        \includegraphics[width=0.33\linewidth]{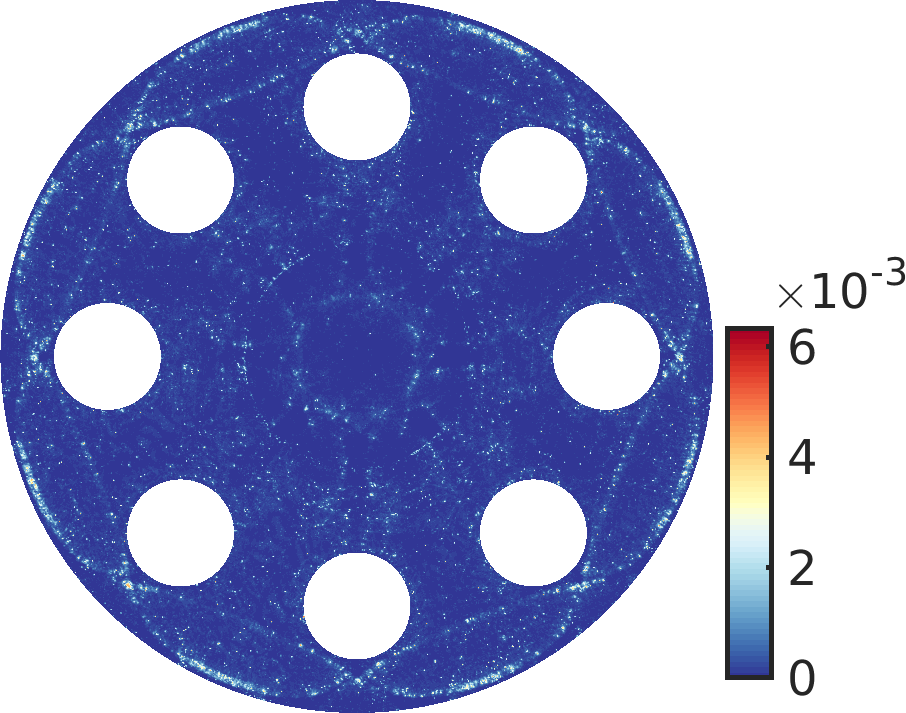} \\
        \hspace{-0.8cm}$\mathbf{t=3.61}$ \vspace{0.2cm} & \hspace{-0.8cm}$\mathbf{t=4.25}$ \\
    \end{tabular}
    \caption{Spatial distribution of the residual viscosity coefficient used to compute the solution of the benchmark Explosion in a domain with cylinders. 
    The coefficient is displayed at different points in time $t$. 
    The corresponding solution is computed obtained using $N=337356$ nodes. 
    The oversampling parameter is $q=5$. The polynomial degree used to construct 
    the stencil-based approximation is $p=3$.}
    \label{fig:euler:explosion_rvCoefficients}
\end{figure}

We use $N=337356$ nodes obtained using Gmsh \cite{Gmsh}. 
Other parameters are: $p=3$, $C_{\text{RV}} = 3$ and $\text{CFL} = 0.3$. 
The results are given in Figure \ref{fig:euler:explosion}, where we observe that RV performs well, even in cases where 
the numerical solution exhibits small details. The spatial distribution of the residual viscosity coefficient 
is displayed in Figure \ref{fig:euler:explosion_rvCoefficients}, where we see that the action of viscosity is localized.

\section{Final discussion}
\label{sec:finaldiscussion}
In this paper we covered several aspects of discretizing conservation laws using the RBF-FD method. 

We observed that oversampling by itself does not improve the time stability of the RBF-FD method. 
A theoretical argumentation for this behavior and further analysis is given in \cite{tominec2021_stability}.

Another observation is that the collocation RBF-FD method is nearly as accurate as the oversampled RBF-FD method 
when using Dirichlet boundary conditions, which is in line with the observations made in \cite{ToLaHe21}. 
\igor{As investigated in Section \ref{sec:experiments:linearadvection}, the residual viscosity stabilization framework 
is applicable to both, the collocation setting and the oversampled setting of the RBF-FD method.} 

The RV constant $C_{\text{RV}}$ which is defined by the user is not a sensitive parameter. We did, however, notice that 
$C_{\text{RV}}$ had to be chosen larger than $1$, when the magnitude of the numerical solution was large. This was not observed 
when RV was combined with the finite element methods, for example in \cite{NazarovLarcher17}.
We did not fully explore the reasons behind that. A speculative explanation is that we used 
a different residual definition compared to the residual definition used in \cite{NazarovLarcher17}.

We found that a symmetric hyperviscosity operator \eqref{eq:discretization:hypervi} is an effective stabilization in time. 
The residual as defined in the present work turned out to be an excellent indicator of discontinuities. 
The residual viscosity stabilization (RV) is consistent when the solution is smooth. 

Finally, all experiments confirmed that a combination of the RBF-FD method and the RV stabilization provides 
a robust framework for discretizing nonlinear conservation laws in scalar and system settings.

\section*{Acknowledgments}
We thank (in alphabetical order) Lukas Lundgren, Tuan Anh Dao and Vidar Stiernström from Uppsala University for fruitful discussions about 
time-dependent conservation laws.

\appendix
\section{(Appendix) A high-order time derivative approximation}
\label{sec:appendix:timederivatives}
\begin{verbatim}
    function w = timederivatives (t)
        % Input: a vector t, where t(i) is time at which the solution is available.
        % Output: a vector w, where each w(i) is used to multiply u|_{t(i)} in order
        % ... to get a derivative at t(end).
        % Usage: d/dt u(t_end) = w(end)*u(end) + w(end-1)*u(end-1) + ... + w(1)*u(1),
        % ... where t_end is the time at which the last solution point is available.

        scale = 1/max(abs(t));
        t = t*scale;
        t_eval = t(end); % The derivative should be evaluated at t(end).
    
        % Construct the polynomial basis, and differentiate it in a point t_eval.
        A = zeros(size(t,1), size(t,1));
        b_t = zeros(1, size(t,1));
        for k=1:length(t)
            A(:,k) = t.^(k-1);
            b_t(k) = (k-1)*t_eval.^(k-2);
        end
        w = scale*(b_t*inv(A));
    
    end
\end{verbatim}

\bibliographystyle{spmpsci}
\bibliography{refs}

\end{document}